# On the behavior of sequences of solutions to U(1) Seiberg-Witten systems in dimension 4.

Clifford Henry Taubes[†]

Department of Mathematics
Harvard University
Cambridge, MA 02318

(chtaubes@math.harvard.edu)

ABSTRACT

This paper studies the behavior of sequences of solutions to Seiberg-Witten like equations for a pair consisting of a Hermitian connection on a line bundle over a 4-dimensional manifold and a section of the self-dual spinor bundle of a complex Clifford module on the manifold. Examples include the cases where the Clifford module is a direct sum of $\mathbb{C}^2$ bundles associated to $\text{Spin}^{\mathbb{C}}$ structures; and the case of the $SU(2)$ Vafa-Witten equations with an Abelian ansatz.

[†]Supported in part by the National Science Foundation

# 1. Introduction

The upcoming Theorem 1.1 (the principle theorem in this paper) describes the behavior of sequences of solutions to certain generalizations of the Seiberg Witten equations on a compact, oriented 4-dimensional Riemannian manifold. These are equations for a pair consisting of a Hermitian connection on a complex line bundle and a section of the self-dual spinor bundle of a given complex Clifford module over the manifold. The original Seiberg-Witten equations (see [W]) concern self-dual spinor bundles with rank 2, whereas Theorem 1.1 considers bundles with rank greater than 2. Haydys and Walpuski [HW] stated and proved a version of Theorem 1.1 for an analogous set of equations on a 3-dimensional manifold.

## a) Self-dual Clifford modules

This first subsection sets the stage and some of the notation for Theorem 1.1. To start, let X denote a compact, oriented Riemannian, 4-dimensional manifold. Let $\mathbb{S}^+$ and $\mathbb{S}^-$ denote a given pair of complex, Hermitian vector bundles of the same dimension over X; both with metric compatible covariant derivatives. (Both of these covariant derivatives are denoted by $\nabla$.) Let $\mathrm{Iso}(\mathbb{S}^+, \mathbb{S}^-)$ denote the fiber bundle over X whose fiber at any given $p \in X$ are the isometries from $\mathbb{S}^+|_p$ to $\mathbb{S}^-|_p$. Suppose that $\mathfrak{c}$ is a homomorphism from T*X to $\mathrm{Iso}(\mathbb{S}^+, \mathbb{S}^-)$. The data consisting of $\mathbb{S}^+$, $\mathbb{S}^-$, $\mathfrak{c}$, and the covariant derivatives on T*X, $\mathbb{S}^+$ and $\mathbb{S}^-$ is said to define a Clifford module over X when the following conditions are met:

- *The homomorphism $\mathfrak{c}$ is covariantly constant with respect to the Levi-Civita covariant derivative on T*X and the given covariant derivatives on $\mathbb{S}^+$ and $\mathbb{S}^-$.*
- *If p is any given point in X, and if v and w are any given elements in $T^*X|_p$, then*

$$\mathfrak{c}^\dagger(v)\mathfrak{c}(w) + \mathfrak{c}^\dagger(w)\mathfrak{c}(v) = 2\langle v, w \rangle \mathbb{I}$$

(1.1)

Here $\langle\, ,\,\rangle$ denotes the inner product on T*X and $\mathbb{I}$ denotes the identity automorphism. (In what follows, $\langle\, ,\,\rangle$ is used to denote the inner product on any given Hermitian or Riemannian vector bundle.)

A homomorphism cl: $\wedge^2 T^*X \to \mathrm{End}(\mathbb{S}^+)$ is defined by its action on decomposable covectors as follows: Supposing that $p \in X$ and $v, w \in T^*X|_p$, then

$$\mathrm{cl}(v \wedge w) = \mathfrak{c}^\dagger(v)\mathfrak{c}(w) - \mathfrak{c}^\dagger(w)\mathfrak{c}(v) \,.$$

(1.2)

The bundle $\wedge^2 T^*X$, which is the domain of $\mathrm{cl}(\cdot)$, splits orthogonally as $\Lambda^+ \otimes \Lambda^-$ with $\Lambda^+$ and $\Lambda^-$ being the respective +1 and -1 eigenbundles for the metric's Hodge star operator. The given Clifford module is said (in this paper) to be *self-dual* when $\mathrm{cl}(\cdot)$ annihilates $\Lambda^-$.

It is assumed henceforth that the Clifford module in question is, in fact, self dual. Examples are given momentarily. Supposing this self-duality, then the adjoint of cl (to be denoted by $cl^\dagger$) is a homomorphism from $End(\mathbb{S}^+)$ to $\Lambda^+ \otimes_\mathbb{R} \mathbb{C}$. Note that this adjoint homomorphism $cl^\dagger$ sends Hermitian endomorphisms of $\mathbb{S}^+$ to $i\mathbb{R}$-valued, self dual 2-forms. As a consequence, the adjoint $cl^\dagger$ defines a quadratic map from $\mathbb{S}^+$ to $i\Lambda^+$ by the rule $x \to cl^\dagger(x \otimes x^\dagger)$. (Given that the Clifford module is self dual, the respective bundles $\mathbb{S}^+$ and $\mathbb{S}^-$ will be called the self dual and anti-self dual spinor bundles.)

Let $E \to X$ denote a given complex, Hermitian line bundle; and let A denote a Hermitian connection on E. This connection with $\nabla$ induces a covariant derivative for sections of $\mathbb{S}^+ \otimes E$ and $\mathbb{S}^- \otimes E$ to be denoted by $\nabla_A$. Since $\mathfrak{c}$ defines $\mathbb{C}$-linear homomorphisms from $\mathbb{S}^+$ to $\mathbb{S}^-$, it canonically defines a self-dual, Clifford module structure with the bundles $\mathbb{S}^+ \otimes E$ and $\mathbb{S}^- \otimes \mathbb{E}$ and their versions of $\nabla_A$. This data (the Clifford module data $\{\mathbb{S}^+, \mathbb{S}^-, \mathfrak{c}, \nabla\}$ with the bundle E and the connection A) defines a Dirac operator, $D_A$, that maps sections of $\mathbb{S}^+ \otimes E$ to sections of $\mathbb{S}^- \otimes E$ according to the rule

$$\hat{a} \to D_A\hat{a} = \mathfrak{c}(\nabla_A)\hat{a} \ .$$

(1.3)

By way of more notation, $F_A$ is used to denote the curvature 2-form of A (an $i\mathbb{R}$-valued 2-form), and $F_A^+$ is used to denote the self-dual part of $F_A$.

**b) U(1) Seiberg-Witten systems**

Fix an $i\mathbb{R}$ valued self-dual 2 form on X to be denoted by $\omega$ (it can be zero). This 2-form plus the data giving the self-dual Clifford module (this is $(\mathbb{S}^+, \mathbb{S}^-, \nabla, \mathfrak{c})$) and the complex Hermitian line bundle E defines a set of differential equations, a *U(1) Seiberg-Witten system*. This is the system of differential equations for a pair $(A, \hat{a})$ with A being a Hermitian connection on E and with $\hat{a}$ being a section of $\mathbb{S}^+ \otimes E$ that demand

$$F_A^+ = -\tfrac{1}{4} cl^\dagger(\hat{a} \otimes \hat{a}^\dagger) + \omega \quad and \quad D_A\hat{a} = 0 \ .$$

(1.4)

These equations are invariant under the action of the group $C^\infty(X; S^1)$ via the following action: A map $g: X \to S^1$ acts on a connection, A, on E by the rule $A \to g(A) = A - g^{-1}dg$; and it acts on a section $\hat{a}$ of $\mathbb{S}^+ \otimes E$ by the rule $\hat{a} \to g\hat{a}$. The map g then acts on pairs of connection on E and section of $\mathbb{S}^+ \otimes E$ by the product of these separate actions.

The original Seiberg-Witten equations [W] has $\mathbb{S}^+$ being the self-dual spinor bundle (a $\mathbb{C}^2$ vector bundle) that is associated to a $Spin^\mathbb{C}(4)$ lift of the principal $SO(4)$ bundle of oriented, orthonormal frames for TX with a connection that is induced by the Riemannian metric connection on TX and a chosen Hermitian connection on $\wedge^2 \mathbb{S}^+$. (The self-dual spinor bundle from a given Spin$^{\mathbb{C}}$(4) lift of the frame bundle is here by $\mathbb{S}_{\Diamond}^+$; and the corresponding version of $\mathbb{S}^-$ is denoted by $\mathbb{S}_{\Diamond}^-$.) As explained in [W] (see also [Mo], [KM]), this instance of (1.4) has the remarkable property that the space of $C^{\infty}(X; S^1)$ orbits of its solutions is sequentially compact.

  This sequential compactness property (modulo the action of $C^{\infty}(X; S^1)$ is not expected for other instances of (1.4). Here are some to keep in mind: Let n denote for the moment a positive integer. For each $k \in \{1, 2, \ldots, n\}$, let $\mathbb{S}_k^+$ denote a version of the bundle $\mathbb{S}_{\Diamond}^+$ that was just described. Define $\mathbb{S}_k^-$ accordingly. (Thus, $\mathbb{S}_k^+$ is the self-dual spinor bundle associated to a Spin$^{\mathbb{C}}$(4) lift of the principle bundle of oriented, orthonormal frames in TX; it is a $\mathbb{C}^2$ bundle.) Take $\mathbb{S}^+$ to be $\oplus_{k=1,\ldots,n} \mathbb{S}_k^+$ and take $\mathbb{S}^-$ to be $\oplus_{k=1,2,\ldots,n} \mathbb{S}_k^-$. There are n ≥ 2 examples (see [BW]) where sequential compactness fails.

  A different sort of example takes $\mathbb{S}^+$ to be the bundle $(\Lambda^+ \otimes_{\mathbb{R}} \mathbb{C}) \oplus \underline{\mathbb{C}}$ with $\underline{\mathbb{C}}$ denoting the product complex line bundle over X. The bundle $\mathbb{S}^-$ in this case is $T^*X \otimes_{\mathbb{R}} \mathbb{C}$. The respective connections on these bundles are the Riemannian metric connections. The Clifford module structure in this case comes by writing $\mathbb{S}^+$ as $\mathbb{S}_{\Diamond}^+ \otimes (\mathbb{S}_{\Diamond}^+)^*$, and writing $\mathbb{S}^-$ as $\mathbb{S}_{\Diamond}^- \otimes (\mathbb{S}_{\Diamond}^+)^*$ with $(\mathbb{S}_{\Diamond}^+)^*$ denoting the dual bundle to $\mathbb{S}_{\Diamond}^+$. The equations in (1.4) (with $\omega = 0$) for this example are equivalent to the Vafa-Witten equations [VW] with an extra ansatz. (The original Vafa-Witten equations as written in [VW] are an instance of (1.4) with $\mathbb{S}^+$ and $\mathbb{S}^-$ as just described, but with E being the $\mathbb{C}^3$ bundle that is obtained by complexifying the associated lie algebra bundle of a principle SU(2) or SO(3) bundle over X. The connection A in this case is supposed to be a unitary connection on E. See [VW] and the more recent references, Ben Mares's Ph.D. thesis [Ma] and [GGP].) A second example along the same lines is an Abelian version of one of the equations that Kapustin and Witten wrote down in [KW]. (See also [T1] and [GU]). This system has $\mathbb{S}^+$ being the $T^*X \otimes_{\mathbb{R}} \mathbb{C}$ and $\mathbb{S}^- = (\Lambda^- \otimes_{\mathbb{R}} \mathbb{C}) \oplus \underline{\mathbb{C}}$. The covariant derivatives on these spinor bundles are again those from the Riemannian structure; and the Clifford module structure is defined by writing spinor bundle $T^*X \otimes_{\mathbb{R}} \mathbb{C}$ as $\mathbb{S}_{\Diamond}^+ \otimes (\mathbb{S}_{\Diamond}^-)^*$ and $(\Lambda^- \otimes_{\mathbb{R}} \mathbb{C}) \oplus \underline{\mathbb{C}}$ as $\mathbb{S}_{\Diamond}^- \otimes (\mathbb{S}_{\Diamond}^-)^*$.

  The preceding examples are special cases of a general form whereby $\mathbb{S}^+$ is $\mathbb{S}_{\Diamond}^+ \otimes \mathcal{E}$ and $\mathbb{S}^-$ is $\mathbb{S}_{\Diamond}^- \otimes \mathcal{E}$ with $\mathcal{E}$ being a fixed complex vector bundle (with Hermitian metric and connection) on X.

  Arguments that differ little from those in [Ma] prove that any sequence of solutions to (1.4) with an a priori bound on the $L^2$ norm of the $\hat{a}$ component of each pair from the sequence can be altered term-wise by automorphisms of E (in the manner described above) so that the resulting sequence has a subsequence that converges in the $C^\infty$ topology to a solution to (1.4). The upcoming Proposition 2.2 makes a formal statement to this effect.

The preceding observation is cause to rewrite the equations in (1.4) so as to highlight the role that is played by the $L^2$ norm of the section $\hat{a}$. This is done by introducing a parameter $r \in [0, \infty)$. Given $r$ consider now the equations written directly below for a pair $(A, a)$ of a connection A on E and a section $a$ of $\mathbb{S}^+$:

- $F_A^+ = -\frac{1}{4} r^2 \text{cl}^\dagger(a \otimes a^\dagger) + \omega$.
- $D_A a = 0$.
- $\int_X |a|^2 = 1$.

(1.5)

If $(A, a)$ obeys these equations, then $(A, \hat{a} = ra)$ obeys the equations in (1.4). Conversely, if $(A, \hat{a})$ obeys (1.4) and $\hat{a} \neq 0$, then $(A, a = r^{-1}\hat{a})$ obeys (1.5) with $r$ being the $L^2$ norm of the section $\hat{a}$ on X. As just noted, any sequence $\{r_n, (A_n, a_n)\}_{n=1,2,\ldots}$ obeying (1.5) with $\{r_n\}_{n \in \{1,2,\ldots\}}$ containing bounded sequence has a $C^\infty$ convergent subsequences up to the action of the automorphisms of E.

The following theorem (it is the centerpiece of this paper) says something about what happens to sequence $\{r_n, (A_n, a_n)\}_{n=1,2,\ldots}$ that obey (1.5) when $\{r_n\}_{n=1,2,\ldots}$ lacks bounded subsequences. By way of notation, the theorem uses $\mathcal{R}^\nabla$ to denote curvature of the connection that defines the given covariant derivative $\nabla$ on $\mathbb{S}^+$. (This is a $\wedge^2 T^*X$ valued, anti-Hermitian endomorphism of $\mathbb{S}^+$.)

**Theorem 1.1**: *Fix a self-dual Clifford module over X, a complex, Hermitian line bundle $E \to X$, and an $i\mathbb{R}$ valued, self-dual 2-form on X. Let $\{(r_k, (A_k, a_k))\}_{k=1,2,\ldots}$ denote a sequence of solutions to the corresponding version of (1.5) with the sequence $\{r_k\}_{k=1,2,\ldots}$ increasing and unbounded. Given this sequence, there is the following data:*
- *A closed, nowhere dense subset of X to be denoted by Z.*
- *A smooth, Hermitian connection on $E|_{X-Z}$ (to be denoted by $\hat{A}$) and a smooth section of $(\mathbb{S}^+ \otimes E)|_{X-Z}$ (to be denoted by $v$) that are described by the three items below.*
  a) *The norm of $v$ extends to the whole of X as a Hölder continuous function that is zero on the set Z.*
  b) *The norm of $\nabla_{\hat{A}} v$ extends to the whole of X as a square integrable function.*
  c) $D_{\hat{A}} v = 0$ *and* $\text{cl}^\dagger(v \otimes v^\dagger) = 0$ *and* $v^\dagger \nabla_{\hat{A}} v = (\nabla_{\hat{A}} v)^\dagger v$.
  d) $F_{\hat{A}} = -\frac{1}{|v|^2}((\nabla_{\hat{A}} v)^\dagger \wedge \nabla_{\hat{A}} v + v^\dagger \mathcal{R}^\nabla v)$.

- *A subsequence $\Theta \subset \{1, 2, \ldots\}$ and a corresponding sequence of Hermitian automorphisms of $E|_{X-Z}$ indexed by $\Theta$. This sequence is denoted by $\{g_n\}_{n \in \Theta}$.*

*This data is such that $\{g_n(A_n)\}_{n \in \Theta}$ converges to $\hat{A}$ in the $L^2_1$ weak topology on compact subsets of $X-Z$; and $\{g_n a_n\}_{n \in \Theta}$ converges to $\nu$ in the $L^2_2$ topology on compact subsets of $X-Z$. In addition, the sequence of norms $\{|a_n|\}_{n \in \Theta}$ converges to $|\nu|$ in the $L^2_1$ and $C^0$ topologies on the whole of $X$, and the sequence $\{|\nabla_{A_n} a_n|\}_{n \in \Theta}$ converges to $|\nabla_{\hat{A}} \nu|$ in the $L^2$ topology on the whole of $X$.*

The proof of this theorem is given in the upcoming Section 10 of this paper. The intervening sections establish various auxilliary lemmas and propositions that are used in the proof. Section 1e has a brief sketch of the arguments for the theorem.

      Note that the 3-dimensional analog of Theorem 1.1 was proved by Haydys and Walpuski [HW]; and their sequential compactness assertion is essentially the same as Theorem 1.1. Many of the construction that follow have antecedents in [HW].

      It is likely that the set Z from Theorem 1.1 has Hausdorff dimension at most 2 in all cases (and, perhaps, it is always rectifiable with finite 2-dimensional Hausdorff measure.) In the special case when $\hat{A}$ is flat with $\mathbb{Z}/2$ holonomy, then the Hausdorff dimension of Z is known to be at most 2. (The theorems in [T2] apply in this case and these theorems give 2 as the upper bound for Z's Hausdorff dimension.) As explained in Section 1d, the theorems in [T2] apply in certain cases when $\mathbb{S}^+$ has dimension 4; the Abelian Vafa-Witten equations being one instance. (These equations are the version of (1.5) with $\mathbb{S}^+ = \Lambda^+ \otimes \mathbb{C} \oplus \underline{\mathbb{C}}$). Examples are known where the dimension of Z is, in fact, 2. See e.g. [BW]. Explicit examples can also be obtained by taking X to be the product of two Riemann surfaces and taking an input sequence for Theorem 1.1 that depends only on one of them.

c) **A very slight generalization**

      The same arguments (but for some very straightforward cosmetic points and some notation) that are used in what follows to prove Theorem 1.1 can be used to prove an analogous theorem for a slight generalization of Theorem 1.1. To state this generalization, first fix a positive integer to be denoted by n. For each $i \in \{1, 2, \ldots, n\}$, let $m_i$ denote a non-zero integer; and let $(\mathbb{S}_i^+, \mathbb{S}_i^-, \nabla, \mathfrak{c})$ denote the data for a complex, self-dual Clifford module on X. Let E denote, as before, a complex, Hermitian line bundle over X. Keep in mind that if A is a connection on E, then A defines a connection on any power of E (supposing that m is an integer, then the induced connection on the m'th power of E is denoted by $A^{(m)}$.) Therefore, if $i \in \{1, \ldots, n\}$, then the connection A and the chosen covariant derivative on $C^\infty(X; \mathbb{S}_i^+)$ defines a Dirac operator on $C^\infty(\mathbb{S}_i^+ \otimes E^{m_i})$ using the formula in (1.3).

The generalization of (1.5) are equations for a data set $(A, a = (a_1, \ldots, a_n))$ with A being a Hermitian connection on E; and with each $i \in \{1, 2, \ldots, n\}$ version of $a_i$ being a section of the corresponding $\mathbb{S}_i^+ \otimes E^{m_i}$. The equations that generalize (1.5) are written below. (The first of these equations uses $\text{I}a$ as shorthand for $(m_1 a_1, \ldots, m_n a_n)$.)

- $F_A^+ = -\frac{1}{4} r^2 \text{cl}^\dagger(\text{I}a \otimes a^\dagger) + \omega$ ,
- $D_{A^{(m_i)}} a_i = 0$ ,
- $\int_X |a|^2 = 1$ .

(1.6)

Here, $\text{I}a = (m_1 a_1, \ldots, m_n a_n)$. Note that if all of the $m_i$'s are equal, then the equations in (1.6) are an instance of (1.5) up to a redefinition of A. The analog of Theorem 1.1 for the equations in (1.6) is identical in almost all respect to Theorem 1.1 with the only difference being in Items c) and d) of the second bullet. When each $k \in \{1, 2, \ldots\}$ version of $(A_k, a_k)$ obeys (1.6), then the second bullet of the Theorem 1.1 analog asserts the existence of:

*A smooth, Hermitian connection on $E|_{X-Z}$ (to be denoted by $\hat{A}$) and a smooth section of $\mathbb{S}^+ \otimes E$ (to be denoted by $v$) that are described by the three items below.*

a) *The norm of $v$ extends to the whole of X as a Hölder continuous function that is zero on the set Z.*
b) *The norm of $\nabla_{\hat{A}} v$ extends to the whole of X as a square integrable function.*
c) $D_{\hat{A}} v = 0$ *and* $\text{cl}^\dagger(\text{I}v \otimes v^\dagger) = 0$ *and* $v^\dagger \text{I} \nabla_{\hat{A}} v = (\nabla_{\hat{A}} v)^\dagger \text{I} v$.
d) $F_{\hat{A}} = -\frac{1}{|v|^2}((\nabla_{\hat{A}} v)^\dagger \wedge \text{I} \nabla_{\hat{A}} v + v^\dagger \mathcal{R}^\nabla \text{I} v)$ .

(1.7)

The only change in the subsequent arguments to prove this more general version of Theorem 1.1 is to insert factors of $m_i$'s (which is to say I) in various formulas and check that whenever signs are relevant, either $\text{cl}^\dagger(\text{I}a \otimes a^\dagger)$ appears as $|\text{cl}^\dagger(\text{I}a \otimes a^\dagger)|^2$ or any given $m_i$ appears as $m_i^2$. All of this is straightforward; and this is why only Equation (1.5)'s version of Theorem 1.1 is proved here.

### d) The case when $\dim_\mathbb{C} \mathbb{S}^+ = 4$

If $\dim_\mathbb{C} \mathbb{S}^+ = 2$, then $\mathbb{S}^+$ is the associated self-dual spinor bundle to some $\text{Spin}^\mathbb{C}(4)$ lift of the oriented, orthonormal frame bundle of X. As noted previously, the equations in (1.4) in this case are the original Seiberg-Witten equations and the space of $C^\infty(X; S^1)$ equivalence classes of solutions is sequentially compact.

The next simplest case to consider has $\dim_{\mathbb{C}} \mathbb{S}^+ = 4$. The upcoming Proposition 1.2 talks about this case. The digression that follows directly defines various notions that appear in this proposition.

Let $(\mathbb{S}^+, \mathbb{S}^-, \nabla, cl)$ denote the data for a complex Clifford module. Let $\mathcal{L}$ denote a complex, Hermitian line bundle with a metric compatible connection (denoted by $\mathcal{A}$). An anti-complex homomorphism $C \colon \mathbb{S}^+ \to \mathbb{S}^+ \otimes \mathcal{L}$ is said here to be a *complex conjugation* if it has the following properties:

- $\nabla_{\mathcal{A}} C = C \nabla$ .
- $cl(\cdot) C = C cl(\cdot)$ .
- $\langle C(a), C(b) \rangle = \langle b, a \rangle$ *for all $a, b$ in the same fiber of* $\mathbb{S}^+$.

(1.8)

The first bullet says that $C$ is covariantly constant with respect to the connections on $\mathbb{S}^+$ and on $\mathbb{S}^+ \otimes \mathcal{L}$. The second bullet asserts that $C$ intertwines the respective $\mathbb{S}^+$ and $\mathbb{S}^+ \otimes \mathcal{L}$ versions of the map $cl(\cdot)$ that is defined in (1.2). The third bullet implies that $C$ preserves norms. These bullets imply that $C$ has an inverse that defines a complex conjugation homomorphism from $\mathbb{S}^+ \otimes \mathcal{L}$ to $\mathbb{S}^+$ (which is $(\mathbb{S}^+ \otimes \mathcal{L}) \otimes \mathcal{L}^{-1}$).

If, by chance, the bundle $\mathcal{L}$ has a square root (which will be denoted by $\mathcal{L}^{1/2}$), then $C$ induces a complex conjugation obeying (1.8) from $\mathbb{S}^+ \otimes \mathcal{L}^{1/2}$ to itself with square one; it is an anti-complex linear involution $\mathbb{S}^+ \otimes \mathcal{L}^{1/2}$. (Keep in mind in this regard that $\mathcal{L}^{1/2}$ inherits a canonical Hermitian structure and metric compatible connection from $\mathcal{L}$.)

By way of an example, the Clifford module described previously in Section 1a with $\mathbb{S}^+$ being $(\Lambda^+ \otimes \mathbb{C}) \oplus \underline{\mathbb{C}}$ has the tautological complex conjugation homomorphism to itself which is given by changing the number $i \in \mathbb{C}$ (a chosen square root of -1) to -i. By way of a second example, let $\mathbb{S}_1^+$ and $\mathbb{S}_2^+$ denote the respective self-dual spinor bundles that are associated to two $\mathrm{Spin}^{\mathbb{C}}(4)$ lifts of the principal $SO(4)$ bundle of oriented, orthonormal frames in $TX$. (Give both the sort of connection that was described previously in Section 1a.) As it turns out, the Clifford module with $\mathbb{S}^+ = \mathbb{S}_1^+ \oplus \mathbb{S}_2^+$ has a complex conjugation homomorphism with the bundle $\mathcal{L}$ being the line bundle $(\wedge^4 \mathbb{S}^+)^{-2}$. (This is due two facts: First, any given version of $\mathbb{S}_\Diamond^+$ is isomorphic to the tensor product of any other version with a suitable complex line bundle. Second, the defining, $\mathbb{C}^2$ representation of $SU(2)$ is isomorphic (over $\mathbb{C}$) to its dual.)

The following observation is the 4-dimensional analog of an observation by Haydys and Walpuski [HW] about the dimension 3 version of (1.5).

**Proposition 1.2**: *Let $(Z, \hat{A}, \nu)$ denote the data given by the first and second bullets of an instance of Theorem 1.1 when $\dim_{\mathbb{C}}\mathbb{S}^+ = 4$. The following is true:*

- $(\nabla_{\hat{A}}\nu)^{\dagger} \wedge \nabla_{\hat{A}}\nu = 0$ *and therefore* $F_{\hat{A}} + \frac{1}{|\nu|^2} \nu^{\dagger}\mathcal{R}^{\nabla}\nu = 0$

*Assume that there is a complex conjugation map $C: \mathbb{S}^+ \to \mathbb{S}^+ \otimes \mathcal{L}$ with $\mathcal{L}$ being a complex line bundle over X with Hermitian metric and connection (to be denoted by $\mathcal{A}$).*

- *Let $\hat{A}^{(2)}$ denote the connection that is induced on the bundle $E^2|_{X-Z}$ by the connection $\hat{A}$ on $E|_{X-Z}$. There is an isometric isomorphism from $\mathcal{L}|_{X-Z}$ to $E^2|_{X-Z}$ that identifies the connection $\mathcal{A}$ with the connection $\hat{A}^{(2)}$.*
    a) *The bundle $E|_{X-Z}$ is therefore a square root of the bundle $\mathcal{L}|_{X-Z}$.*
    b) *$F_{\hat{A}} = \frac{1}{2} F_{\mathcal{A}}$.*
- *The set Z has Hausdorff dimension at most 2.*

Section 11 of this paper has the proof of Proposition 1.2. (The proof differs little in substance from the proof in [HW]; the difference is, for the most part, that more background is included here.)

Supposing that $\dim_{\mathbb{C}}(\mathbb{S}^+) = 4$ and that there is a complex conjugation map, the spinor $\nu$ from Proposition 1.2 can be viewed (locally on X) as a $\mathbb{Z}/2$ *harmonic spinor* which is a notion that is defined in [T2]. In particular, theorems in [T2] describe the structure of Z and the behavior of $\nu$. More is said about this in Section 11.

### e) A table of contents for the paper and a look at the proof of Theorem 1.1

This paper has the 11 sections listed directly below. The subsequent preview of the proof of Theorem 1.1 says more about the contents of these sections.

The proof of Theorem 1.1 has much in common with arguments from [T1] and [HW] which, in turn, owe some allegiance to [T3]. These papers discuss solutions to equations that have much the same form as (1.5): They are equations for a connection (to be denoted by A) and a section of a Clifford module bundle (to be denoted by $a$) whose $L^2$ norm is set equal to 1. The common theme in these works is to exploit four observations about any given solution of the relevant version of (1.5). Here is the first:

OBSERVATION 1: *The Bochner-Weitzenboch formula for the operator $D_A$ leads to a priori $L^2$ bounds for the functions $|\nabla_A a|$ and $r^{-1}|F_A^+|$. The Bochner-Weitzenboch formula also leads to $L^\infty$ bounds for $|a|$.*

This observation holds in the context of (1.5) also. Section 2 states and proves precise statements along the lines of Observation 1.

Observation 1 has the following consequence: If $\{r_n, (A_n, a_n)\}_{n \in \{1,2,...\}}$ is a sequence with each $n \in \{1, 2, ...\}$ version of $(A_n, a_n)$ obeying that $r=r_n$ version of (1.5), then a subsequence can be extracted (to be denoted by $\Lambda$) such that $\{|a_n|\}_{n \in \Lambda}$ converges weakly in the $L^2_1$ topology. (What is denoted by $L^2_1$ is the Sobolev space of square integrable functions with square integrable differential.) The limit function is denoted by $|v|$; it is a priori bounded and it can be defined *pointwise* by the rule

$$|v|(p) = \limsup_{n \in \Lambda} |a_n|(p) .$$

(1.8)

To set the notation for the second observation, fix for the moment $r > 1$ and let $(A, a)$ denote a solution to the corresponding version of (1.5). (The Clifford module and the bundle E are fixed once and for all.) Suppose that $p \in X$. A rough measure of the size of $a$ near p at a given (small) length scale r is the average of $|a|^2$ over the boundary of the radius r ball centered at p. A positive function of r to be denoted by $K_p$ is defined by the rule whereby its square is $2\pi^2$ times this average (times a function that is 1 to order $r^2$). A key observation is that the derivative of the function $K_p$ can be written as

$$\tfrac{d}{dr} K_p = \tfrac{N_p}{r} K_p$$

(1.9)

with $N_p$ denoting the function of r given by the rule

$$N_p(r) = \tfrac{1}{r^2 K_p^2} \int_{\text{dist}(\cdot,p)\le r} (|\nabla_A a|^2 + 2r^{-2}|F_A^+|^2) .$$

(1.10)

Since $N_p$ is positive, $K_p$ is increasing. This implies in particular that

$$K_p \geq K_p(0) = \sqrt{2\pi}|a|(p) .$$

(1.11)

This function $N_p$ is the analog of the frequency function that Almgren [Al] and other analysts (for example, [DF], [HHL] and the recent survey [DeL]) have used with great success over the years to study various regularity questions such as harmonic maps, eigenfunctions of the Laplacian, etc.)

OBSERVATION 2: *The pair $(A, a)$ are very well behaved on a small radius ball centered at any given point $p \in X$ where both the $L^2$ norm of $F_A$ is small and $N_p$ is small. In particular, if $r$ is the radius of the ball in question, then*

- *The function $|a|$ is nearly constant on the concentric, radius $\frac{99}{100} r$ ball (its value is nearly $\frac{1}{\sqrt{2\pi}} K_p$).*
- *There are uniform $L^2_1$ bounds for the connection on the concentric, radius $\frac{99}{100} r$ ball (after the application of a suitable automorphism of E.)*
- *After application of the same automorphism, there are uniform $L^2_2$ bounds for $K_p(r)^{-1} a$ on this same concentric ball.*

Section 3 of this paper makes precise statements with regards to this observation.

An important consequence of Observation 2 follows: Return to the sequence $\{r_n, (A_n, a_n)\}_{n \in \Lambda}$. Fix $p \in X$ and suppose that there is a ball centered at $p$ where the small $L^2$ norm condition holds for each $F_{A_n}$; and suppose that the small $N_p$ condition also holds for each $(A_n, a_n)$. Then, by virtue of the second and third bullets of Observation 2, there is a subsequence of $\{(A_n, a_n)\}_{n \in \Lambda}$ that converges nicely on a slightly smaller radius ball after the term-wise application of automorphisms of E.

The preceding begs the following question: Why should the small norm condition for the curvature of $A_n$, and the small $N_p$ condition for $(A_n, a_n)$ hold on a ball with radius *independent* of n? The third and fourth observations are used to answer this question.

OBSERVATION 3: *Fix $p \in X$ and a ball of small radius. Denote the radius by $r$. If the $L^2$ norm of $F_A$ on this ball is small, and if $N_p(r)$ is also small, then the $L^2$ norm of $F_A$ on the concentric smaller radius $\frac{99}{100} r$ ball is a factor $10^{-2}$ or more smaller than it is on the original, radius $r$ ball.*

Section 4 of this paper makes a precise statement of Observation 3. (This observation is a mathematical manifestation of what a physicist might call the 'Higgs mechanism'; which is to say that the effect of the connection A (i.e its curvature) is very short range where $|a|$ is non-zero.) The condition in Observations 2 and 3 that the $L^2$ norm of $F_A$ is 'small' can be quantified by the assertion that the $L^2$ norm is at most $c^{-1}$ with $c$ being a certain number greater than 1. The condition that $N_p$ is small can be likewise quantified.

It proves useful to introduce a radius $r_c(p)$ which is the radius of the ball around p where the $L^2$ norm of $F_A$ is exactly $c^{-1}$. Such a radius exists because A is a smooth connection (unless the integral of $F_A$ over X is also small). Observation 3 leads to some serious tension by virtue of the fact that the $L^2$ norm of $F_A$ on the ball of radius $r_c(p)$ is exactly $c^{-1}$ and it is much less than $c^{-1}$ on a slightly smaller ball. Thus, most of the $L^2$ norm must come from very near the boundary unless $r_c(p)$ is $\mathcal{O}(1)$ or $N_p(r_c(p))$ is $\mathcal{O}(1)$.

The option that $N_p(r_c(p))$ is $\mathcal{O}(1)$ can also lead to tension because of the next observation.

OBSERVATION 4: *Fix* $\varepsilon \in (0, 1]$. *If* r *is small (an upper bound determined by* $\varepsilon$*), and if* $N_p(r)$ *is greater than* $\varepsilon$*, then* $K_p(s)$ *for* $s \geq r$ *obeys* $K_p(s) \leq \kappa\, s^{\varepsilon^3/\kappa}$ *with* $\kappa$ *being* $\mathcal{O}(1)$.

This observation and Observation 3 and (1.11) and the third bullet of (1.5) lead to an a priori positive *lower bound* in terms of $|a|(p)$ for a radius that satisfies the conditions in Observation 2. It also results in a priori Hölder type bounds for $|a|$ along the locus where $|a|$ is zero (which is needed to prove the assertion in Theorem 1.1 to the effect that Z is closed.) How all of this happens is explained in detail in Section 10; but truth be told, the logic in Section 10 is common to [T1], [T3] and [HW] so maybe it can be called a routine argument at this point.

The analog of Observation 4 in [T1], [T3] and [HW] is proved using a monotonicity theorem for the function $N_p$. These monotonicity theorems make assertions to the effect that if $\varepsilon \in (0, 1]$ and $N_p(r) \geq \varepsilon$, then $N_p(s)$ is greater than an $\mathcal{O}(1)$ fraction of $\varepsilon$ if s is greater than r. This then leads to Observation 4 via (1.9). The monotonicity of $N_p$ is proved in these references using a formula for the derivative of their respective versions of $N_p$ whose form was suggested by an analogous formula in [Al].

The formula for the derivative of the current version of $N_p$ is derived in Section 6 of this paper; but it is not evident from this formula that the version here of $N_p$ obeys the required monotonicity. This said, the really hard and novel work in this paper is the proof that Observation 4 holds without a monotonicity formula for $N_p$. Proposition 7.1 makes a formal assertion to this effect. Meanwhile, Sections 7, 8 and 9 of this paper (and also Section 5 and part of Section 6) contain (in total) the proof of Proposition 7.1, which is to say Observation 4.

The U(1) assumption on the connection in (1.5) is a key input to the proof of Proposition 7.1. A short digression is worth taking here to point out the crucial difference between: On the one hand, the 3-dimensional version of (1.5) (the case treated in [HW]) and 4-dimensional, non-Abelian versions that are treated in [T3]; and on the other hand, the case of (1.5) in dimension 4 (with A being a U(1) connection or non-Abelian). The difference is this: The full curvature tensor $F_A$ appears in the top bullet of (1.5) in 3-dimensions, where as only $F_A^+$ appears in 4-dimensions. Meanwhile, the full curvature in dimension 4 is not determined by $F_A^+$ (even non-locally); by way of an example, there are connections with $F_A^+ = 0$ that are not flat. The missing piece to a proof that $N_p$ is monotonic in dimension 4 is control over the anti-self dual part of $F_A$. For the equations in [T1] (which are in 4-dimensions), this missing piece is provided by a Chern-Weil curvature identity that leads to an a priori bound for the $L^2$ norm of the difference between the anti-self dual part of $F_A$ and a quadratic expression in $ra$.

The control of the anti-self dual part of $F_A$ in the U(1) case (which is the case in this paper) exploits a simple monotonicity formula that always holds for closed anti-self dual 2-forms. (This monotonicity formula says (to a first approximation) that the $L^2$ norm of the 2-form on a radius r ball is no greater than $\frac{1}{4}$ times its $L^2$ norm on the concentric 2r ball.) This anti-self dual 2-form monotonicity can be used to say something about $F_A$ when the $L^2$ norm of its anti-self dual part on any given ball is very much greater than the $L^2$ norm of its self dual part (which is controlled in any event by $N_p$). The point here is that $F_A$ is very well approximated by a closed, anti-self dual 2-form when its anti-self dual part on a given ball is much bigger than its self dual part (as measured by $L^2$ norms).

However, be forwarned that the derivation of Observation 4 from this observation about closed, anti-self dual 2-forms (plus the derivative formula for $N_p$) is a lot of work; it occupies more than half of the paper (Sections 5, part of Section 6 and Sections 7-9.)

All of the pieces (Observations 1-4) are put together in Section 10 to prove Theorem 1.1. Section 11 proves Proposition 1.2.

**f) Conventions and notation**

The notation will be simpler if $\mathbb{S}_E^+$ and $\mathbb{S}_E^-$ are used henceforth to denote what the previous sections denoted by $\mathbb{S}^+ \otimes E$ and $\mathbb{S}^- \otimes E$.

Many computations are best done using a local orthonormal frame. Supposing that $\{e^i\}_{i=1,2,3,4}$ is an oriented orthonormal frame for T*X defined on a neighborhood of a point p, then a covector $v \in T^*X|_p$ is written with respect to this frame as $v_i e^i$ with the summation over the repeated indices to be understood implicitly. This convention whereby repeated indices are implicitly summed is used throught this article; summation over repeated indices should be assumed unless stated specifically to the contrary. With regards to indices, covariant derivatives in the direction of the basis vectors are denoted by either $\nabla_i$ or $\nabla_{A,i}$, the latter depending on the indicated connection A and the former being independent of A. Another point with indices: An oriented, orthonormal frame for T*X at a given point, $\{e^i\}_{i=1,2,3,4}$ determines a canonical, oriented, orthonormal frame for $\Lambda^+$ at that point. This is the frame $\{\omega^a\}_{a=1,2,3}$ given by the rule

$$\omega^1 = \tfrac{1}{\sqrt{2}}(e^2 \wedge e^3 + e^1 \wedge e^4), \quad \omega^2 = \tfrac{1}{\sqrt{2}}(e^3 \wedge e^1 + e^2 \wedge e^4), \quad \omega^3 = \tfrac{1}{\sqrt{2}}(e^1 \wedge e^2 + e^3 \wedge e^4).$$
(1.12)

(Note that the volume 4-form is $e^1 \wedge e^2 \wedge e^3 \wedge e^4$.) A section $s$ of $\Lambda^+$ at the point can be written using this frame as $s = s_a \omega^a$ with the summation over the repeated indices to be understood implicitly.

It also proves convenient to introduce $\langle\,,\,\rangle$ to signify the inner product on all Hermitian vector bundles and the metric inner product on all tensor bundles; and their mutual tensor products. For example, $v^\dagger \mathcal{R}^\nabla v$ appears in Item d) of the second bullet of Theorem 1.1, and this is written henceforce as $\langle v, \mathcal{R}^\nabla v\rangle$. But, be forwarned that there is some opportunity for confusion because, $v^\dagger \nabla_A v$ in Item c) of the second bullet of Theorem 1.1 is written below in shorthand as $\langle v, \nabla_A v\rangle$ which should be interpreted as

$$\langle v, \nabla_A v\rangle = \langle v, \nabla_{A,i} v\rangle e^i$$
(1.13)

(referring to the convention of the previous paragraph). Also, $(\nabla_{\hat{A}} v)^\dagger \wedge \nabla_{\hat{A}} v$ from Item d) of Theorem 1.1 will be written in shorthand as $\langle \nabla_{\hat{A}} v \wedge \nabla_{\hat{A}} v\rangle$ which means

$$\langle \nabla_{\hat{A}} v \wedge \nabla_{\hat{A}} v\rangle = \langle \nabla_{\hat{A},i} v, \nabla_{\hat{A},j} v\rangle e^i \wedge e^j.$$
(1.14)

Another convention that is used below involves numbers that arise in various inequalities whose precise value has no real bearing on the matters at hand; the important thing being that the number is independent of a given value of $r$ and a given pair $(A, a)$ from (1.5) and, if relevant, a given point in X and a given distance from this point. Numbers of this sort are denoted by $\kappa$ in propositions and lemmas, and by $c_0$ in the text. These numbers are always greater than 1 and they can be assumed to increase between successive appearances. Given a number $\mu \in (0, \tfrac{1}{100})$, then numbers $\kappa_\mu$ and $c_\mu$ are like $\kappa$ and $c_0$ except that they can also depend on the choice for $\mu$. In particular, the incarnations of $\kappa_\mu$ and $c_\mu$ are greater than 1 and they also increase between successive appearances.

By way of motivation for the final bit of notation, a smooth 'bump' function is needed in many places below; this being a function with values in [0, 1] that is equal to 1 on a given set and equal to zero on a slightly larger set. In the contexts below, such a function is constructed from a universal bump function on $\mathbb{R}$ that is denoted by $\chi$. This function $\chi$ is non-increasing, equal to 1 on $(-\infty, \tfrac{1}{2}]$ and equal to 0 on $[\tfrac{3}{4}, \infty)$.

## 2. Fundamental identities

This section summarizes some basic identities and inequalities that play key roles in the proof of Theorem 1.1.

**a) Function space inequalities**

The results here and in subsequent sections make use of the basic function space inequalities that are stated below in (2.1). To set the stage and notation, fix $p \in X$ and given $r \in (0, c_0^{-1}]$ with the upper bound here being far less than the radius of a Gaussian coordinate chart centered at p. Use $B_r$ to denote the open, radius r ball centered at p. The notation has $L_1^2$ denoting the space of square integrable functions on $B_r$ with square integrable differential. Supposing that $f$ is an $L_1^2$ function on $B_r$, then both $f^2$ and $\text{dist}(\cdot, p)^{-1} f$ are square integrable on $B_r$ with norms that obey the following inequalities:

- $(\int_{B_r} f^4)^{1/2} < c_0 \int_{B_r} (|df|^2 + \frac{1}{r^2}|f|^2)$ .

- $\int_{B_r} \frac{1}{\text{dist}(p, \cdot)^2} f^2 < c_0 \int_{B_r} (|df|^2 + \frac{1}{r^2}|f|^2)$ .

- *If f has compact support in* $B_r$ *or if* $\int_{B_r} f = 0$, *then both of the preceding inequalities hold with the right hand side being* $c_0 \int_{B_r} |df|^2$ .

(2.1)

The top bullet in (2.1) is a Sobolev inequality and the second bullet is Hardy's inequality. For $k \in [1, \infty)$, let $L^k$ denote the space of measurable functions on $B_r$ with k'th power of the norm being integrable. It follows from (2.1) that $L_1^2$ functions are in $L^k$ for $k \leq 4$ and the resulting map from $L_1^2$ to $L^k$ is bounded. This map is also compact when $k < 4$. (See, e.g. [AF] for a general reference on Sobolev inequalities.)

Suppose that $|\cdot|$ is the induced norm of a metric on a vector bundle over $B_r$ and that $\nabla$ is a metric compatible covariant derivative. Supposing that $\mathfrak{w}$ is a measurable section of the bundle, then $|\mathfrak{w}|$ is an $L_1^2$ function on $B_r$ if $\mathfrak{w}$ and $\nabla \mathfrak{w}$ are square integrable with it understood that the norm of $\nabla \mathfrak{w}$ is defined using the norm on the vector bundle and the norm defined by the Riemannian metric. This is so because the norm of $d|\mathfrak{w}|$ is bounded almost everywhere by that of $|\nabla \mathfrak{w}|$. By way of a consequence, it follows that the inequalities in (2.1) have implications for sections of a vector bundles over $B_r$ with a metric and metric compatible covariant derivative.

The notation in what follows uses $\|\cdot\|_2$ to indicate the $L^2$ norm of a given function or given section of a vector bundle. In general, for $p \geq 1$, the notation has $\|\cdot\|_p$ denoting the $L^p$ norm of the function or section.

**b) Algebraic identities**

Certain algebraic identities involving Clifford multiplication lie behind many of the identities that are presented in the subsequent subsections. To say more, fix $p \in X$ and let $\{\omega^a\}_{a=1,2,3}$ denote an oriented, orthonormal basis for $\Lambda^+|_p$. Since the Clifford module is self-dual, the homomorphims cl whose adjoint appears in (1.4) and (1.5) is determined at the point p by its values on the basis elements $\{\omega^a\}_{a=1,2,3}$. With this fact in mind, define the endomorphisms $\{\tau^a\}_{a=1,2,3}$ of $\mathbb{S}^+$ (and thus $\mathbb{S}_E^+$) by the rule

$$\tau^a = \tfrac{1}{2\sqrt{2}} cl(\omega^a) .$$

(2.2)

It follows now from (1.1) that these endomorphisms act like imaginary quaternions:

$$\tau^a \tau^b = -\delta^{ab} - \varepsilon^{abc} \tau^c ,$$

(2.3)

where $\delta^{ab}$ is to equal 1 when $a = b$ and 0 otherwise, and where $\varepsilon^{abc}$ is completely anti-symmetric and such that $\varepsilon^{123} = 1$. The formula in (2.2) implies that the term $cl^\dagger(a \otimes a^\dagger)$ that appears in (1.4) and (1.5) can be written as

$$cl^\dagger(a \otimes a^\dagger) = 2\sqrt{2} \langle a, \tau^a a \rangle \omega^a .$$

(2.4)

Note that when $F_A^+$ is written in terms of its components as $F_{A\,a}^+ \omega^a$, then the equation in the top bullet of (1.5) asserts that

$$F_{A\,a}^+ = \tfrac{1}{\sqrt{2}} r^2 \langle a, \tau^a a \rangle + \mathfrak{w} \quad \text{for each } a \in \{1, 2, 3\}.$$

(2.5)

By way of notation, the self-dual 2-form $\langle a, \tau^a a \rangle \omega^a$ is written henceforth as $\langle a, \tau a \rangle$.

A second useful identity comes from the middle equation in (1.5): Fix an index i from $\{1, 2, 3, 4\}$ and act on this equation by $\mathfrak{c}^\dagger(e^i)$ to see that $a$ obeys

$$\nabla_{A,i} a + \tfrac{1}{2} cl(e^i \wedge e^j) \nabla_{A,j} a = 0 ,$$

(2.6)

which can be written in turn as

$$\nabla_{A,i} a + \sqrt{2} \langle \omega^a, e^i \wedge e^j \rangle \tau^a \nabla_{A,j} a = 0 .$$

(2.7)

The preceding identities will be used at times with little by way of comment.

**c) Weitzenboch formula for *a***

Fix a connection on E (to be denoted by A) so as to define the Dirac operator $D_A$ using (1.3). Let $D_A^\dagger$ denote its formal $L^2$ adjoint. The Bochner-Weitzenbock formula for the operator $D_A^\dagger D_A$ can be written as

$$D_A^\dagger D_A a = \nabla_A^\dagger \nabla_A a - \tfrac{1}{2} \mathrm{cl}(F_A) a + \mathfrak{R}^\diamond a$$

(2.8)

with $\nabla_A^\dagger$ denoting the formal $L^2$ adjoint of the covariant derivative $\nabla_A$, and with $\mathfrak{R}^\diamond$ denoting an automorphism of $\mathbb{S}^+$ that is independent of the chosen connection A and the section $a$ and the bundle E. With (2.8) understood, suppose now that (1.5) holds for a given $r$. Because $D_A a = 0$, so $D_A^\dagger D_A a = 0$. This being the case, (2.8) with the top bullet in (1.5) and (2.2) and (2.5) lead to a second order equation for $a$ having the form

$$\nabla_A^\dagger \nabla_A a - r^2 \langle a, \tau^a a \rangle \tau^a a + \mathfrak{R} a = 0$$

(2.9)

with $\mathfrak{R}$ denoting $\mathfrak{R}^\diamond - \tfrac{1}{2} \mathrm{cl}(\omega)$. Taking the inner product of this equation with $a$ (and using the fact that the covariant derivatives are metric compatible) leads to the following differential identity for $|a|^2$:

$$\tfrac{1}{2} d^\dagger d |a|^2 + |\nabla_A a|^2 + r^2 |\langle a, \tau a \rangle|^2 + \langle a, \mathfrak{R} a \rangle = 0.$$

(2.10)

This last equation in particular leads to the two inequalities that are stated by the lemma below.

**Lemma 2.1**: *There exists* $\kappa > 1$ *with the following significance: Fix* $r > 0$ *and let* $(A, a)$ *denote a pair of connection on* E *and section of* $\mathbb{S}_E^+$ *obeying (1.5).*

- $\int_X (|\nabla_A a|^2 + r^2 |\langle a, \tau a \rangle|^2) \le \kappa$.
- $\sup_{p \in X} |a|(p) \le \kappa$.
- *If* $p \in X$, *then* $\int_X \frac{1}{\mathrm{dist}(p,\cdot)^2} (|\nabla_A a|^2 + r^2 |\langle a, \tau a \rangle|^2) \le \kappa$.

Noting that $F_A^+ = \tfrac{1}{\sqrt{2}} r^2 \langle a, \tau a \rangle + \omega$, the integral bounds in this lemma give a priori bounds for integrals of $r^{-2} |F_A^+|^2$. In particular, these bound for integrals over X imply in turn that

$$\int_X |F_A|^2 \le c_0 (r^2 + 1).$$

(2.11)

Indeed, this follows by virtue of two facts: First, if A is any Hermitian connection on E, then the integral of $-\tfrac{1}{4\pi^2} F_A \wedge F_A$ gives the pairing between the square of the first Chern class of E and the fundamental class of X. Second, $-*(F_A \wedge F_A) = |F_A^+|^2 - |F_A^-|^2$ for any connection on E. Granted these facts, then the $L^2$ bounds for $r^{-1}F_A^+$ leads to (2.11).

*Proof of Lemma 2.1*: To obtain the first bullet's inequality, integrate the equation in (2.10) over X, and then use the fact that $|\Re| \le c_0$ with the third bullet of (1.5) to bound the norm of the integral of $\langle a, \Re a \rangle$. To prove the second and third inequalities, let $G_p$ denote the Green's function on X for the operator $d^\dagger d + 1$ with pole at p. This $G_p$ is a positive, smooth function on X−p that obeys the following bounds:

$$c_0^{-1} \frac{1}{\mathrm{dist}(p,\cdot)^2} < G_p \le c_0 \frac{1}{\mathrm{dist}(p,\cdot)^2} \quad and \quad |\nabla G_p| \le c_0 \frac{1}{\mathrm{dist}(p,\cdot)^3} \ .$$

(2.12)

Multiply both sides of (2.10) by $G_p$ and then integrate the resulting identity over X. Use integration by parts and the upper and lower bounds for $G_p$ in (2.12) with fact that $|\Re| \le c_0$ to obtain the inequality

$$\tfrac{1}{2} |a|^2(p) + \int_X \frac{1}{\mathrm{dist}(p,\cdot)^2} (|\nabla_A a|^2 + r^2 |\langle a, \tau a \rangle|^2) \le c_0 \int_X \frac{1}{\mathrm{dist}(p,\cdot)^2} |a|^2 \ .$$

(2.13)

To see that the right hand integral in this last inequality is less than $c_0$, first invoke the middle bullet of (2.1) on the radius $c_0^{-1}$ ball centered at p to bound the right hand integral by $c_0$ times the integral over X of $|d|a||^2 + |a|^2$. The latter integral is itself bounded by $c_0$ because $|d|a|| \le |\nabla_A a|$ and because of what is said by the first bullet of the lemma.

Lemma 2.1 has the following immediate corollary:

**Proposition 2.2**: *Let $\{r_n, (A_n, a_n)\}_{n=1,2,\ldots}$ denote a sequence such that for each positive integer n, what is denoted by $r_n$ is a positive number and what is denoted by $(A_n, a_n)$ is a pair of connection on E and section of $\mathbb{S}_E^+$ obeying the $r = r_n$ version of (1.5). Suppose that the sequence of numbers $\{r_n\}_{n=1,2,\ldots}$ has a bounded subsequence. There is a subsequence $\Theta \subset \{1, 2, \ldots\}$ and a corresponding sequence $\{g_n\}_{n \in \Theta}$ of hermitian automorphisms of the line bundle E such that the sequence $\{(A_n - g_n^{-1}dg_n, g_n a_n)\}_{n \in \Theta}$ converges in the $C^\infty$ topology to a solution to some $r > 0$ version of (1.5).*

*Proof of Proposition 2.2*: The assumptions about $\{r_n\}_{n=1,2,\ldots}$ imply that there is a subsequence $\Xi \subset \{1, 2, \ldots\}$ such that the corresponding subsequence $\{r_n\}_{n \in \Xi}$ converges. Let $r$ denote the limit. It follows from (2.11) that there is an a priori bound on the $L^2$ norms of the $\{F_{A_n}\}_{n\in\Xi}$. One can now argue from this using the top bullet in Lemma 2.1 (along the lines of the argument in Ben Mare's thesis [Ma]), and using the inequalities in (2.1), that there is a subsequence $\Theta \subset \Xi$ and a sequence $\{u_n\}_{n\in\Theta}$ of Hermitian automorphisms of E with the following property: The associated sequence of connections $\{A_n - g_n^{-1}dg_n\}_{n\in\Theta}$ converges weakly in the $L^2_1$ topology to an $L^2_1$ connection on E (to be denoted by A), and the sequence $\{g_n a_n\}$ converges weakly in the $L^2_1$ topology to a section of $\mathbb{S}_E^+$ (to be denoted by $a$). Moreover, this pair $(A, a)$ obeys the equations in (1.5) with $r$ as just described. (These are first order equations with quadratic nonlinearities and so they make sense for $L^2_1$ pairs of connection on E and section of $\mathbb{S}_E^+$.) At this junction, rather standard elliptic bootstrap arguments (see again Mare's thesis [Ma]) can be used to conclude that the pair $(A, a)$ is smooth and that the convergence of $\{(A_n - g_n^{-1}dg_n, g_n a_n)\}_{n\in\Theta}$ to $(A, a)$ is in the $C^\infty$ topology. The details for all of this are straightforward and therefore omitted.

**d) Differential equations for the curvature 2-form and the covariant derivative of $a$**

There are other equations that play roles in the subsequent story. The first set is an equation for the curvature 2-form $F_A$:

$$*d*F_A = r^2(\langle a, \nabla_A a\rangle - \langle \nabla_A a, a\rangle) + *d*\omega \, .$$

(2.14)

By way of notation, $*$ denotes the Riemannian metric's Hodge star. By way of a reminder, $\langle a, \nabla_A a\rangle$ denotes the $\mathbb{C}$-valued 1-form on X that is defined as follows: Let $p \in$ X be given and let $\{e^i\}_{i=1,2,3,4}$ denote an oriented, orthonormal frame for $T^*X|_p$. Then $\langle a, \nabla_A a\rangle$ at p is $\langle a, \nabla_{A,i} a\rangle e^i$ when written using this frame. The identity in (2.14) is derived by first noting that $d*F_A = 2dF_A^+$ because $dF_A$ is zero. The 2-form $dF_A^+$ (and thus $dF_A$) can be written in terms of $a$ and $\nabla_A a$ using (2.5) and (2.7). This results in (2.14).

Taking the exterior derivative of (2.14) equation leads in turn to a second order differential equation for $F_A$ that has the form

$$d*d*F_A - 2r^2|a|^2 F_A = 2r^2 \langle \nabla_A a \wedge \nabla_A a\rangle + \mathfrak{G}.$$

(2.15)

The notation here uses $\langle \nabla_A a \wedge \nabla_A a\rangle$ shorthand for the $i\mathbb{R}$ valued 2-form on X that can be written using an orthonormal frame $\{e^i\}_{i=1,2,3,4}$ for $T^*X$ as $\langle \nabla_{A,i} a, \nabla_{A,j} a\rangle e^i \wedge e^j$. Meanwhile, $\mathfrak{G}$ denotes $d*d*\omega$. Because $dF_A = 0$, and thus $*d*dF_A = 0$, the identity in (2.15) can also be written using the Laplacian on 2-forms as

$$\nabla^\dagger \nabla F_A + 2r^2|a|^2 F_A + \mathcal{R} \cdot F_A = -2r^2 \langle \nabla_A a \wedge \nabla_A a \rangle + \mathfrak{S}$$
(2.16)

with $\mathcal{R}$ denoting an endomorphism that comes from the Riemann curvature tensor. In particular, $|\mathcal{R}| \leq c_0$. Note also that $\mathcal{R}$ commutes with the Hodge star and thus it maps $\Lambda^+$ to $\Lambda^+$, and it maps $\Lambda^-$ to $\Lambda^-$.

A useful second order differential equation for $\nabla_A a$ is obtained from (2.9) by taking the A-covariant derivative of both sides. Commuting derivatives leads to the next equation. This is an equation for the directional derivative of $\nabla_A a$ along a basis vector of a local orthonormal frame for TX.

$$\nabla_A^\dagger \nabla_A (\nabla_{A,i} a) - F_{A,ij} \nabla_{A,j} a - \nabla_{A,j}(F_{A,ij} a) - r^2 \nabla_i \langle a, \tau^a a \rangle) \tau^a a - r^2 \langle a, \tau^a a \rangle) \tau^a \nabla_{A,i} a$$
$$+ \mathfrak{R}_{0,i} a + \mathfrak{R}_{1,ij} \nabla_{A,j} a = 0 .$$
(2.17)

Here, $\mathfrak{R}_0$ and $\mathfrak{R}_1$ are homomorphisms that are independent of the pair $(A,a)$ and $r$. (Their norms, in particular, are bounded by $c_0$.) This identity in (2.16) is used after taking the inner product of both sides with $\nabla_A a$ to obtain an inequality that has the form

$$\tfrac{1}{2} d^\dagger d |\nabla_A a|^2 + |\nabla_A \nabla_A a|^2 + \tfrac{1}{2} r^2 |\nabla \langle a, \tau^a a \rangle|^2 - 2F_{A,ij}, \langle \nabla_{A,i} a, \nabla_{A,j} a \rangle - \langle \nabla_{A,i} a, a \rangle \nabla_j F_{A,ij}$$
$$- r^2 \langle a, \tau^a a \rangle \langle \nabla_{A,i} a, \tau^a \nabla_{A,i} a \rangle \leq c_0 (|\nabla_A a| |a| + |\nabla_A a|^2)$$
(2.18)

This inequality and (2.16) can be used to derive a priori bounds for the integrals over X of $|\nabla_A \nabla_A a|^2$ and $|\nabla F_A|^2$. However, such global bounds don't play much of a role in what follows. The equations in (2.16) and (2.18) are used instead (in subsequent sections) to bound norms of $\nabla_A \nabla_A a$ and $\nabla F_A$ over small radius balls in X.

### 3. Implications of small covariant derivative and small curvature

Fix $r > 1$ and let $(A, a)$ denote a pair of connection on E and section of the bundle $\mathbb{S}_E^+$ that obeys (1.5). The propositions in this section describe the behavior of the pair $(A, a)$ on balls in X where the integral of the curvature of A is relatively small.

### a) The number $r_{cF}$

A number that is defined momentarily will be used to measure the size of $F_A$ near a given point in X. The definitions of this number requires the a priori specification of a positive number to be denoted by $c$. The convention takes $c > 100$. With $c$ chosen, fix $p \in U$. Supposing that $r \in (0, c_0^{-1}]$, define $B_r$ to be the (open) ball of radius r centered at p. (The version of the number $c_0$ is chosen so that $B_r$ is well inside a Gaussian coordinate chart centered at p.) With this notation understood, the number $r_{cF}$ is defined to be the largest of the numbers $r \in (0, c_0^{-1}]$ such that

$$\int_{B_r} |F_A|^2 \leq c^{-2} .$$

(3.1)

The number $r_{cF}$ depends on the chosen point p but since p is fixed in most of what follows, this dependence is not indicated by the notation. Some of the upcoming propositions in Section 3 assert a priori bounds for $a$ and $F_A$ on balls centered at p with radius less than $r_{cF}$.

**b) The functions K and N**

The constructions in this section mimic constructions in Section 3a of [T1]. To set the stage, suppose as before that $r > 1$ and that $(A, a)$ are a pair of connection on E and section of the bundle $\mathbb{S}_E^+$ that obey (1.5). This data is used with a given point $p \in U$ to define a pair of positive functions on $[0, c_0^{-1}]$ to be denoted by K and N. The notation in what follows has $\partial B_r$ for $r \in (0, c_0^{-1}]$ denoting the boundary of the closed, radius r ball centered at the chosen point p.

The definition of K starts with the definition of the function h on $[0, c_0^{-1}]$ whose value at any given $r \in [0, c_0^{-1}]$ is the integral of $|a|^2$ on $\partial B_r$; this is to say that

$$h(r) = \int_{\partial B_r} |a|^2 .$$

(3.2)

Arguments that differ only cosmetically from those used by Aronszajn [Ar] can be used with (1.5) to prove that $a$ can not vanish on an open set. This implies in particular that h is positive on $(0, c_0^{-1}]$.

The definition of K requires a second function, this denoted by ∂. The definition of ∂ is in the upcoming equation (3.3). What is denoted by M in this equation is defined by writing the trace of the second fundamental form of $\partial B_r$ for any given $r \in (0, c_0^{-1}]$ as $\frac{3}{r} + M$. The function ∂ is defined by the rule

$$r \to \partial(r) = \int_0^r (\frac{1}{h(s)} (\int_{B_s} \langle a, \Re a \rangle + \tfrac{1}{2} \int_{\partial B_s} M|a|^2)) \, ds .$$

(3.3)

The norm of ∂ is bounded by $c_0 r^2$. The proof that this is so differs only in notation from the proof in Section 3a of [T1] that ∂'s namesake in [T1] has norm bounded by $c_0 r^2$. Note for future reference that the arguments for the assertion that $|\partial| \leq c_0 r^2$ lead as they did in [T1] to the inequality $h(r) \geq (\frac{r}{s})^3 e^{-c_0(r^2-s^2)} h(s)$ when r is greater than s and both are from $(0, c_0^{-1})$. As was the case in [T1], the latter inequality implies in turn that

$$\int_{B_r} |a|^2 \leq 4 e^{c_0 r^2} r h(r) .$$

(3.4)

when $r \in [0, c_0^{-1}]$.

Define $K$ to be the positive square root of $\frac{1}{r^3} e^{-2\mathfrak{d}} h$, thus

$$K(r) = \frac{1}{r^{3/2}} e^{-\mathfrak{d}} \sqrt{h(r)} .$$

(3.5)

And, with $K$ in hand, the function $N$ is then defined by the rule

$$r \to N(r) = \frac{1}{r^2 K(r)^2} \int_{B_r} (|\nabla_A a|^2 + r^2 |\langle a, \tau a \rangle|^2) .$$

(3.6)

Note that $N(r)$ can also be written (because of the top bullet to (1.5)) as

$$N(r) = \frac{1}{r^2 K(r)^2} \int_{B_r} (|\nabla_A a|^2 + 2r^{-2} |F_A^+|^2) .$$

(3.7)

In any event, the definition of $\mathfrak{d}$ is chosen so that $N$ and $K$ are related via the identity

$$\frac{d}{dr} K = \frac{1}{r} N K .$$

(3.8)

The proof of (3.8) is identical to the proof an analogous identity in Section 3a of [T1]. Therefore, supposing that $r \in (0, c_0^{-1})$ and $s \in (0, r]$, then $K(r)$ and $K(s)$ are related via

$$K(r) = \exp(\int_s^r \frac{N(\tau)}{\tau} d\tau) K(s)$$

(3.9)

which is the integral form of (3.8).

By way of a convention: Let $\underline{K}$ denote the positive square root of the function on $(0, c_0^{-1}]$ whose value at any given r is the average of $|a|^2$ on the radius r sphere centered at p. This function $\underline{K}$ differs little from $K$ in as much as $|K - \underline{K}|$ and $|r \frac{d}{dr} \underline{K} - N K|$ are both bounded by $c_0 r^2 K$. As a consequence, the function $K$ can be used as a proxy for $\underline{K}$ in subsequent arguments in this paper if one remembers to add an appropriately placed error term of order $c_0 r^2 K$. In particular, the convention in what follows is to use $\underline{K}$ and $K$ interchangebly without notational distinction or further comment (for the most part). The appearance of $\mathcal{O}(c_0 r^2)$ error terms in subsequent equations serves as the (implicit) justification. The function K is preferred over <u>K</u> because the frequency function N that appears in (3.7) is manifestly positive.

Analogs of the function N were introduced in [Al], [DF] and [HHL] where they were used to study singular level sets of solutions to elliptic equations. The applications of N here are guided by the applications in the latter references.

**c) The pointwise norm of $\mathfrak{a}$**

The proposition below says that the pointwise norm of $|\mathfrak{a}|$ on a ball in X is bounded by a multiple of the value of K on a slightly larger ball. It also says (roughly) that $|\mathfrak{a}|$ can not be much larger than this value of $\frac{1}{\sqrt{2\pi}}$ K when the corresponding value of N is small, and that it can't be much less except (possibly) on a set of small volume.

**Proposition 3.1**: *There exists $\kappa > 1$; and given $\mu \in (0, \frac{1}{100}]$, there exists $\kappa_\mu > \kappa$; these numbers having the following significance: Fix $r > 1$ and suppose that $(A, \mathfrak{a})$ is a connection on E and section of $\mathbb{S}_E^+$ that obeys (1.5). Fix $p \in X$ so as to define the functions K and N. Fix $r \in (0, \kappa^{-1}]$, then*

- $|\mathfrak{a}| \leq \kappa_\mu \frac{1}{\sqrt{2\pi}} K(r)$ *on* $B_{(1-\mu)r}$.

*Moreover, if $N(r) \leq 1$, then*

- $|\mathfrak{a}| \leq (1 + \kappa_\mu \sqrt{N(r)}) \frac{1}{\sqrt{2\pi}} K(r)$ *on* $B_{(1-\mu)r}$.
- *The fraction of the volume of $B_{(1-\mu)r}$ where $|\mathfrak{a}| < \frac{3}{4} \frac{1}{\sqrt{2\pi}} K(r)$ is at most $\kappa_\mu N(r)^{3/2}$.*

The rest of this subsection is occupied with the proof of this proposition.

By way of notation, keep in mind that $c_\mu$ is used to denote a number that is greater than 1 and depends only on $\mu$ and the data giving the Clifford module; and that it can be assumed to increase between successive appearances. The proof also introduces a non-negative function on $B_r$ denoted by $\chi_\mu$. This function $\chi_\mu$ is 0 where $\mathrm{dist}(\cdot, p) \geq (1 - \frac{5}{8}\mu)r$, it equals 1 where $\mathrm{dist}(\cdot, p) \leq (1 - \frac{3}{4}\mu)r$; and its derivatives obey the bounds $|d\chi_\mu| \leq c_0 \frac{1}{\mu r}$, and $|\nabla d\chi_\mu| \leq c_0 \frac{1}{\mu^2 r^2}$. The function $\chi_\mu$ can be constructed using the fiducial function $\chi$.

*Proof of Proposition 3.1*: The proof has three parts. Part 1 of the proof derives the generic upper bound for $|\mathfrak{a}|$ in the first bullet of the proposition. The upper bound that is asserted by the proposition's second bullet is proved in Part 2; and Part 3 explains why the third bullet of the proposition is true

*Part 1*: There are three steps in this part of the proof. These steps go from an $L^2$ bound on $|a|$ on the ball $B_r$ to an $L^2_1$ bound on $B_{(1-\mu/8)r}$ to a pointwise bound on $B_{(1-\mu/2)r}$.

Step 1: The promised $L^2$ bound asserts

$$\int_{B_r} |a|^2 \leq \tfrac{1}{4}(1+c_0 r^2) r^4 \kappa(r)^2 \ .$$

(3.10)

To prove this, do the integral by first integrating $|a|^2$ over the constant $s \in (0, r]$ spheres with center p and then integrating over s. The integral of $|a|^2$ over a constant s sphere differs from $\kappa^2(s)$ by at most $c_0 s^2 \kappa(s)^2$. The bound in (3.10) follows from the fact that $\kappa(s)$ is less than $\kappa(r)$ when $s < r$.

Step 2: Multiply both sides of the identity in (2.10) by $\chi_{\mu/64}$ and then integrate the result over $B_r$. An integration by parts on the left leads to the bound

$$\int_{B_r} \chi_{64\mu}(|\nabla_A a|^2 + r^2|\langle a, \tau a\rangle|^2) \leq c_\mu \tfrac{1}{r^2} \int_{B_r} |a|^2 \ .$$

(3.11)

Noting that $|\nabla_A a| \leq |d|a||$ (because $\nabla_A$ is metric compatible), the bounds in (3.11) and (3.10) leads to a $c_\mu r^2$ bound for the integral of $|d|a||^2$ on $B_{(1-\mu/8)r}$.

Step 3: Supposing that $q \in B_{(1-\mu/2)r}$, let $G_q$ denote the Dirichelet Green's function for the Laplacian (this is $d^\dagger d$) on $B_r$ with pole at the point q. This function is positive on the interior of $B_r$, it vanishes on the boundary and it obeys

$$|G_q| \leq c_0 \frac{1}{\text{dist}(q,\cdot)^2} \quad \text{and} \quad |\nabla G_q| \leq c_0 \frac{1}{\text{dist}(q,\cdot)^3} \ .$$

(3.12)

Multiply both sides of (2.10) by $\chi_{\mu/4} G_q$ and then integrate by parts on the left hand side of the resulting identity to see that

$$\tfrac{1}{2}|a|^2(q) + \int_{B_r} \chi_{\mu/4} G_q (|\nabla_A a|^2 + r^2|\langle a, \tau a\rangle|^2) = -\int_{B_r} \chi_{\mu/4} G_q \langle a, \Re a\rangle$$
$$-\tfrac{1}{2}\int_{B_r} (d^\dagger d\chi_{\mu/4} G_q + 2\langle d\chi_{\mu/4}, dG_q\rangle) |a|^2 \ .$$

(3.13)

The absolute value of the integral of $\chi_{\mu/4} G_q \langle a, \Re a\rangle$ that appears on the right hand side of (3.13) is no greater than $c_\mu$ times the integral of $|d|a||^2 + \tfrac{1}{r^2}|a|^2$ on $B_{(1-\mu/8)r}$ (this follows from the left most bound in (3.12) and (2.1)). Granted such a bound, then it follows from what is said by Steps 1 and 2 that the absolute value of the left most integral on the right hand side of (3.13) is no greater than $c_\mu r^2 K(r)^2$. Meanwhile, (3.12) and (3.10) lead directly to a $c_0 K(r)^2$ bound for the norm of the right most integral on the right hand side of (3.13). Therefore, the right hand side of (3.13) is no greater than $c_\mu K(r)^2$. Since $G_q$ is positive, it follows that $|a|^2$ on $B_{(1-\mu/2)r}$ is also bounded by $c_\mu K^2(r)$. This is the generic bound asserted by the proposition's top bullet.

*Part 2*. This part of the proof discusses the second bullet of Proposition 3.1. The proof of this bullet has three steps also.

<u>Step 1</u>: Multiplying both sides of the equation $d^\dagger d(1) = 0$ by $\chi_{\mu/4} G_q$ leads to the identity

$$1 = \int_{B_r} (d^\dagger d\chi_{\mu/4} G_q + 2\langle d\chi_{\mu/4}, dG_q\rangle)$$

(3.14)

Multiply both sides of this by $\frac{1}{4\pi^2} K(r)^2$ and subtract the result from (3.13) to see that

$$\tfrac{1}{2}(|a|^2(q) - \tfrac{1}{2\pi^2} K(r)^2) + \int_{B_r} \chi_{\mu/4} G_q (|\nabla_A a|^2 + r^2|\langle a, \tau a\rangle|^2) = -\int_{B_r} \chi_{\mu/4} G_q \langle a, \Re a\rangle$$
$$- \tfrac{1}{2} \int_{B_r} (d^\dagger d\chi_{\mu/4} G_q + 2\langle d\chi_{\mu/4}, dG_q\rangle)(|a|^2 - \tfrac{1}{2\pi^2} K(r)^2).$$

(3.15)

Granted this equation and granted the bounds in (3.12), and granted that the absolute value of the integral of $\chi_{\mu/4} G_q \langle a, \Re a\rangle$ is at most $c_\mu r^2 K(r)^2$, then the assertion in the top bullet of the lemma follows from the bound

$$\int_{B_{(1-\mu/8)r} - B_{(1-\mu/4)r}} \left||a|^2 - \tfrac{1}{2\pi^2} K(r)^2\right| \leq c_\mu \sqrt{N(r)} K(r)^2 r^4.$$

(3.16)

This bound is established in the next two steps.

<u>Step 2</u>: Let $\nabla_{A,r}$ denote the directional derivative along the unit length outward pointing tangent vector to the geodesic arcs through the point p. Supposing that $s \in (0, r)$, the numbers $K^2(s)$ and $K^2(r)$ differ by

$$K(r)^2 - K(s)^2 = 2\int_s^r \frac{1}{\rho} N(\rho) K(\rho)^2 \, d\rho$$

(3.17)

this being an application of (3.8) and the fundamental theorem of calculus. It follows from the definition of N in (3.6) that the integrand on the right hand side of (3.17) is no greater than $\rho^{-3} r^2 K^2(r) N(r)$; and therefore

$$K(r)^2 - K(s)^2 \leq \frac{r^2}{s^2} K(r)^2 N(r) .$$

(3.18)

Since $K^2(r) \geq K^2(s)$, this bound leads to the inequality

$$\int_{B_{(1-\mu/8)r} - B_{(1-\mu/4)r}} ||a|^2 - \tfrac{1}{2\pi^2} K(r)^2| \leq \int_{B_{(1-\mu/8)r} - B_{(1-\mu/4)r}} ||a|^2 - \tfrac{1}{2\pi^2} K(\cdot)^2| \leq c_\mu K(r)^2 N(r) .$$

(3.19)

with K viewed on the left hand side of (3.19) as a function on $B_r$–p that depends only on the radial coordinate.

Step 3: Let $\underline{K}^2$ denote for the moment the function on $(0, c_0^{-1}]$ that gives the average of the function $|a|^2$ over any given radius sphere with center p. (This exception to the rule noted at the end of Section 3a of treating $K^2$ and $\underline{K}^2$ is being made in this case to exhibit explicitly how the appropriate $c_0 r^2$ error term appears.) Since $K^2$ differs from $\underline{K}^2$ by at most $c_0 r^2 K^2(r)$, the inequality in (3.19) implies in turn that

$$\int_{B_{(1-\mu/8)r} - B_{(1-\mu/4)r}} ||a|^2 - \tfrac{1}{2\pi^2} K(r)^2| \leq \int_{B_{(1-\mu/8)r} - B_{(1-\mu/4)r}} ||a|^2 - \tfrac{1}{2\pi^2} \underline{K}(\cdot)^2| \leq c_\mu K(r)^2 (r^2 + N(r)) .$$

(3.20)

Meanwhile, the Cauchy-Schwarz inequality bounds the integral on the right hand side of the preceding inequality by

$$c_\mu r^2 \Big( \int_{B_{(1-\mu/8)r} - B_{(1-\mu/4)r}} ||a|^2 - \tfrac{1}{2\pi^2} \underline{K}(\cdot)^2|^2 \Big)^{1/2} ,$$

(3.21)

which is, in turn, bounded by $c_\mu r^3$ times the $L^2$ norm of $\nabla(|a|^2)$. This is so because only the constant functions are in the kernel of the Laplacian on the constant $s \in (0, c_0^{-1}]$ spheres, and because the next smallest eigenvalue is no less than $c_0 s^{-2}$. Now, the $L^2$ norm of $\nabla(|a|^2)$ on $B_{(1-\mu.8)r}$ is no less than the sup-norm of $|a|$ on this sphere times the $L^2$ norm of $|\nabla_A a|$ and thus no less than $c_\mu r K(r)^2 \sqrt{N(r)}$. (The sup norm of $|a|$ on $B_{(1-\mu/8)r}$ is bounded by $c_\mu K(r)$; this is the assertion of the top bullet of Proposition 3.1 with $\mu$ replaced by $\mu/8$.)

Putting all of the preceding together gives a $c_\mu\, \kappa^2(r)\sqrt{N(r)}\, r^4$ bound for (3.21) which implies what is asserted by (3.16).

*Part 3*: This part of the proof explains why the third bullet of Proposition 3.1 is true. The argument that follows for this has three steps. Supposing that $z > 1$ has been specified, these steps use $V_z$ to denote the subset of $B_{(1-\mu)r}$ where $|a| \leq (1-z^{-1})\frac{1}{\sqrt{2\pi}}\kappa(r)$. The fraction of the volume of $B_{(1-\mu)r}$ that is accounted for by $V_z$ is denoted by $\Omega_z$. (The second bullet of Lemma 3.1 concerns the case when $z = 4$.) The volume of $B_{(1-\mu)r}$ differs from $\frac{1}{2}\pi^2(1-\mu)^4 r^4$ by at most $c_0 r^6$.

<u>Step 1</u>: Multiply both sides of (3.18) by $s^3$ and then integrate the resulting inequality on the domain $(0, (1-\mu)r]$. The result, after some rearrangement says that

$$\int_{B_{(1-\mu)r}} |a|^2 \geq (1 - c_0 r^2)\tfrac{1}{4}(1-\mu)^4 r^4 (1 - 2N)\kappa^2 \ .$$

(3.22)

Write $B_r$ as $(B_r - V_z) \cup V_z$. Because the integral of $|a|^2$ on $V_z$ is at most $(1 - z^{-1})(\tfrac{1}{4} r^4\, \Omega)\kappa^2$ when $r \leq c_0^{-1}$, the preceding inequality implies in turn that

$$\int_{B_{(1-\mu)r} - V_z} |a|^2 + (1 - z^{-1})(\tfrac{1}{4}(1-\mu)^4 r^4\, \Omega)\kappa^2 \geq (1 - c_0 r^2)\tfrac{1}{4}(1-\mu)^4 r^4 (1 - N)\kappa^2 \ .$$

(3.23)

Meanwhile, the integral of $|a|^2$ over $B_{(1-\mu)r} - V_z$ is no greater than the maximum of $|a|^2$ on $B_{(1-\mu)r}$ times the volume of $B_{(1-\mu)r}$. Thus, (3.23) and the bound in the second bullet of the proposition lead to the following (after dividing both sides by $\tfrac{1}{4}(1-\mu)^4 r^4 \kappa^2$):

$$(1 + c_\mu \sqrt{N})(1 - \Omega) + (1 - z^{-1})\Omega \geq (1 - c_0 r^2)(1 - N) \ .$$

(3.24)

This inequality can be rearranged to obtain an inequality of the form $\Omega_z < c_\mu\, z(r^2 + \sqrt{N})$ when $r \leq c_0^{-1}$ and when $N \leq 1$.

<u>Step 2</u>: The bound $\Omega_z < c_\mu(r^2 + \sqrt{N})$ from the previous step (when $z = 4$) is weaker than the bound that is asserted by the third bullet of the proposition. This step and Step 3 describe an iterative process that gets the stronger bound. To start, use the bump function $\chi$ to create a non-increasing function on $[0, \infty)$ that equals 1 on $[0, (1 - \tfrac{1}{2z})]$ and equals 0 on $[(1 - \tfrac{1}{4z}), \infty)$. Denote this function by $\beta$. Require that $|d\beta| \leq c_0 z^{-1}$. Use $\beta$ to define the function to be denoted by $f$ on $B_r$ by the rule $f(\cdot) = \beta(\sqrt{2\pi}K(r)^{-1}|a|(\cdot))$. Thus, the function $f$ is 1 where $|a| \leq (1 - \frac{1}{2z})\frac{1}{\sqrt{2\pi}}K(r)$; and it is zero where $|a| \geq (1 - \frac{1}{4z})\frac{1}{\sqrt{2\pi}}K(r)$.

The integral of $f^2$ over $B_{(1-\mu)r}$ obeys

$$\int_{B_{(1-\mu)r}} f^2 \geq (1 - c_0 r^2)\frac{1}{2\pi^2}(1-\mu)^4 r^4 \Omega_z$$

(3.25)

because $f = 1$ on the set $V_z$. Meanwhile, Hölder's inequality bounds the left hand side by

$$c_0 r^2 \sqrt{\Omega_{4z}} \, (\int_{B_{(1-\mu)r}} f^4)^{1/2}$$

(3.26)

because the support of $f$ is in the set $V_{4z}$. The top bullet of (2.1) bounds (3.26) by

$$c_0 r^2 \sqrt{\Omega_{4z}} \int_{B_{(1-\mu)r}} (|df|^2 + r^{-2}f^2);$$

(3.27)

and since $f \leq 1$ and since $|df| \leq c_0 z^{-1} K(r)^{-1}|\nabla_A a|$, this last bound implies in turn that

$$\int_{B_{(1-\mu)r}} f^2 \leq c_0 \sqrt{\Omega_{4z}}\, r^4 (N(r) + \Omega_{4z}).$$

(3.28)

(This bound invoked the formula in (3.7) for N.) Now, an appeal to the previous step (with $z$ replaced there by $4z$) leads from (3.25) and (3.28) to the bound

$$\Omega_z \leq c_\mu z^2 (r^2 + N(r)^{3/4})$$

(3.29)

which is an improvement over the bound in Step 1.

Step 3: Redo Step 2 with $z$ replaced by $4z$ to see that $\Omega_{4z} \leq c_\mu z^2(r^2 + N(r)^{3/4})$. With this bound in hand, return to Step 2 and use this bound to replace the appearances of $\Omega_{4z}$ in (3.28); and then use the resulting inequality with (3.25) to replace (3.29) by

$$\Omega_z \leq c_\mu z^3 (r^2 + N(r)^{9/8}).$$

(3.30)

One more iteration using the $z \to 4z$ version of (3.30) in (3.28) replaces the right hand side of (3.30) by $c_\mu z^3 (r^2 + N(r)^{27/16})$. This gives the bound in the third bullet of Proposition 3.1 because $\Omega_4$ is the fraction of the volume of $B_{(1-\mu)r}$ where $|a| \leq \frac{3}{4}\frac{1}{\sqrt{2\pi}}K(r)$.

### d) The size of $\nabla_A a$ and $r \langle a, \tau a \rangle$ on a ball of radius less than $r_{cF}$

Fix $p \in U$. The next two propositions state a priori bounds that hold on balls centered at $p$ with radii less than $r_{cF}$. These are the analogs for the equations in (1.5) of Propositions 3.1 and 3.2 in [T1]. With a number $r_c \in (0, c_0^{-1}]$ specified, the upcoming propositions use $r_{\ddagger}$ to denote $\kappa(r_c) r$ and they use $a_{\ddagger}$ to denote $\frac{1}{\kappa(r_c)} a$.

**Proposition 3.2**: *There exists $\kappa > 1$; and given $\mu \in (0, \frac{1}{100}]$, there exists $\kappa_\mu > \kappa$; these numbers having the following significance: Fix $r > 1$ and suppose that $(A, a)$ is a connection on $E$ and section of $\mathbb{S}_E^+$ that obeys (1.5). Fix $c \geq \kappa$ and $p \in X$ so as to define $r_{cF}$ using $(A, a)$. Fix positive $r_c$ no greater than $r_{cF}$.*

- $\int_{B_{r_c}} \frac{1}{\text{dist}(q, \cdot)^2} (|\nabla_A a_{\ddagger}|^2 + r_{\ddagger}^2 |\langle a_{\ddagger}, \tau a_{\ddagger} \rangle|^2) \leq \kappa_\mu (c^{-2} + r_c^4 + N(r_c))$ *if* $q \in B_{(1-\mu)r_c}$.

- $\int_{B_{(1-\mu)r_c}} (|\nabla_A (\nabla_A a_{\ddagger})|^2 + r_{\ddagger}^2 |\nabla \langle a_{\ddagger}, \tau a_{\ddagger} \rangle|^2) \leq \kappa_\mu (c^{-2} + r_c^4 + N(r_c))$.

This proposition is proved in Section 3e.

The next proposition talks about the variation across $B_{(1-\mu)r_c}$ of $|a_{\ddagger}|^2$ and $\langle a_{\ddagger}, \tau a_{\ddagger} \rangle$.

**Proposition 3.3**: *There exists $\kappa > 1$, and given $\mu \in (0, \frac{1}{100}]$, there exists $\kappa_\mu > \kappa$; these numbers having the following significance: Fix $r > 1$ and suppose that $(A, a)$ is a connection on $E$ and section of $\mathbb{S}_E^+$ that obeys (1.5). Fix $c \geq \kappa$ and then a point $p \in X$ and use this data with $(A, a)$ to define $r_{cF}$. If $r_c$ is positive but no greater than the minimum of $r_{cF}$ and $c^{-1}$, and if $N(r_c)$ is less than $c^{-2}$, then $||a_{\ddagger}| - \frac{1}{\sqrt{2\pi}}| < \kappa_\mu c^{-1}$ on $B_{(1-\mu)r_c}$.*

The proof of Proposition 3.3 is in Section 3f.

### e) Proof of Proposition 3.2

The first bullet of the proposition follows from the second and the inequality in the second bullet of (2.1). The proof of the second bullet of the proposition differs in only one respect from that of the second bullet of Proposition 3.1 in [T1]. The full argument is provided because of this difference.

The argument for the proposition's second bullet is given below in five steps. These steps exploit the fact that $(A, a_{\ddagger})$ obeys (1.5) with $r$ replaced by $r_{\ddagger}$.

Step 1: Since $(A, a_{\ddagger})$ obeys (1.5) with $r$ replaced by $r_{\ddagger}$, the pair $(A, a_{\ddagger})$ also obeys (2.18) with $r$ replaced by $r_{\ddagger}$. Multiply both sides of the $(r_{\ddagger}, (A, a_{\ddagger}))$ version of (2.18) by $\chi_{\mu}^2$ and then integrate over $B_{r_c}$. An integration by parts then leads to inequality that can be written schematically as

$$\int_{B_{r_c}} (|\nabla_A (\chi_{\mu} \nabla_A a_{\ddagger})|^2 + r_{\ddagger}^2 \chi_{\mu}^2 |\nabla \langle a_{\ddagger}, \tau a_{\ddagger} \rangle|^2) + \alpha + \beta + \gamma \leq c_{\mu} r_c^{-2} \int_{B_{r_c}} (|\nabla_A a_{\ddagger}|^2 + c_0 \int_{B_{r_c}} |a_{\ddagger}|^2$$

(3.31)

with the terms $\alpha$, $\beta$ and $\gamma$ are as follows: The $\alpha$ term is the integral over $B_{r_c}$ of the function $-\chi_{\mu}^2 r_{\ddagger}^2 \langle a_{\ddagger}, \tau^a a_{\ddagger} \rangle \langle \nabla_{A,i} a_{\ddagger}, \tau^a \nabla_{A,i} a_{\ddagger} \rangle$. The term denoted by $\beta$ is the integral over $B_{r_c}$ of $-\chi_{\mu}^2 \langle \nabla_{A,i} a, a \rangle \nabla_j F_{A,ij}$; and the term that is denoted by $\gamma$ is the integral over $B_{r_c}$ of the function $\chi_{\mu}^2 \langle F_{A,ij}, \langle \nabla_{A,i} a, \nabla_{A,j} a \rangle \rangle$.

Step 2: The absolute value of $\alpha$ obeys

$$|\alpha| \leq c_0 c^{-1} \left( \int_{B_{r_c}} |\chi_{\mu} \nabla_A a_{\ddagger}|^4 \right)^{1/2}$$

(3.32)

because the $L^2$ norm of $r^2 |\langle a, \tau a \rangle|$ on $B_{r_c}$ is less than $c^{-1}$ (which is because $r_c < r_{cF}$).

Step 3: To bound the absolute value of $\beta$, first use the Bianchi identity to write $d_A * F_A$ as $2 d_A * F_A^+$. This writes $\beta$ as the integral over $B_{r_c}$ of $-2\chi_{\mu}^2, \langle \nabla_{A,i} a_{\ddagger}, a_{\ddagger} \rangle \nabla_i F_A^+{}_{ij}$. Next, an integration by parts writes $\beta$ as a sum of three terms which are the integrals of

- $2 \nabla_j (\chi_{\mu}^2) F_A^+{}_{ij} \langle \nabla_{A,i} a_{\ddagger}, a_{\ddagger} \rangle$,
- $2\chi_{\mu}^2 F_A^+{}_{ij} \langle \nabla_{A,i} a_{\ddagger}, \nabla_{A,j} a_{\ddagger} \rangle$,
- $-2\chi_{\mu}^2 |F_A^+|^2 |\alpha_{\ddagger}|^2$.

(3.33)

To say something about the size of the integrals of these three expressions, introduce by way of notation $\kappa_{\mu/16}$ to denote the version $\mu/16$ version of what is denoted by $\kappa_{\mu}$ in the top bullet of Proposition 3.1. Keep in mind that this number $\kappa_{\mu/16}$ is greater than the norm of $|a_{\ddagger}|$ on the support of $\chi_{\mu}$.

The norm of the integral of the term in the top bullet of (3.33) is bounded by $\kappa_{\mu/16} \mu^{-1} c^{-1} N(r_c)^{1/2}$ because the $L^2$ norm of $\nabla_A a_{\ddagger}$ on the support of $\chi_{\mu}$ is bounded by $r_c N(r_c)^{1/2}$ and the $L^2$ norm of $F_A^+$ over the support of $\chi_{\mu}$ is no greater than $c^{-1}$. The norm of the integral of the term in the middle bullet of (3.33) is bounded by the expression on the right hand side of (3.32) because of the $c^{-1}$ bound for the $L^2$ norm of $F_A^+$ on the support of the function $\chi_\mu^2$. The norm of the integral of the term in the third bullet of (3.33) is bounded by $\kappa_{\mu/16}^2 c^{-2}$ because of the $\kappa_{\mu/16}$ bound for $|a_\ddagger|$ on the support of $\chi_\mu$ and this same $c^{-1}$ bound for the $L^2$ norm of $F_A^+$ on the support of $\chi_\mu$.

   Step 4: An integration by parts writes the term $\gamma$ as the sum of three terms which are the integrals over $B_{r_c}$ of the following:

- $2\nabla_j(\chi_\mu^2) F_{A,ij}\langle \nabla_{A,i} a_\ddagger, a_\ddagger \rangle$,
- $-2\chi_\mu^2 |F_A|^2 |a_\ddagger|^2$,
- $2\chi_\mu^2 (\nabla_{A,i} F_A)_{ij} \langle a_\ddagger, \nabla_{A,j} a_\ddagger \rangle$.

(3.34)

The integral of the term in the top bullet of (3.34) is no greater than $c_\mu \mu^{-1} c^{-1} N(r_c)^{1/2}$. This is because the $L^2$ norm of $\nabla_A a_\ddagger$ on the support of $\chi_\mu$ is bounded by $r_c N(r_c)^{1/2}$ and the $L^2$ norm of $|F_A||a_\ddagger|$ on the support of $\chi_\mu$ is no greater than $c_\mu c^{-1}$ since $r_c \leq r_{cF}$ (and since $|a_\ddagger|$ has an apriori bound by $c_\mu$). This is also why the norm of the integral of the term in the second bullet of (3.34) is less than $c_\mu c^{-2}$. The integral of the term in the third bullet of (3.34) is $-2\alpha$.

   Step 6: Use the bounds from the preceding steps in (3.31) to conclude that

$$\int_{B_{r_c}} (|\nabla_A(\chi_\mu \nabla_A a_\ddagger)|^2 + r_\ddagger^2 \chi_\mu^2 |\nabla\langle a_\ddagger, \tau a_\ddagger\rangle|^2) \leq c_0 c^{-1} (\int_{B_{r_c}} |\chi_\mu \nabla_A a|^4)^{1/2} + c_\mu(c^{-2} + r_c^4) + c_\mu N(r_c).$$

(3.35)

To finish the proof, use the fact that the norm of $d|\chi_\mu \nabla_A a|$ is at most $|\nabla_A(\chi_\mu \nabla_A a)|$; and then use the Sobolev inequalities from the first and third bullets in (2.1) to bound the square of the $L^4$ norm of $|\chi_\mu \nabla_A a|$ that appears on the right hand side (3.35) by $c_0$ times the left hand side of (3.35). Insert the latter bound to obtain the bound that is claimed by the second bullet of Proposition 3.2 bound when $c > c_0$.

### f) Proof of Proposition 3.3

   Before starting, take note that there is much borrowing of arguments for the analogous Proposition 3.2 in [T1]. This borrowing is often not noted.
   The proof invokes the conclusions of the Propositions 3.1 and 3.2 with $\mu$ replaced by $\mu/64$. These propositions can be used in the event that $c > c_\mu$ and $N(r_c) < 1$; and these conditions are assumed in what follows.
   To start the proof, invoke the second bullet of the version of Proposition 3.1 with $\mu$ replaced by $\mu/64$ to obtain half of what needs to be proved since Proposition 3.1 finds

$$|a_{\ddagger}| - 1 \leq c_\mu \sqrt{N(r_c)}$$

(3.36)

on the radius $(1-\mu/64)\,r_c$ ball centered at p. This understood, it remains to prove that

$$|a_{\ddagger}| - 1 \geq -c_\mu (\sqrt{N(r_c)} + r_c^2 + c^{-2}) \ .$$

(3.37)

To prove (3.37), return to the identity in (3.15). As noted in Step 3 of Part 1 of the proof of Proposition 3.1, the absolute value of the left most integral on the right hand side of (3.15) (this is integral of $\chi_{\mu/4} G_q \langle a, \Re a \rangle$) is at most $c_\mu r_c^2 K(r_c)^2$. Meanwhile, the inequality in (3.16) leads to a $c_0 \sqrt{N(r_c)}\, K(r_c)^2$ bound for the absolute value of the left most integral in (3.15). As a consequence of these bounds for the absolute value of the integrals on right hand side of (3.15), the latter identity implies the following inequality:

$$|a|^2(q) - \tfrac{1}{2\pi^2} K(r)^2 \geq -2 \int_{B_r} \chi_{\mu/4} G_q (|\nabla_A a|^2 + r^2 |\langle a, \tau a \rangle|^2) - c_\mu (\sqrt{N(r_c)} + r_c^2) K(r_c)^2 \ .$$

(3.38)

Meanwhile, the integral that appears on the right hand side of (3.38) is no greater than $c_\mu c^{-2} K(r_c)^2$, this being a consequence of the first bullet of Proposition 3.2. Therefore, (3.28) implies the inequality

$$|a|^2(q) - \tfrac{1}{2\pi^2} K(r)^2 \geq - c_\mu (\sqrt{N(r_c)} + r_c^2 + c^{-2}) K(r_c)^2 \ .$$

(3.39)

This with (3.36) lead directly to (3.37).

### 4. The size of $F_A$ on small radius balls

This section states and proves a proposition about the size of $F_A$ on balls of radius less than the minimum of $c^{-2}$ and $r_{cF}$. This proposition say in effect that $F_A$ is unexpectedly small when $N$ is small. The proposition is an analog for the equations in (1.5) of Proposition 4.1 in Section 4 of [T1]. It is also a direct analog of a central result in [HW] which concerns the dimension 3 version of (1.5).

**Proposition 4.1**: *Fix $\mu \in (0, \tfrac{1}{100}]$ and there exists $\kappa_\mu > 1$ such that the subsequent assertion is true. Fix $r > 1$ and assume that $(A, a)$ is a pair of connection on $E$ and section of $\mathbb{S}^+ \otimes E$ that obeys (1.5). Fix $p \in X$ and $c > \kappa_\mu$ to define $r_{cF}$ using $(A, a)$. Given $\varepsilon \in (0, 1)$, suppose that $r_c$ is positive but less than the smaller of $c^{-1}$, $r_{cF}$ and $\kappa_\mu^{-1} c^{-1} \varepsilon^2$ and that $N(r_c)$ is less than $\kappa_\mu^{-1} c^{-2} \varepsilon^3$. Then $\int_{B_{(1-\mu) r_c}} |F_A|^2 < \varepsilon c^{-2}$.*

The remainder of this section is occupied with the proof of this proposition. The arguments assume (sometimes implicitly) that $r_c$ and $N(r_c)$ are such that the conclusions of the versions of Propositions 3.1-3.3 can be invoked when $\mu$ is replaced by $\mu/64$. The arguments that follow also use $r_{\ddagger}$ as shorthand for $r\kappa(r_c)$ and $a_{\ddagger}$ as shorthand for $\frac{1}{\kappa(r_c)} a$.

### a) Proof of Proposition 4.1 when $r\, r_c\, \kappa(r_c)$ is large

This subsection proves Proposition 4.1 with the additional assumption that $r_{\ddagger} r_c$ (which is $r r_c \kappa(r_c)$) is large (which means in this case that it is greater than $c_\mu^{-1} \varepsilon^{-1}$). To start, take the inner product of both sides of (2.16) with $\chi_{\mu/8}^2 F_A$ and integrate the resulting equation over $B_{r_c}$. An integration by parts leads to the inequality

$$\int_{B_{r_c}} |\nabla(\chi_{\mu/8} F_A)|^2 + r_{\ddagger}^2 \int_{B_{r_c}} \chi_{\mu/8}^2 |a_{\ddagger}|^2 |F_A|^2 \le c_\mu \frac{1}{r_c^2} \int_{B_{(1-\mu/16) r_c}} |F_A|^2 + r_{\ddagger}^2 \int_{B_{r_c}} \chi_{\mu/8}^2 |\nabla_A a_{\ddagger}|^2 |F_A| + c_0 r_c^6 .$$

(4.1)

Hölder's inequality bounds the right most integral in (4.1) by a product $\alpha \beta \gamma$ with $\alpha$, $\beta$ and $\gamma$ being positive and being defined by the following rules:

$$\alpha^2 = \int_{B_{r_c}} |\nabla_A a_{\ddagger}|^2 , \quad \beta^4 = \int_{B_{(1-\mu/16) r_c}} |\nabla_A a_{\ddagger}|^4 , \quad \gamma^4 = \int_{B_{r_c}} \chi_{\mu/8}^2 |F_A|^4 .$$

(4.2)

What is denoted above by $\alpha^2$ is no greater than $r_c^2 N(r_c)$ as can be seen from (3.7). Meanwhile, it follows from the second bullet of Proposition 3.2 (with $\mu$ replaced by $\mu/64$) and the top bullet in (2.1) that the number $\beta^4$ is no greater than $c_\mu (r_c^2 + N(r_c)^2 + c^{-2})^2$. As for $\gamma^4$, it follows from (2.1) that $\gamma^2 \le c_\mu \int_{B_{r_c}} |\nabla(\chi_{\mu/8} F_A^+)|^2$. With these bounds, (4.2) gives

$$\int_{B_{r_c}} |\nabla(\chi_{\mu/8} F_A)|^2 + r_{\ddagger}^2 \int_{B_{r_c}} \chi_{\mu/8}^2 |a_{\ddagger}|^2 |F_A|^2 \le c_\mu \frac{1}{r_c^2} \int_{B_{(1-\mu/16) r_c}} |F_A|^2 + c_\mu r_{\ddagger}^2 N(r_c)(r_c^2 + N(r_c)^2 + c^{-2}) + c_0 r_c^6 .$$

(4.3)

To proceed from here, note that the curvature integral on the right hand side of (4.3) is no greater than $c^{-2}$ because $r_c$ is assumed to be less than $r_{cF}$. Meanwhile, the version of Proposition 3.3 with $\mu$ replaced by $\mu/64$ can be brought to bear to see that the integral of $\chi_{\mu/8}^2 |a_{\ddagger}|^2 |F_A|^2$ that appears on the left hand side of (4.3) is no smaller than $c_\mu^{-1}$ times the integral of $\chi_{\mu/8}^2 |F_A|^2$. (This is assuming that $N(r_c) < 1$ and that $r_c$ is less than the minimum of $c^{-1}$ and $r_{cF}$). With the preceding bounds understood, then (4.3) implies the following:

$$\frac{1}{r_\ddagger^2} \int_{B_{(1-\mu)r_c}} |\nabla F_A|^2 + \int_{B_{(1-\mu)r_c}} |F_A|^2 \leq c_\mu (\frac{1}{r_\ddagger^2 r_c^2} c^{-2} + N(r_c)(r_c^2 + N(r_c)^2 + c^{-2})).$$

(4.4)

This last inequality leads directly to the bound asserted by Proposition 4.1 if $r_\ddagger^2 r_c^2 \geq c_\mu \varepsilon^{-1}$ and if $N(r_c) < c_\mu^{-1} c^{-2} \varepsilon$.

**b) Proof of Proposition 4.1 when $r\, r_c\, \kappa(r_c)$ is small**

This subsection gives the proof of Proposition 4.1 when $r_\ddagger r_c$ is less than $c_\mu \varepsilon^{-1/2}$. These arguments are given in five parts.

*Part 1*: What with (3.7), it follows that

$$\int_{B_{(1-\mu/16)r_c}} |F_A^+|^2 \leq r_\ddagger^2 r_c^2 N(r_c).$$

(4.5)

Therefore, if $z > 1$ and $\varepsilon > 0$ are given, and supposing that $r_\ddagger r_c \leq z \varepsilon^{-1}$ and $N(r_c) < \varepsilon^3 z^{-2} c^{-2}$, then right hand side of (4.5) is no greater than $\varepsilon c^{-2}$.

*Part 2*: Any chosen isomorphism between the bundle E on $B_{(1-\mu/16)r_c}$ and the product bundle $B_{(1-\mu/16)r_c} \times \mathbb{C}$ will write the connection A as $\theta_0 + A$ with $\theta_0$ being the product connection and with A being an $i\mathbb{R}$-valued 1-form. A particular isomorphism is chosen momentarily.

Writing A as $\theta_0 + A$ using an isomorphism from E on $B_{(1-\mu/16)r_c}$ to the product bundle writes $F_A$ as $dA$ and $F_A \wedge F_A$ as $dA \wedge dA$. Since $*(F_A \wedge F_A) = |F_A^-|^2 - |F_A^+|^2$, the writing of $F_A$ as $dA \wedge dA$ using a chosen isomorphism can be used (with an integration by parts) to see that

$$\int_{B_{(1-\mu/16)r_c}} \chi_{\mu/2} |F_A^-|^2 = \int_{B_{(1-\mu/16)r_c}} \chi_{\mu/2} |F_A^+|^2 - \int_{B_{(1-\mu/16)r_c}} d\chi_{\mu/2} \wedge A \wedge F_A.$$

(4.6)

Therefore, by virtue of (4.5),

$$\int_{B_{(1-\mu/16)r_c}} \chi_{\mu/2} |F_A^-|^2 \leq r_\ddagger^2 r_c^2 N(r_c) + c_\mu r_c^{-1} (\int_{B_{(1-\mu/4)r_c}} |A|^2)^{1/2} (\int_{B_{(1-\mu/4)r_c}} |F_A|^2)^{1/2}.$$

(4.7)

Since the integral of $|F_A|^2$ over the ball $B_{r_c}$ is less than $c^{-2}$ by assumption, the inequality in (4.7) implies that

$$\int_{B_{(1-\mu/16)r_c}} \chi_{\mu/2} |F_A^-|^2 \leq r_\ddagger^2 r_c^2 N(r_c) + c_\mu r_c^{-1} c^{-1} \left( \int_{B_{(1-\mu/4)r_c}} |A|^2 \right)^{1/2}.$$

(4.8)

Save this bound for Part 4 of the proof. Part 3 of the proof describes an isomorphism of E with the product bundle over $B_{(1-\mu/16)r_c}$ whose corresponding $\mathbb{A}$ has small $L^2$ norm.

*Part 3*: Fix a unit length element in $\mathbb{S}^+|_p$ and use parallel transport (with the connection that defines $\nabla$) along the geodesic arcs through p to define a unit length section of $\mathbb{S}^+$ on the whole of $B_{1/c_0}$. The chosen unit length element in $\mathbb{S}^+|_p$ is denoted in what follows by $e$ as is the corresponding section of $\mathbb{S}^+$ on $B_{1/c_0}$. The Hermitian pairing $\langle e, \cdot \rangle$ defines a homomorphism from $\mathbb{S}^+ \otimes E$ to E over $B_{1/c_0}$. Thus $\langle e, a_\ddagger \rangle$ defines a section of E over $B_{(1-\mu/16)r_c}$. The next lemma talks about this section.

**Lemma 4.2**: *Fix $\mu \in (0, \frac{1}{100}]$ and there exists $\kappa_\mu > 1$ with the following significance: Fix $r > 1$ and assume that $(A, a)$ is a pair of connection on E and section of $\mathbb{S}^+ \otimes E$ that obeys (1.5). Fix $p \in X$ and fix $c > \kappa_\mu^{-1}$ to define $r_{cF}$ using $(A, a)$. Suppose that $r_c$ is positive but less than the minimum of $c^{-1}$, $r_{cF}$ and $\kappa_\mu^{-1}$. Suppose also that $N(r_c)$ is less than the minimum of $\kappa_\mu^{-1}$ and $c^{-2}$. Granted these constraints, then there is a unit normed section, e, of $\mathbb{S}^+$ over $B_{1/c_0}$ as just described such that $|\langle e, a_\ddagger \rangle| \geq \kappa_\mu^{-1}$ on $B_{(1-\mu/16)r_c}$.*

This lemma is proved momentarily. Accept it as true for the time being.

Let $e$ denote the section from Lemma 4.2. Since $\langle e, a_\ddagger \rangle$ is nowhere zero on $B_{(1-\mu/16)r_c}$, there is an isomorphism between E and the product complex line bundle on this ball that identifies the section $\langle e, a_\ddagger \rangle$ with the positive real number $|\langle e, a_\ddagger \rangle|$. Use this isomorphism to write A as $\theta_0 + \mathbb{A}$ with $\mathbb{A}$ being an $i\mathbb{R}$-valued 1 form. The identification of E with the product bundle identifies $\nabla_A \langle e, a_\ddagger \rangle$ with $d|\langle e, a_\ddagger \rangle| + \mathbb{A} |\langle e, a_\ddagger \rangle|$. Since the first term in this sum is an $\mathbb{R}$-valued 1-form and the second term is an $i\mathbb{R}$-valued 1-form, the norm of $|\nabla_A \langle e, a_\ddagger \rangle|$ obeys

$$|\nabla_A \langle e, a_\ddagger \rangle| \geq |\langle e, a_\ddagger \rangle| |\mathbb{A}|.$$

(4.9)

Therefore, because $|\nabla_A \langle e, a_\ddagger \rangle| \leq c_0 (|a_\ddagger| + |\nabla_A a_\ddagger|)$ and because $|\langle e, a_\ddagger \rangle| \geq \kappa_\mu^{-1}$, the inequality in (4.9) leads to the $L^2$ bound

$$\int_{B_{(1-\mu/16)r_c}} |A|^2 \le c_\mu r_c^2 (N(r_c) + r_c^2) \ .$$

(4.10)

Use the bound from (4.10) in (4.8) to obtain an $L^2$ bound for $F_A^-$:

$$\int_{B_{(1-\mu)r_c}} |F_A^-|^2 \le r_\ddagger^2 r_c^2 N(r_c) + c_\mu c^{-1}(r_c + \sqrt{N(r_c)}) \ .$$

(4.11)

It follows as a consequence of this bound that if $r_\ddagger r_c \le z^2 \varepsilon$, then the right hand side of (4.11) will be smaller that $\varepsilon c^{-2}$ if $r_c < c_\mu^{-1} c^{-1} \varepsilon$ and if $N(r_c)$ is less than $c_\mu^{-1} z^{-2} c^{-2} \varepsilon^3$. After this last observation, the proof of Proposition 4.1 requires only a proof of Lemma 4.2

*Part 4*: The upcoming proof of Lemma 4.2 requires the bound that follows directly for the $L^2$ norm of $\nabla F_A^+$:

$$\int_{B_{(1-\mu/4)r_c}} |\nabla F_A^+|^2 \le c_\mu (r_\ddagger^2 r_c^2) \frac{1}{r_c^2} N(r_c) + c_0 r_c^6 \ .$$

(4.12)

The derivation starts with (2.16) and its inner product with $\chi_{\mu/8}^2 F_A^+$. Integrating the result leads to a modified version of (4.1) with $F_A$ replaced in each occurrence with $F_A^+$. The latter version of (4.1) leads to a modified version of (4.3) with, again, $F_A$ replaced in all occurrences with $F_A^+$. This modfied version of (4.3) leads directly to the following:

$$\int_{B_{(1-\mu/4)r_c}} |\nabla F_A^+|^2 \le c_\mu \frac{1}{r_c^2} \int_{B_{(1-\mu/16)r_c}} |F_A^+|^2 + c_\mu r_\ddagger^2 N(r_c)(r_c^2 + N(r_c)^2 + c^{-2}) + c_0 r_c^6.$$

(4.13)

Since the integral of $F_A^+$ that appears in (4.13) is no greater than $r_\ddagger^2 r_c^2 N(r_c)$, the inequality in (4.13) implies the one in (4.12).

*Part 5*: This part of the proof has the task of proving Lemma 4.2.

**Proof of Lemma 4.2**: The equation in (2.9) holds with $r$ replaced by $r_\ddagger$ and with $a$ replaced by $a_\ddagger$. This understood, take the inner product of both sides of the $(r_\ddagger, (A, a_\ddagger))$ version of (2.9) with $e$ to obtain a differential equation for the section $\langle e, a_\ddagger \rangle$ of E on $B_{1/c_0}$ that has the schematic form

$$\nabla_A{}^\dagger \nabla_A \langle e, a_\ddagger \rangle + r_\ddagger{}^2 \langle a_\ddagger, \tau^a a_\ddagger \rangle \langle e, \tau^a a_\ddagger \rangle + \mathfrak{Q}_0 a + \mathfrak{Q}_1 \nabla_A a = 0$$
(4.14)

with $\mathfrak{Q}_0$ and $\mathfrak{Q}_1$ being homomorphism on $B_{1/c_0}$ from $\mathbb{S}^+ \otimes E$ to E and from $(\mathbb{S}^+ \otimes E) \otimes T^*X$ to E with norms bounded by $c_0$. Let $w$ denote $|\langle e, a_\ddagger \rangle|^2$. Taking the inner product of both sides of (4.14) with $\langle e, a_\ddagger \rangle$ leads to a differential inequality for the function $w$:

$$d^\dagger d w \geq -c_\mu (1 + |\nabla_A a_\ddagger|^2 + |F_A{}^+|) .$$
(4.15)

(The derivation of this inequality uses the a priori bound $|a_\ddagger| \leq c_\mu$ from the $\mu/64$ version of Proposition 3.1. It also uses the formula in (2.5).)

To make something of (4.15), introduce $\underline{w}$ to denote the average of the function $w$ on the ball $B_{(1-\mu/16)r_c}$. Since $\underline{w}$ is constant, the inequality in (4.15) also holds if $w$ is replaced by $w - \underline{w}$. Keeping this in mind, let q denote a point in $B_{(1-\mu)r_c}$ and let $G_q$ again denote the Dirichlet Green's function for the operator $d^\dagger d$ on $B_{r_c}$ with pole at q. Multiply both sices of (4.15) by $\chi_{\mu/2} G_q$ and then integrate the resulting inequality over $B_{r_c}$. An integration by parts then leads to the inequality

$$w - \underline{w} \geq -c_\mu r_c{}^2 - c_\mu \int_{B_{r_c}} \chi_{\mu/2} G_q |\nabla_A a|^2 - c_\mu \int_{B_{r_c}} \chi_{\mu/2} G_q |F_A{}^+| - c_\mu r_c{}^{-4} \int_{B_{r_c}} \chi_{\mu/4} |w - \underline{w}|$$
(4.16)

The next paragraph supplies suitable bounds for the three integrals that appear on the right hand side of this inequality.

The bound in (3.12) and the bound in the top bullet of Proposition 3.2 can be brought to bear on the left most integral in (4.16) to bound it by $c_\mu(c^{-2} + r_c{}^2 + N(r_c))$. To bound the middle integral on the right hand side of (4.16), first use the Cauchy-Schwartz inequality with (3.12) to bound it by

$$c_\mu r_c \Big( \int_{B_{r_c}} \chi_{\mu/2} \frac{1}{\operatorname{dist}(q,\cdot)^2} |F_A{}^+|^2 \Big)^{1/2} .$$
(4.17)

Use the top bullet in (2.1) with (4.12) to bound this expression by $c_\mu((r_\ddagger r_c) \sqrt{N(r_c)} + r_c{}^4)$. To bound the right most integral in (4.16), first use the Cauchy-Schwartz inequality to bound it by

$$c_\mu r_c{}^2 \Big( \int_{B_{(1-\mu/16)r_c}} |w - \underline{w}|^2 \Big)^{1/2}$$
(4.18)

Because $\underline{w}$ is the average of $w$ over $B_{(1-\mu/16)r_c}$, the integral that appears in (4.18) is no larger than the integral over $B_{(1-\mu/16)r_c}$ of $c_\mu r_c^2 |dw|^2$. Now, the integral of $|dw|^2$ is no larger than $c_0$ times the sum of $r_c^4$ and the integral over $B_{(1-\mu/16)r_c}$ of $|\nabla_A a_\ddagger|^2$, and the latter integral is no larger than $c_\mu r_c^2 N(r_c)$. Thus, (4.18) is at most $c_\mu r_c^4(r_c+\sqrt{N(r_c)})$; and the right most term in (4.16) is at most $c_\mu(r_c+\sqrt{N(r_c)})$.

The bound in (4.16) with the bounds in the preceding paragraph give directly that

$$w \geq \underline{w} - c_\mu (c^{-2} + r_c + (1 + r_\ddagger r_c) \sqrt{N(r_c)}).$$

(4.19)

Now, it follows from Proposition 3.3 that there is *some* unit section $e$ for which $\underline{w} \geq c_0^{-1}$. This fact with (4.19) lead directly to what is said by Lemma 4.2.

## 5. The size of $F_A$ on large radius balls

Integral bounds for $F_A$ on balls of relatively large radius (but still less than $c_0^{-1}$) are also needed in what follows. Note in this regard that if $p \in X$ and $r \in (0, c_0^{-1}]$ are given, then the $L^2$ norm of $F_A^+$ on $B_r$ is bounded a priori by $r r K \sqrt{N}$, this being a consequence of the formula for $N$ in (3.7). The proposition that follows talks about the $L^2$ norm of the full curvature tensor $F_A$ on balls such as $B_r$.

**Proposition 5.1**: *There exists $\kappa > 100$ with the following significance: Fix $r > 1$ and suppose that $(A, a)$ is a solution to (1.5). Fix a point $p \in X$.*

- *If $r \in (0, \kappa^{-1})$, then $\int_{B_r} |F_A|^2 \leq \kappa\, r^2 r^2$.*

*Let M denote a increasing function on $(0, \kappa^{-1})$ with $M(s) \geq \kappa \int_{B_s} |F_A^+|^2$ for all $s \in (0, \kappa^{-1})$.*

- *Suppose that $r_1 \in (0, \kappa^{-1})$ and $r_0 \in (0, \frac{1}{4} r_1)$ obey $\int_{B_s} |F_A|^2 \geq M(s)$ for all points $s \in [r_0, r_1]$ with $\int_{B_{r_1}} |F_A|^2 = M(r_1)$. Then $\int_{B_{r_0}} |F_A|^2 \leq 64 (\frac{r_0}{r_1})^{2+1/\kappa} M(r_1)$.*

- *If $r_1 = \kappa^{-1}$ and $r_0 \in (0, \frac{1}{4} r_\cdot)$ are such that $\int_{B_s} |F_A|^2 \geq M(s)$ for all $s \in [r_0, r_\cdot]$, then*

$$\int_{B_{r_0}} |F_A|^2 \leq \kappa^2 (\frac{r_0}{r_1})^{1/\kappa} r_0^2\, r^2.$$

By way of a look ahead, the applications of the second and third bullets will be using for the most part functions $M$ that are constant multiples of $r^2 r^2 \kappa^2 N$. It follows from (3.7) that such a function is increasing and it has the required lower bound if the multiplying constant is large.

The proof of this proposition is deferred to Section 5c. The intervening Section 5a states and proves a proposition about the $L^1$ norm of $F_A$. This proposition gives $L^1$ bounds for $|F_A^+|$ on $B_r$ that are better than the Hölder bound (which is $c_0 r^2$ times the $L^2$ norm of $F_A^+$ on $B_r$) when $rr\kappa$ is large. (This proposition is invoked to prove Proposition 5.1.) Meanwhile, Section 5b states and proves a proposition that writes the connection A as a sum of a connection with anti-self dual curvature and an $i\mathbb{R}$ valued 1-form with controlled $L_1^2$ Sobolev norm.

A final subsection (this being Section 5d) states and then proves an assertion to the effect that the $L^2$ norms of $F_A$ on $B_r$ and the numbers $\kappa(r)$ and $N(r)$ determine an a priori upper bound for the $L_1^2$ norm of $F_A$ on $B_{(1-\mu)r}$.

### a) The $L^1$ norm of $F_A^+$

The desired $L^1$ norm bound for the self-dual part of the curvature is stated below in Proposition 5.2.

**Proposition 5.2**: *There exists $\kappa > 2$; and given $\mu \in (0, \frac{1}{100}]$, there exists $\kappa_\mu > \kappa$; and they have the following significance: Fix $r > 1$ and suppose that $(A,a)$ is a solution to (1.5). Fix $p \in X$ and $r \in (0, \kappa^{-1}]$ to define the ball $B_{(1-\mu)r}$. Then*

$$\int_{B_{(1-\mu)r}} |F_A^+| \le \kappa_\mu \, (rr\kappa)^{1-1/\kappa} (1+N)(\sqrt{N}+r)\, r^2$$

The remainder of this subsection is occupied with the proof.

***Proof of Proposition 5.2***: Since the $L^1$ norm of $|F_A^+|$ is bounded via the Cauchy-Schwarz inequality by $c_0 r^2$ times the $L^2$ norm of $F_A^+$, and since the latter is bounded by $rr\kappa\sqrt{N}$ (by virtue of (3.7)), the proposition follows directly if $rr\kappa \le 100$. Therefore, assume in what follows that $rr\kappa \ge 100$. The proof in this case has seven steps.

<u>Step 1</u>: Supposing that n is a positive integer, let $U_n$ denote the subset of $B_{(1-\mu)r}$ where the function $|a|$ obeys the bound $|a| \le 2^{-n} \frac{1}{\sqrt{2\pi}} \kappa(r)$. This step and Steps 2 and 3 prove the following claim about the contribution from $U_n$ to the $L^1$ norm of $F_A^+$:

*There exists $\gamma \in (\frac{1}{2}, 1)$ that is independent of the integer n, the radius r (if $r < c_0^{-1}$), and the chosen point p in X; and which has the following significance: If $r > 1$ and if $(A,a)$ is a solution to (1.5), then*

$$\int_{U_n} |F_A^+| \leq c_\mu \gamma^n (rrK)(1+\sqrt{N})(\sqrt{N}+r) r^2 .$$

(5.1)

To prove this claim, first introduce by way of notation $\Omega$ to denote the fraction of the volume of the ball $B_{(1-\mu)t}$ where $|a| \leq \frac{1}{2\sqrt{2}\pi} K$. This is the fraction of the volume of the ball that is accounted for by $U_1$. The proof starts with the Cauchy-Schwarz inequality $\int_{U_n} |F_A^+| \leq c_0 r^2 \Omega^{1/2} (\int_{U_n} |F_A^+|^2)^{1/2}$; and this leads (via the top bullet of (2.2)) to

$$\int_{U_n} |F_A^+| \leq c_0 r^2 \Omega^{1/2} (\int_{U_n} (|\nabla_A a|^2 + r^2 |\langle a, \tau a\rangle|^2))^{1/2} .$$

(5.2)

Since the third bullet of Proposition 3.1 bounds $\Omega$ by $c_\mu (N^{3/2} + r^2)$, the inequality in (5.1) follows with a proof that

$$\int_{U_n} (|\nabla_A a|^2 + r^2 |\langle a, \tau a\rangle|^2) \leq c_0 \gamma^{2n} r^2 K^2 (1+N)$$

(5.3)

with $\gamma \in (\frac{1}{2}, 1)$ being independent of n, r (if $r \leq c_0^{-1}$), the point $p \in X$ and (A, a) and r. Steps 2 and 3 prove (5.3).

Step 2: To prove (5.3), introduce by way of notation $f_n$ to denote the integral that appears on the left hand side of the inequality. As explained in the next step, the collection $\{f_n\}_{n=1,2,\ldots}$ obey a recursion inequality

$$f_n \leq \frac{x}{1+x} f_{n-1} + c_\mu 2^{-2n} r^2 K^2$$

(5.4)

with $x \in (2, c_0)$. This inequality can be iterated to see that

$$f_n \leq (\frac{x}{1+x})^n f_1 + c_\mu (\sum_{1 \leq m \leq n} (\frac{x}{1+x})^{n-m} 2^{-2m}) r^2 K^2 .$$

(5.5)

Noting that $f_1 \leq r^2 K^2 N$ and that $\frac{x}{1+x} > \frac{1}{2}$, this inequality leads to (5.3) with $\gamma^2 = \frac{x}{1+x}$.

Step 3: To prove (5.4), define the function $\beta_n$ by the rule

$$\beta_n = \chi(2^n \sqrt{2}\pi K^{-1}|a| - 1).$$

(5.6)

This is a non-negative function that is equal to 1 on $U_n$ and equal to 0 on the complement of $U_{n+1}$. It obeys $|d\beta_n| \leq c_0 2^n K^{-1} |\nabla_A a|$. Multiply both sides of (2.9) by $\chi_\mu \beta_n$ and integrate both sides over the ball $B_r$. Integration by parts then leads to the inequality

$$\int_{U_n} (|\nabla_A a|^2 + r^2 |\langle a, \tau a \rangle|^2) \leq c_0 \int_{U_{n-1}-U_n} |\nabla_A a|^2 + c_\mu r^{-2} \int_{U_{n-1}} |a|^2 .$$

(5.7)

Noting that the integral of $|\nabla_A a|^2$ on $U_{n-1} - U_n$ is no greater than $f_{n-1} - f_n$, and noting that $|a|$ on $U_{n-1}$ is no greater than $c_0 2^{-n} K$, this last inequality implies (5.4) if $x$ is no less than the version of $c_0$ that appears in (5.7).

Step 4: With (5.1) now understood, this step and Step 4 derive the bound for the $L^1$ norm of $F_A^+$ where $|a| \geq 2^{-n} \frac{1}{\sqrt{2\pi}} K(r)$ that is stated below in (5.8). By way of notation, this bound involves to the number $\gamma$ that appears in (5.1).

*Suppose that $r > 1$ and that $(A, a)$ is a solution to (1.5). Then*

$$\int_{B_{(1-\mu)r} - U_n} |F_A^+| \leq c_\mu (1+N)(2^{2n} + \gamma^n 2^{3n} \tfrac{1}{rrK})(\sqrt{N}+r) r^2$$

(5.8)

The derivation of this bound starts with the self-dual part of (2.16). Taking the inner product of the latter with $F_A^+$ leads to a differential inequality for $|F_A^+|$ that has the form

$$d^\dagger d |F_A^+| + c_0^{-1} r^2 |a|^2 |F_A^+| \leq c_0 (1 + r^2 |\nabla_A a|^2) .$$

(5.9)

Since $|a| \geq 2^{-n} \frac{1}{\sqrt{2\pi}} |a|^2$ on $B_r - U_{n+1}$ (and hence on $B_r - U_n$ which is inside $B_r - U_{n+1}$), the inequality in (5.9) implies that

$$d^\dagger d |F_A^+| + c_0^{-1} 2^{-2n} r^2 K^2 |F_A^+| \leq c_0 (1 + r^2 |\nabla_A a|^2 + |F_A^+|^2)$$

(5.10)

on the domain $B_r - U_{n+1}$.

Step 5: To exploit (5.10), define the function $\tilde{\beta}_n$ to be $\chi(2 - 2^{n+1} \frac{\sqrt{2\pi}|a|}{K})$. This function is equal to 1 where $|a| \geq 2^{-n} \frac{1}{\sqrt{2\pi}} K$, and it is equal to zero where $|a| \leq 2^{-n-1} \frac{1}{\sqrt{2\pi}} K$ (and thus on $U_{n+1}$). Note in particular that (5.10) holds on the support of $\tilde{\beta}_n$. Multiply both sides of (5.10) by $\chi_{\mu/4} \tilde{\beta}_n$ and integrate the resulting inequality over $B_r$. The result of doing this integration can be written (after some integration by parts and keeping in mind the identity $F_A^+ = \tfrac{1}{2} r^2 [a; a]$) as

$$\tfrac{1}{\pi^2} 2^{-2n-4} r^2 K^2 \int_{B_{(1-\mu)r} - U_n} |F_A^+| \le c_0 r^2 \int_{B_r} (|\nabla_A a|^2 + r^2 |\langle a, \tau a \rangle|^2) + \mathcal{T}_1 + \mathcal{T}_2 + \mathcal{T}_3$$

(5.11)

where the terms $\mathcal{T}_1$, $\mathcal{T}_2$ and $\mathcal{T}_3$ are written below.

- $\mathcal{T}_1 = -\int_{B_r} \chi_{\mu/4} (d^\dagger d\chi_{\mu/2}) \tilde{\beta}_n |F_A^+|$.
- $\mathcal{T}_2 = \int_{B_r} \chi_{\mu/4} \langle \nabla \tilde{\beta}_n, \nabla \chi_{\mu/4} \rangle |F_A^+|$.
- $\mathcal{T}_3 = -\int_{B_r} \chi_{\mu/4} \langle \nabla \tilde{\beta}_n, \nabla |F_A^+| \rangle$.

(5.12)

With (5.11) understood, the next order of business is to derive suitable upper bounds for the norms of the terms that appear on its right hand side.

Step 6: The left most term on the right hand side of (5.11) (the explicit integral) is $c_0 r^2 r^2 K^2 N$ which is just a restating of the definition of N in (3.7). With regards to $\mathcal{T}_1$, use the fact that $|\nabla\nabla \chi_{\mu/4}|$ is at most $c_\mu \frac{1}{r^2}$ and Hölder's inequality to bound $|\mathcal{T}_1|$ by $c_\mu$ times the $L^2$ norm of $F_A^+$, and thus by $c_\mu r r K \sqrt{N}$. To bound $|\mathcal{T}_2|$ (and, later, $|\mathcal{T}_3|$), note that

$$|\nabla \tilde{\beta}_n| \le c_0 2^n K^{-1} |\nabla a|$$

(5.13)

and that this is supported in $U_n$. Since $|\nabla \chi_{\mu/4}| \le c_\mu \frac{1}{r}$, these observations lead to the bound $|\mathcal{T}_2| \le c_\mu 2^n \frac{1}{rK} r f_n$ with $f_n$ being the integral on the left hand side of (5.3). Given what is said by (5.3), it follows that $|\mathcal{T}_2| \le c_\mu \gamma^{2n} 2^n \sqrt{\mathfrak{f}}(rrK)(1+N)(\sqrt{N}+r)$. As for $|\mathcal{T}_3|$, use (3.7), the bound $|\nabla F_A^+| \le c_0 r^2 |a| |\nabla_A a|$ and the fact that $|a| \le c_0 2^{-n} K$ on the support of $\nabla \tilde{\beta}_n$ to see that $\mathcal{T}_3$ is no greater than $c_0 r^2 r^2 K^2 N$ also.

The inequality in (5.11) and the bounds from the preceding paragraph lead directly to the desired bound in (5.8) for the integral of $|F_A^+|$ over $B_{(1-\mu)r} - U_n$.

Step 7: The bounds in (5.1) and (5.8) lead to a bound for the $L^2$ norm of $|F_A^+|$ on $B_{(1-\mu)r}$ that has the integer n appearing as a parameter:

$$\int_{B_{(1-\mu)r}} |F_A^+| \le c_\mu (\gamma^n (rrK) + 2^{2n} + \gamma^n 2^{3n} \tfrac{1}{rrK})(1+N)(\sqrt{N}+r) r^2$$

(5.14)

The task now is choose the integer n to obtain the bound that is claimed by the proposition. Supposing that $r r \kappa \geq 2$, consider, for example, taking n to be the greatest integer that is less than $\frac{\ln(r r \kappa)}{3 \ln(2)}$. This choice leads from (5.14) to the bound

$$\int_{B_{(1-\mu)r}} |F_A^+| \leq c_\mu (r r \kappa)^{1-\sigma}(1+N)(\sqrt{N}+r)r^2$$

(5.15)

with $\sigma$ given by $\frac{|\ln(\gamma)|}{3 \ln(2)}$.

### b) Writing A as $\hat{A} + \mathfrak{b}$ with $\hat{A}$ having anti-self dual curvature

The proof of a proposition in the next subsection, one in Section 5d, and the proof of a proposition in Section 6 exploit the writing of the connection A given below.

**Proposition 5.3**: *There exists $\kappa > 4$ and, given $\mu \in (0, \frac{1}{100})$, there exists $\kappa_\mu > \kappa$; and these numbers have the following significance: Fix $r > 1$ and suppose that $(A, a)$ is a solution to (1.5). Fix $p \in X$ and $r \in (0, \kappa^{-1})$. The connection A can be written on $B_r$ as $\hat{A} + \mathfrak{b}$ where $\hat{A}$ and $\mathfrak{b}$ are described by the subsequent bullets.*

- *$\hat{A}$ is a Hermitian connection on E over $B_r$ and $\mathfrak{b}$ is an $i\mathbb{R}$ valued 1-form on $B_r$.*
- *The curvature 2-form of $\hat{A}$ is anti-self dual.*
- *Let $\kappa$ denote $\kappa(r)$ and let $\tilde{\kappa}$ denote $\kappa(\frac{1}{1-\mu}r)$. By the same token, let N denote $N(r)$ and let $\tilde{N}$ denote $N(\frac{1}{1-\mu}r)$. Then $\mathfrak{b}$ has the properties in the list below.*
  
  a) $\int_{B_r} |\mathfrak{b}|^2 \leq \kappa (r r \kappa)^2 (N + r^2) r^2$.
  
  b) $\int_{B_r} |\mathfrak{b}|^2 \leq \kappa_\mu (r r \tilde{\kappa})^{2-1/\kappa}(1+\tilde{N})(\sqrt{\tilde{N}}+r)(\sqrt{N}+r)r^2$.

- *With $\kappa$ and N as in the previous bullet:*
  
  a) $\int_{B_{(1-\mu)r}} |\nabla \mathfrak{b}|^2 \leq \kappa_\mu \int_{B_r} |F_A^+|^2$.
  
  b) $\int_{B_{(1-\mu)r}} |\nabla \mathfrak{b}|^2 \leq \kappa_\mu (r r \kappa)^2 (N+r^2)$.
  
  c) $\int_{B_{(1-\mu)r}} |\nabla \nabla \mathfrak{b}|^2 \leq \kappa_\mu \frac{1}{r^2}(r r \kappa)^4(r^2 + N)$.

*Proof of Proposition 5.3*: The proof has four parts.

*Part 1*: As explained momentarily, there is a unique, smooth section of $\Lambda^+$ on $B_r$ (to be denoted by u) that obeys the equations

- $(dd^\dagger u)^+ = F_A^+$
- $u = 0$ *on* $\partial B_r$.

(5.16)

Let $\mathfrak{b} = d^\dagger u$. Then $\hat{A} = A - \mathfrak{b}$ has anti-self dual curvature 2-form. Thus, if u exists, then it remains only to verify that $\mathfrak{b}$ is described by bullets three and four of Proposition 5.3.

The existence of a (unique) section of $\Lambda^+$ obeying (5.16) follows from the fact that the operator $(dd^\dagger(\cdot))^+$ on $C^\infty(\Lambda^+|_{B_r})$ with Dirichelet boundary conditions is uniformly positive in the following sense: If $r < c_0^{-1}$ and if w is a smooth section of $\Lambda^+$ on $B_r$ vanishing on $\partial B_r$, then

$$\int_{B_r} |d^\dagger w|^2 \geq c_0^{-1} \left( \int_{B_r} |\nabla w|^2 + \frac{1}{r^2} \int_{B_r} |w|^2 \right)$$

(5.17)

Keep in mind also that there is a Bochner-Weitzenboch formula that writes $(dd^\dagger w)^+$ as

$$(dd^\dagger w)^+ = \nabla^\dagger \nabla w + \mathcal{R} w$$

(5.18)

with $\mathcal{R}$ being an endomorphism whose components are linear functions of the Riemann curvature tensor.

*Part 2*: Suppose now that u obeys (5.16). This part of the proof derives a bound for the $L^2$ norm of $d^\dagger u$ (and thus the $L^2$ norm of $\mathfrak{b} = (d^\dagger u)\sigma$):

$$\int_{B_r} |d^\dagger u|^2 \leq c_0 \, r r_K (\sqrt{N} + r) \int_{B_r} |F_A^+| \,.$$

(5.19)

Granted (5.19) for the moment, then Item a of the third bullet in Proposition 5.3 follows from this and the Cauchy-Schwarz inequality (and (3.7)); and Item b) of the third bullet follows from (5.19) and Proposition 5.2. (Keep in mind that $\kappa$ is an increasing function on $(0, c_0^{-1})$ which is a consequence of (3.8) and (3.9). Thus, $\kappa \leq \tilde{\kappa}$.)

To derive (5.19), take the inner product of both sides of the equation in the top bullet of (5.16) with u and integrate over $B_r$. An integration by parts leads to the identity

$$\int_{B_r} |d^\dagger u|^2 = \int_{B_r} \langle u, F_A^+ \rangle \,.$$

(5.20)

Bringing Hölder's inequality to bear on the right hand side of (5.20) leads to

$$\int_{B_r} |d^\dagger u|^2 \le \left(\int_{B_r} |u|^4\right)^{1/4} \left(\int_{B_r} |F_A^+|^{4/3}\right)^{3/4}.$$

(5.21)

Since the $L^4$ norm in (5.21) is no greater than $c_0$ times the $L^2$ norm of $\nabla u$ (courtesy of (2.1)) and the latter is no greater than $c_0$ times the $L^2$ norm of $d^\dagger u$ (courtesy of (5.17)), the inequality in (5.21) bounds the $L^2$ norm of $d^\dagger u$ by the $L^{4/3}$ norm of $F_A^+$. The latter bound leads (using Hölder's inequality) to the bound

$$\int_{B_r} |d^\dagger u|^2 \le c_0 \left(\int_{B_r} |F_A^+|^2\right)^{1/2} \int_{B_r} |F_A^+|.$$

(5.22)

Since the $L^2$ norm of $F_A^+$ is at most $r r_K \sqrt{N}$ (see (3.7)), the bound in (5.22) gives (5.19).

*Part 3*: The argument that follows proves that $\mathfrak{b}$ (which is $(d^\dagger u)$) obeys Item a) of the fourth bullet in Proposition 5.3. To this end, note first that the left hand side of (5.22) is no greater than $c_0 r^2$ times the square of the $L^2$ norm of $F_A^+$ on $B_r$ (by an appeal to Hölder's inequality). Therefore, this is also a bound for the square of the $L^2$ norm of $\mathfrak{b}$ (which is $d^\dagger u$) on $B_r$. Keep this fact for later use.

By way of a reminder, $\chi_\mu$ is a non-negative function on $B_r$ with compact support that is equal to 1 on $B_{(1-\mu)r}$ and whose differential has norm bounded by $c_0 \frac{1}{\mu r}$. Integration by parts and the fact that $d^\dagger \mathfrak{b} = 0$ (because $\mathfrak{b}$ is $d^\dagger u$) leads to

$$\int_{B_r} \chi_\mu |\nabla \mathfrak{b}|^2 \le \int_{B_r} \chi_\mu |(d\mathfrak{b})^+|^2 + c_0 \left(1 + \frac{1}{\mu^2 r^2}\right) \int_{B_r} |\mathfrak{b}|^2.$$

(5.23)

The left most term on the right hand side of (5.23) is no greater than the square of the $L^2$ norm of $F_A^+$ on $B_r$ because $(d\mathfrak{b})^+ = F_A^+$; and granted what was said in the preceding paragraph about the $L^2$ norm of $\mathfrak{b}$, this is also the case for the right most term on the right hand side of (5.23). These remarks lead directly to Item b) of the third bullet of the proposition since the square of the $L^2$ norm of $F_A^+$ on $B_r$ is no greater than $r^2 r_K^2 N$.

Item b) of the fourth bullet of the proposition follows from Item a) using (3.7).

*Part 4*: To derive Item c) of the fourth bullet, first differentiate the equations $(d\mathfrak{b})^+ = F_A^+$ and $d^\dagger \mathfrak{b} = 0$ to obtain an elliptic, first order equation for $\nabla \mathfrak{b}$ of the form $D(\nabla \mathfrak{b}) = Q$ with $|Q| \le c_0(r^2|a||\nabla_A a| + |\mathfrak{b}|)$. Multiply both sides of the latter equation by $\chi_\mu$ and then integrate the square of the norm of both sides of the resulting identity over $B_r$. The bound in question then follows from this integral identity after a straightforward integration by parts using the defintion of N in (3.6) and using the following already established results: First, the bound for $|a|$ in Proposition 3.1 with the number r replaced by $r' = (1 - \frac{1}{64}\mu)r$ and, having made this replacement, then replacing $\mu$ with $\mu' = \frac{1}{1024}\mu$. Second, Item a) of the third bullet of Proposition 5.3 and Items a) and b) of the fourth bullet, but with the number r replaced by $r'$ and the number $\mu$ replaced by $\mu'$ in each of these three items.

### c) Proof of Proposition 5.1

The proposition's first bullet follows from the second bullet or the third bullet with the function M given by the rule $s \to M(s) = c_0 r^2 s^2 K(s)^2 N(s)$. To see that this is so, note first that the integral of $|F_A^+|^2$ over a ball of radius s is, in any event, no greater than $r^2 s^2 K(s)^2 N(s)$. (This follows from the definition of N in (3.7)). Therefore, the second or third bullets of Proposition 5.1 can be invoked $r_1 = r$ and $r_0 = 8r$ using the preceding version of the function M. The corresponding version of the third bullet of Proposition 5.1 leads directly to the assertion made by the first bullet. If the second bullet of Proposition 5.1 is invoked, then the first bullet follows because (as explained below) the value of $r_0^2 K(r_0)^2 N(r_0)$ can't be greater then $\frac{2}{\ln(2)} r_0^2 K(2r_0)^2$ which is no greater than $c_0 r_0^2$. Therefore, $M(r_0)$ with $r_0 = 8r$ is no greater than $c_0 r^2 r^2$.

To restate the claim made in the preceding paragraph: For any $r \in (0, c_0^{-1}]$,

$$r^2 K(r)^2 N(r) \le \frac{2}{\ln(2)} r^2 K(2r)^2 .$$

(5.24)

To prove this, note first that the function $r \to r^2 K(r)^2 N(r)$ is increasing (look at (3.6)), and as a consequence

$$r^2 K(r)^2 N(r) \le \frac{1}{\ln(2)} \int_r^{2r} s^2 K(s)^2 N(s) \frac{ds}{s} .$$

(5.25)

This implies in turn what is asserted by (5.24) because the factor of $s^2$ in the integrand is less than $4r^2$ and $K^2 \frac{N}{s}$ is $\frac{1}{2} \frac{d}{ds} K^2$ due to (3.8).

The six steps that follow prove the second and third bullets of the lemma.

<u>Step 1</u>: Reintroduce the section $\mathfrak{b}$ from Proposition 5.3. Since $(d\mathfrak{b})^+ = F_A^+$, the connection $\hat{A} = A - \mathfrak{b}$ has anti-self dual curvature tensor. Fix $\delta > 0$ and $s \in (0, \frac{3}{4} r]$. Use the triangle inequality and Item b) of the third bullet of Proposition 5.3 (with $\mu$ small and r replaced by $(1-\mu)r$) to see that

- $\int_{B_s} |F_A^-|^2 \leq (1-\delta) \int_{B_s} |F_{\hat{A}}|^2 + c_\mu \delta^{-1} \int_{B_r} |F_A^+|^2$.

- $\int_{B_s} |F_A^-|^2 \geq (1-\delta) \int_{B_s} |F_{\hat{A}}|^2 - c_\mu \delta^{-1} \int_{B_r} |F_A^+|^2$.

(5.26)

Note in particular how the $L^2$ norm of $F_A^-$ is determined for the most part by the $L^2$ norm of $F_{\hat{A}}$ when the $L^2$ norm of $F_A^+$ is small.

Step 2: With the preceding in mind, the following basic fact about $\mathbb{R}$-valued closed, anti-self dual 2-forms will be brought to bear: Suppose that $r \in (0, c_0^{-1}]$, that $s \in (0, r]$, and that $\omega$ is anti-self dual on $B_r$. Then

$$\int_{B_s} |\omega|^2 \leq (\tfrac{s}{r})^{4-c_0 r^2} \int_{B_r} |\omega|^2 .$$

(5.27)

This bound is proved in Step 6. Accept it as true for now.

Since $iF_{\hat{A}}$ is an $\mathbb{R}$-valued, closed and anti-self dual 2-form, the inequality in (5.27) can be used with $\omega = iF_{\hat{A}}$. With this version of (5.27) available, invoke (5.26) (once using the top bullet and again using the lower) to see that

$$\int_{B_s} |F_A^-|^2 \leq (1+c_0\delta)(\tfrac{s}{3r/4})^{4-c_0 r^2} \int_{B_{3r/4}} |F_A^-|^2 + c_0 \delta^{-1} \int_{B_r} |F_A^+|^2 .$$

(5.28)

when $s \in (0, \tfrac{3}{4} r]$. Of particular interest is the $s = \tfrac{1}{2} r$ version of (5.28):

$$\int_{B_{r/2}} |F_A^-|^2 \leq (1+c_0\delta)(\tfrac{2}{3})^{4-c_0 r^2} \int_{B_{3r/4}} |F_A^-|^2 + c_0 \delta^{-1} \int_{B_r} |F_A^+|^2 .$$

(5.29)

Step 3: Let $L(r) = M(r)(\int_{B_r} |F_A^+|^2)^{-1}$. Assume that $L$ is in any event larger than $c_0 \delta^{-2} 2^4$ and use it to distinguish two possible cases (designated CASE$_r$- and CASE$_r$+):

- CASE$_r$- *occurs when* $\int_{B_r} |F_A^-|^2 \leq L(r) \int_{B_r} |F_A^+|^2$.

- CASE$_r$+ *occurs when* $\int_{B_r} |F_A^-|^2 > L(r) \int_{B_r} |F_A^+|^2$.

(5.30)

The distinction has the following significance: If CASE$_r$- occurs, then with no further ado,

$$\int_{B_{r/2}} |F_A|^2 \leq (L(r)+1) \int_{B_r} |F_A^+|^2 \ .$$

(5.31)

(This is because the integral of $|F_A|^2$ is the sum of the integrals of $|F_A^+|^2$ and $|F_A^-|^2$ and their integrals over the half-radius ball are less than their integrals over the full radius ball.) Anyway, nothing more will be said if CASE$_r$- occurs. In the event that CASE$_r$+ happens (and if $r < c_0^{-1}\delta$), then (5.29) leads to

$$\int_{B_{r/2}} |F_A|^2 \leq (1+c_\diamond\delta) \, (\tfrac{2}{3})^{4-\delta} \int_{B_r} |F_A|^2$$

(5.32)

with $c_\diamond$ a number greater than 1 that is independent of r, p, $\delta$, the pair (A, $a$) and $r$.

Step 4: Supposing that (5.32) holds, if CASE$_{2r}$- is relevant (the top bullet of (5.30) with r replaced by 2r), then (5.32) with the version of (5.31) that has r replaced by 2r leads to the bound

$$\int_{B_{r/2}} |F_A|^2 \leq (1+c_\diamond\delta) \, (\tfrac{2}{3})^{4-\delta} (L(2r)+1) \int_{B_{2r}} |F_A^+|^2 \ .$$

(5.33)

Nothing more will be said if CASE$_{2r}$- is relevant. On the otherhand, if $2r < c_\diamond^{-1}\delta$ and if CASE$_{2r}$+ is relevant, then (5.32) implies in turn

$$\int_{B_{r/2}} |F_A|^2 \leq (1+c_\diamond\delta)^2 \, (\tfrac{2}{3})^{2(4-\delta)} \int_{B_{2r}} |F_A|^2 \ .$$

(5.34)

Given (5.34), ask which of CASE$_{4r}$- or CASE$_{4r}$+ is relevant to say more about its right hand side; and if it is CASE$_{4r}$+, then one can ask about CASE$_{8r}$- or CASE$_{8r}$+, and so on. Continuing in this vein produces a positive integer N with the following two properties: The first is that CASE$_\rho$+ occurs for each $\rho \in \{r, 2r, \ldots, 2^{N-1}r\}$; and the second is that either CASE$_\rho$- occurs for $\rho = 2^N r$ or else $c_\diamond^{-1}\delta \in [2^{N-1}r, 2^N r]$. If the former and $2^N r \leq c_\diamond^{-1}\delta$, then

$$\int_{B_{r/2}} |F_A|^2 \leq (1+c_\diamond\delta)^N \, (\tfrac{2}{3})^{N(4-\delta)} (L(2^N r)+1) \int_{B_{2^N r}} |F_A^+|^2 \ .$$

(5.35)

On the other hand, if $c_\diamond^{-1}\delta \in [2^{N-1}r, 2^N r]$; then

$$\int_{B_{r/2}} |F_A|^2 \leq c_0 (1+c_0\delta)^N \left(\tfrac{2}{3}\right)^{N(4-\delta)} r^2 .$$

(5.36)

Indeed, the latter inequality follows because the integral of $|F_A|^2$ over the whole of X is at most $c_0 r^2$. (See (2.11).)

Step 5: To obtain the statement of the second bullet of Proposition 5.1, let $r_0$ be as given by the proposition and set $r = 2r_0$. Then the inequality in (5.35) holds for N such that $2^N r_0 < r_1 < 2^{N+1} r_0$. This understood, set $\kappa \geq c_0 \delta^{-1}$ for $\delta < \tfrac{1}{100} c_0^{-1}$. The conclusions of Proposition 5.1's second bullet follow from (5.35) because $\left(\tfrac{2}{3}\right)^{N(4-\delta)} \leq c_0 \left(\tfrac{r_0}{r_1}\right)^{2+\delta/c_0}$ when $\delta$ is small. The third bullet of Propositon 5.1 follows from (5.36) by again taking $\delta \leq \tfrac{1}{100} c_0^{-1}$ and taking $\kappa > c_0 \delta^{-1}$.

Step 6: The proof of (5.27) starts by writing $\omega$ as $dv$ with $v$ being an $\mathbb{R}$-valued 1-form on $B_r$ (this is doable because $d\omega = 0$.) Since $\omega = dv$ is anti-self dual, the Hodge star of the 4-form $dv \wedge dv$ is $-|\omega|^2$. Therefore, letting $\iota: \partial B_s \to B_r$ denote the tautological inclusion map, Stoke's theorem can be invoked to see that

$$\int_{B_s} |\omega|^2 = - \int_{\partial B_s} \iota^* v \wedge d(\iota^* v)$$

(5.37)

Holding onto (5.37), note that the smallest absolute value of any eigenvector of the operator $*d(\cdot)$ on the space of coclosed 1-forms on $\partial B_s$ differs from $\tfrac{2}{s}$ by at most $c_0 s^2$. This is true when $B_s$ is a Euclidean ball in $\mathbb{R}^4$ without the error term $c_0 s^2$; and so it is true for $B_s$ here, with the error term, because the metric on $\partial B_s$ for small s (for $s \leq c_0^{-1}$) differs by at most $c_0 s^2$ from the induced metric on the Euclidean sphere of radius s in $\mathbb{R}^4$. Use of this eigenvalue bound in (5.37) leads to the inequality

$$\int_{B_s} |\omega|^2 \leq (1+c_0 s^2) \tfrac{1}{2} s \int_{\partial B_s} |\iota^*(dv)|^2 .$$

(5.38)

Since $|\iota^*(dv)|^2 = \tfrac{1}{2} |dv|^2$, this inequality can be written entirely in terms of $|\omega|$ as

$$\int_{B_s} |\omega|^2 \leq (1+c_0 s^2) \tfrac{1}{4} s \int_{\partial B_s} |\omega|^2 .$$

(5.39)

To continue, let $f$ denote the function on $(0, r]$ given by the integral on the left hand side of (5.39). The inequality in (5.39) can be written as

$$\frac{d}{ds} \ln(f) \geq \frac{4}{s} - c_0 s$$

(5.40)

because the integral that appears on the right hand side of (5.39) is the derivative of $f$. Integration of (5.40) leads directly to (5.27).

### d) The $L^2_1$ norm of $F_A$

The upcoming proposition asserts an a priori bound for the $L^2$ norm of $\nabla F_A$ on small radius balls when $(A, a)$ obeys some $r > 1$ version of (1.5).

**Proposition 5.4**: *There exists $\kappa > 4$ and, given $\mu \in (0, \frac{1}{100})$, there exists $\kappa_\mu > \kappa$; and these numbers have the following significance: Fix $r > 1$ and suppose that $(A, a)$ is a solution to (1.5). Fix $p \in X$ and $r \in (0, \kappa^{-1})$. Then*

$$\int_{B_{(1-\mu)r}} |\nabla F_A|^2 \leq \kappa_\mu \frac{1}{r^2} \left( \int_{B_r} |F_A^-|^2 + r^4 r^4 K^4 (N + r^2) \right).$$

***Proof of Proposition 5.4***: Invoke Proposition 5.3 to write $F_A$ as $F_{\hat{A}} + d\flat$ with $F_{\hat{A}}$ being an $i\mathbb{R}$ valued, anti-self dual 2-form. Having done so, then write $\nabla F_A$ as $\nabla d\flat + \nabla F_{\hat{A}}$. The proposition's inequality follows from this rewriting of $\nabla F_A$ using Items a) and b) of the fourth bullet of Proposition 5.3 and the fact that the $L^2_1$ norm on $B_{(1-\mu)r}$ of a closed, anti-self dual 2-form on $B_r$ is bounded by $c_\mu \frac{1}{r}$ times its $L^2$ norm on $B_r$. (The argument for latter bound goes as follows: Let $\omega$ denote an anti-self dual 2-form. A straightforward integration by parts bounds the integral over $B_r$ of the function $\chi_\mu^2 |\nabla \omega|^2$ by the integral over $B_r$ of $(\chi_\mu^2 |d\omega|^2 + c_\mu \frac{1}{r^2} |\omega|^2)$. In the case when $d\omega = 0$, this last bound gives the desired a priori bound for the integral on $B_{(1-\mu)r}$ of $|\nabla \omega|^2$.)

### 6. The derivative of the frequency function

Suppose that $r > 1$ and $(A, a)$ are a pair of connection on P and section of the bundle $\Lambda^+ \otimes (P \times_{SO(3)} \mathfrak{S})$ that obey (1.5). Fix a point $p \in X$ and use $(A, a)$ to define the functions K and N using the formulas in Section 3b. Sections 6a and 6b state formulas for the derivative of N. Section 6c introduces a very useful local average of N (it is denoted by N̲) and Section 6d gives a formula for the derivative of this local average.

### a) The first formula for the derivative of N

A formula for the derivative of N is obtained directly by differentiating (with respect to r) the right hand side of (3.6) and invoking (3.8). This is the result:

$$\frac{d}{dr} N = \frac{1}{r^2 K^2} \int_{\partial B_r} (|\nabla_A a|^2 + r^2 |\langle a, \tau a \rangle|^2) - \frac{2}{r} N(1+N).$$

(6.1)

This formula has the following immediate consequence:

**Lemma 6.1**. *There exists $\kappa > 1$ with the following significance: Fix $r > 1$ and suppose that $(A, a)$ is solution to (1.5). Fix a point in $X$ so as to define the function N. Suppose that $r \in (0, \kappa^{-1})$ and $s \in (0, r]$. Then*

- $N(r) \geq N(s) \left(\frac{s}{r}\right)^2 \frac{1}{1 + N(s)(1 - (\frac{s}{r})^2)}.$
- $N(s) \leq N(r) \left(\frac{r}{s}\right)^2 \frac{1}{1 - N(r)((\frac{r}{s})^2 - 1)}.$
- *If $\delta \in (0, 1)$ and if $N(r) \leq \delta$, and if $s \in (\frac{\sqrt{2\delta}}{\sqrt{1+2\delta}} r, r)$, then $N(s) \leq 2\delta \left(\frac{r}{s}\right)^2$.*

The second bullet says (roughly) that if $N(r)$ is less than 1, then $N(s)$ for $s < r$ can not be all that much larger than $N(r)$ if $s/r$ is not small.

*Proof of Lemma 6.1*: It follows from (6.1) that $\frac{d}{dr} N \geq -\frac{2}{r} N(1+N)$, which is to say that

$$\frac{d}{dr} \left( \ln \left( \frac{N}{N+1} \right) \right) \geq -\frac{2}{r}.$$

(6.2)

This inequality integrates to give the inequality that is asserted by the lemma's first bullet. The second bullet's inequality is just a rewriting of the first bullet's inequality. Meanwhile, the assertion of the third bullet is an application of the second bullet.

By way of a parenthetical remark, the inequality in the second bullet of Lemma 6.1 leads to the following uniform bound for the product $K^2 N$ which comes from the following: Suppose that $r \in (0, c_0^{-1}]$ and that $s \in (0, r)$. Then

$$K(s)^2 N(s) \leq \frac{1}{1 - (\frac{s}{r})^2} (K(r)^2 - K(s)^2).$$

(6.3)

This follows from (3.9) when the first bullet of Lemma 6.1 (with $\tau$ used in lieu of r) is invoked to replace $N(\tau)$ in (3.9) with $N(s) \left(\frac{s}{\tau}\right)^2 \frac{1}{1 + N(s)(1 - (\frac{s}{\tau})^2)}.$

**b) A second formula for the derivative of N**

A formula for the derivative of N is given below in Proposition 6.2. Having fixed a point $p \in X$, the proposition introduces $\nabla_{A,r} a$ to denote the section of $\mathbb{S}_E^+$ on the complement of p in $B_{1/c_0}$ that gives the pairing between $\nabla_A a$ and the unit length vector field that points outward along the geodesic arcs that start at p. The proposition also uses $E_A$ to denote the $i\mathbb{R}$ valued 1-form on the complement of p in $B_{1/c_0}$ that gives the pairing between $F_A$ and the same unit length vector field. Define $B_A$ to be the pairing between $*F_A$ and this outward pointing radial vector field. The inner product between $E_A$ and $B_A$ is denoted by $\langle E_A, B_A \rangle$ which is a function on $B_r - p$.

**Proposition 6.2**: *There exists $\kappa > 1$ with the following significance: Fix $r > 1$ and suppose that $(A, a)$ is a solution to (1.5). Fix a point $p \in X$ so as to define the functions K and N as done in Section 3b. If $r \in (0, \kappa^{-1})$, then*

$$\tfrac{d}{dr} \text{N} = \tfrac{2}{r^2 \text{K}^2} \int_{\partial B_r} |\nabla_{A,r} a - \tfrac{1}{r} \text{N} a|^2 + \tfrac{1}{r^2 r^2 \text{K}^2} \int_{\partial B_r} |E_A|^2 + \tfrac{1}{r^2 r^2 \text{K}^2} \int_{\partial B_r} \langle E_A, B_A \rangle + \mathfrak{q} ,$$

*with $\mathfrak{q}$ being a function on $[0, \kappa^{-1}]$ whose absolute value obeys $|\mathfrak{q}| \leq \kappa r (1 + \text{N})$.*

(Note that $\mathfrak{q}$ is zero if the metric on $B_r$ is flat.)

Analogs of Proposition 6.2 appear in [T1], [T3] and [HW]. The proof of Proposition 6.2 that follows was guided by what is done in [T1] and [T3]. Meanwhile, the proofs of the analog of Proposition 6.2 in [T1] and [T3] was guided by proofs of propositions about analogs of N in [Al], [DF] and [HHL].

*Proof of Proposition 6.2*: The proof starts with the formula in (6.1) and then procedes to completion in seven steps.

<u>Step 1</u>: Fix an orthonormal frame for $T^*X$ on $B_{1/c_0}$. Use this frame to write the 4 directional covariant derivatives of $\nabla_A a$ along the dual vector fields as $\{\nabla_{A,\alpha} a\}_{\alpha \in \{1,2,3,4\}}$ and for each $\alpha \in \{1, 2, 3, 4\}$. Then write the 3 components of $\nabla_{A,\alpha} a$ with respect to a corresponding orthonormal frame for $\Lambda^+$ as $\{(\nabla_{A,\alpha} a)^c\}_{c \in \{1,2,3\}}$. Write the components of the curvature $F_A^+$ with respect to the chosen frame as $\{F_A^+{}_{\alpha\beta}\}_{\alpha,\beta \in \{1,2,3,4\}}$ and do likewise for $F_A^-$. With this notation understood, the next equation defines a symmetric section over $B_{1/c_0}$ of $T^*X \otimes T^*X$ to be denoted by T. The components $\{T_{\alpha\beta}\}_{\alpha,\beta \in \{1,2,3,4\}}$ are

$$T_{\alpha\beta} = \tfrac{1}{2} (\langle \nabla_{A,\alpha} a, \nabla_{A,\beta} a \rangle + \langle \nabla_{A,\beta} a, \nabla_{A,\alpha} a \rangle) + r^{-2} \langle F_A^+{}_{\alpha\nu} F_A^-{}_{\beta\nu} \rangle - \tfrac{1}{2} \delta_{\alpha\beta} (|\nabla_A a|^2 + \tfrac{1}{2} r^2 |\langle a, \tau a \rangle|^2) .$$

(6.4)

The formal, $L^2$ adjoint of the covariant derivative maps sections of $T^*X \otimes T^*X$ to $T^*X$. This operator is denoted by $\nabla^\dagger$. Step 6 proves that $\nabla^\dagger T = 0$ when the metric is flat. When the metric is not flat, Steps 6 and 7 prove that $\nabla^\dagger T$ can be written as

$$\nabla^\dagger T = \mathfrak{r} \tag{6.5}$$

where the notation has $\mathfrak{r}$ denoting a 1-form whose norm obeys $|\mathfrak{r}| \leq c_0 |a| |\nabla_A a|$. A norm bound of this sort exists because $\mathfrak{r}$ can be written as $\mathfrak{r}_{a\beta c} \langle a_a \nabla_{A,\beta} a_c \rangle$ with each a, c $\in \{1,2,3\}$ and $\beta \in \{1, 2, 3, 4\}$ version of $\mathfrak{r}_{a\beta c}$ being a 1-form whose components are proportional to components of the Riemann curvature tensor.

Step 2: Use $\hat{x}$ to denote the differential of the function $\frac{1}{2}\text{dist}(p, \cdot)^2$. The norm of $\hat{x}$ is the function $\text{dist}(p,\cdot)$, and its dual is a vector field that is tangent to the geodesic rays from the point p. The covariant derivative of $\hat{x}$ can be written as $\mathfrak{m} + \mathfrak{z}$ with $\mathfrak{m}$ being the metric tensor and $\mathfrak{z}$ being a tensor with $|\mathfrak{t}|$ bounded by $c_0 \text{dist}(p,\cdot)^2$ and with $|\nabla \mathfrak{t}|$ bounded by $c_0 \text{dist}(p,\cdot)$. The notation also uses $E_A^+$ to denote the $i\mathbb{R}$-valued 1-form on the complement of p in $B_{1/c_0}$ that gives the pairing between $F_A^+$ and the unit length vector field that is tangent and outward pointing along the geodesic arcs that start at p. Define $E_A^-$ to be the pairing between $F_A^-$ and this same outward pointing radial vector field. The inner product between $E_A^+$ and $E_A^-$ is denoted by $\langle E_A^+, E_A^- \rangle$, this being a function on $B_r - p$.

With the preceding understood, take the inner product of both sides of (6.5) with $\hat{x}$ and integrate the resulting identity over $B_r$. Having done so, integrate by parts to remove derivatives from T and $\mathfrak{R}$ to derive the identity

$$\frac{1}{2} \int_{\partial B_r} (|\nabla_A a|^2 + \frac{1}{2} r^2 |\langle a,\tau a \rangle|^2) = \int_{\partial B_r} (|\nabla_{A,r} a|^2 + \frac{1}{r^2} \langle E_A^+, E_A^- \rangle) + \frac{1}{r} \int_{B_r} (|\nabla_A a|^2 + r^2 |\langle a,\tau a \rangle|^2)$$
$$+ \frac{1}{r} \int_{B_r} (\mathfrak{m}(\mathfrak{z},T) - \mathfrak{m}(\hat{x},\mathfrak{r}))$$

(6.6)

with $\mathfrak{m}(\cdot,\cdot)$ denoting the metric inner product on both $T^*X$ and $T^*X \otimes T^*X$.

Step 3: Use (6.6) with the definition of N to rewrite (6.1):

$$\frac{d}{dr} N = \frac{2}{r^2 K^2} \int_{\partial B_r} |\nabla_{A,r} a|^2 - \frac{2}{r} N^2 + \frac{1}{r^2 r^2 K^2} \int_{\partial B_r} (\frac{1}{2} r^4 |\langle a,\tau a \rangle|^2 + 2 \langle E_A^+, E_A^- \rangle)$$
$$- \frac{2}{r^3 K^2} \int_{B_r} \mathfrak{m}(\hat{x},\mathfrak{r}) + \frac{2}{r^3 K^2} \int_{B_r} \mathfrak{m}(\mathfrak{z},T).$$

(6.7)

What are denoted in (6.7) as $E_A^+$ and $E_A^-$ can be written in terms of the vector fields $E_A$ and $B_A$ as $E_A^+ = \frac{1}{2}(E_A + B_A)$ and $E_A^- = \frac{1}{2}(E_A - B_A)$. Use this rewriting with the equation $F_A^+ = \frac{1}{\sqrt{2}} r^2 \langle a, \tau a \rangle$ from (2.5) to write

$$\tfrac{1}{2} r^4 |\langle a, \tau a \rangle|^2 + 2 \langle E_A^+, E_A^- \rangle = |E_A|^2 + \langle E_A, B_A \rangle .$$
(6.8)

This identity is used to rewrite the term in (6.7) that involves the integral over $\partial B_r$ of the function $\tfrac{1}{2} r^4 |[a; a]|^2 + 2\langle E_A^+, E_A^- \rangle$.

The next task is to rewrite the term $-\tfrac{2}{r} N^2$ in (6.7). This is done by using Stokes theorem and the identity in (2.9) to write $N$ as

$$N = \tfrac{1}{r^2 K^2} \int_{\partial B_r} \langle a, \nabla_{A,r} a \rangle - \tfrac{1}{r^2 K^2} \int_{B_r} \Re(\langle a \otimes a \rangle) .$$
(6.9)

Use this identity to replace the two left most terms on the right hand side of (6.7) by

$$\tfrac{2}{r^2 K^2} \int_{\partial B_r} |\nabla_{A,r} a - \tfrac{1}{r} N a|^2 + \mathfrak{e} N$$
(6.10)

with $\mathfrak{e}$ being a function on $(0, c_0^{-1}]$ whose absolute value obeys

$$|\mathfrak{e}| \le c_0 \tfrac{1}{r^3 K^2} \left( \int_{B_r} |a|^2 \right).$$
(6.11)

This in last term is no greater than $c_0 r$ as can be seen by writing the integral over $B_r$ of $|a|^2$ as an integral of the function $s \to s^3 K(s)^2$ on the interval $[0, r]$ and using the fact that $K$ is an increasing function of $r$.

Replacing the left most two terms on the right hand side of (6.7) with (6.10) and using (6.8)'s identity results in an equation for $N$'s derivative that has the form depicted in the proposition with $\mathfrak{q}$ being the sum of $\mathfrak{e} N$ and the terms in (6.7) that are proportional to the integrals over $B_r$ of $\mathfrak{m}(\hat{x}, \mathfrak{r})$ and $\mathfrak{m}(\mathfrak{z}, T)$.

Step 4: This step and Step 5 prove $|\mathfrak{q}| \le c_0 r (1 + N)$. This step derives a preliminary bound for $\mathfrak{q}$ which is the upcoming (6.12). This preliminary bound uses $\langle F_A^+ \circ F_A^- \rangle$ to denote the symmetric bilinear form on $TX$ that assigns the inner product of $F_A^+(v_1)$ with $F_A^-(v_2)$ to any given pair of vectors $v_1, v_2 \in TX$. Here is the promised bound:

$$|q| \leq c_0 r (1 + N) + c_0 \frac{1}{r^2 r K^2} \left| \int_{B_r} \langle F_A^+ \circ F_A^- \rangle(t) \right|$$

(6.12)

with t denoting here a symmetric section of $\otimes_2 TX$ on $B_r$ obeying $|t| \leq c_0$ and $|\nabla t| \leq c_0 \frac{1}{r}$.
Step 6 explains why

$$\frac{1}{r^2 r K^2} \left| \int_{B_r} \langle F_A^+ \circ F_A^- \rangle(t) \right| \leq c_0 r + c_0 \frac{1}{r K^2} \int_{B_r} (|\nabla_A a|^2 + r^2 |\langle a, \tau a \rangle|^2).$$

(6.13)

Since the right hand side of (6.13) is less than $c_0 r(1+N)$ (by virtue of (3.6)), it then follows directly from (6.12) and (6.13) that $|q|$ is less than $c_0 r(1+N)$ as claimed.

The derivation of (6.12) starts with the contribution to q from the $\mathfrak{e} N$ term in (6.10). As noted in the previous step, $|\mathfrak{e}| \leq c_0 r$, and by virtue of this, $|\mathfrak{e}| N$ is bounded by what is written on the right hand side of (6.12).

To continue, consider next the contribution to q of the term

$$\frac{2}{r^3 K^2} \int_{B_r} \mathfrak{m}(\mathfrak{z}, T)$$

(6.14)

from (6.7). The absolute value of the integral of $\mathfrak{m}(\mathfrak{z}, T)$ over $B_r$ is no greater than the sum of $c_0 r^2$ times the integral of $|\nabla_A a|^2 + r^2 |\langle a, \tau a \rangle|^2$ and $c_0 r^2$ times the absolute value of the integral of $r^{-2} \langle F_A^+ \circ F_A^- \rangle (r^{-2} \mathfrak{z})$. This is to say that the absolute value of what is written in (6.14) obeys the asserted bound for $|q|$ in (6.12) with $t = r^{-2} \mathfrak{z}$. Note in this regard that $|\mathfrak{z}| \leq c_0 r^2$ and $|\nabla \mathfrak{z}| \leq c_0 r$ so the norm of this version of t is bounded $c_0$ and the norm of its covariant derivative is bounded by $c_0 r^{-1}$.

The final contribution to q comes from the term in (6.7) with the integral of $\mathfrak{m}(\hat{x}, \mathfrak{r})$. Since $|\hat{x}| \leq c_0 r$ and $|\mathfrak{r}| \leq c_0 |a| |\nabla_A a|$, this contribution is at most

$$c_0 \frac{1}{r^2 K^2} \left( \int_{B_r} |a|^2 \right)^{1/2} \left( \int_{B_r} |\nabla_A a|^2 \right)^{1/2},$$

(6.15)

which is, in turn, at most

$$c_0 \frac{1}{r^3 K^2} \int_{B_r} |a|^2 + c_0 \frac{1}{r K^2} \int_{B_r} |\nabla_A a|^2.$$

(6.16)

The left most term in (6.16) is what appears on the right hand side of (6.11) and, as noted subsequent to (6.11), it is bounded by $c_0 r$. Meanwhile, the right most term in (6.16) is no greater than $c_0 r N$ Thus, (6.15) is bounded by $c_0 r (1 + N)$.

Step 5: The bound that is asserted by (6.13) follows from the assertion

$$r^{-2} \left| \int_{B_r} \langle F_A^+ \circ F_A^- \rangle(t) \right| \leq c_0 \int_{B_r} (|\nabla_A a|^2 + r^2 |\langle a, \tau a \rangle|^2) + c_0 r^{-1} \int_{B_r} |a| \, |\nabla_A a|$$

(6.17)

because of what was said in the last step about (6.15) and (6.16). The derivation of (6.17) is given momentarily. To set the notation, fix an oriented, orthonormal frame for $\Lambda^+$ over $B_r$ to be denoted now as $\{\omega^1, \omega^2, \omega^3\}$. Let $v_r$ denote the unit length vector field on $B_r - p$ that points outward from p along the geodesics through p. Use $v_r$ and $\{\omega^a\}_{a=1,2,3}$ to define orthonormal covectors fields $\{\hat{e}^1, \hat{e}^2, \hat{e}^3\}$ by the rule whereby $\hat{e}^a = \sqrt{2}\omega^a(v_r, \cdot)$. These 1-forms with the differential of the function dist(p,·) define an orthonormal basis for T*X over $B_r - p$.) It is important in what follows that the dual vector fields to $\{\hat{e}^a\}_{a=1,2,3}$ are tangent to the spheres of constant radius in $B_r$ with center at p.

Both of the covectors $E_A$ and $B_A$ are linear combinations of $\{\hat{e}^a\}_{a=1,2,3}$ (because they are obtained from $F_A$ and $*F_A$ by contracting with $v_r$). The coefficients of $B_A$ with respect to this (partial) basis for T*X are denoted in what follows by $\{B_{A,a}\}_{a=1,2,3}$. If $a \in \{1, 2, 3\}$, let $\nabla_{A,a}$ denote the directional covariant derivative on $B_r - p$ in the direction of the dual vector field to $\hat{e}^a$.

The argument for (6.17) starts with the observation that the coefficients of $F_A^-$ are linear combinations of those of $E_A - B_A$ and those of $F_A^+$ are linear combinations of those of $E_A + B_A$. This being the case, each coefficient of $F_A^-$ can be written as a linear combination of coefficients of $F_A^+$ and of $B_A$. This rewriting of $F_A^-$ writes the function $\langle F_A^+ \circ F_A^- \rangle(t)$ schematically as

$$\langle F_A^+ \circ F_A^- \rangle(t) = \mathfrak{X}^+ + r^2 \langle a, \tau^a B_{A,b} a] \rangle \mathfrak{b}_{ab}$$

(6.18)

with the norm of $\mathfrak{X}^+$ obeying $|\mathfrak{X}^+| \leq c_0 |F_A^+|^2$ and each $\mathfrak{b} \in \{\mathfrak{b}_{ab}\}_{1 \leq a,b \leq 3}$ obeying $|\mathfrak{b}| \leq c_0$ and $|\nabla \mathfrak{b}| \leq c_0 r^{-1}$. Of particular note: The term with $B_A$ in (6.18) can be written as

$$\tfrac{1}{2} r^2 \langle a, \tau^a F_{A,cd}, a \rangle \mathfrak{b}_{ab} \varepsilon_{bcd} ,$$

(6.19)

where $\{\varepsilon_{bcd}\}_{1 \leq b,c,d \leq 3}$ are the coefficients of the completely antisymmetric 3×3 tensor having $\varepsilon_{123}$ equal to 1.

The next step rewrites (6.18) using the fact that the product $F_{A,cd} a$ can be written as a commutator of covariant derivatives plus a term that is bounded by $c_0 |a|$. The point being that if $s$ is any given section of $\mathbb{S}_E$, then $F_{A,cd} s$ can be written as

$$F_{A,cd}\, s = \nabla_{A,c}\nabla_{A,d}\, s - \nabla_{A,d}\nabla_{A,c}\, s + t_{cd,e}\nabla_{A,e}s + \mathcal{R}_{c,d}\, s$$

(6.20)

with $\mathcal{R}_{c,d}$ and the functions $\{t_{cd,e}\}_{e \in \{1,2,3\}}$ having $c_0$ norm bounds. (The functions $\{t_{cd,e}\}_{c,d,e \in \{1,2,3\}}$ are defined as follows: Let $\{\hat{e}_b\}_{b=1,2,3}$ denote the dual vector fields to the basis vectors $\{\hat{e}^b\}_{b=1,2,3}$. Since these are tangent to the spheres of constant radius centered at p, they define an involutive system, which is to say that their commutators is in their linear span. The functions $\{t_{cd,e}\}_{c,d,e \in \{1,2,3\}}$ are obtained by writing the vector field commutator $[\hat{e}_c, \hat{e}_d]$ as $t_{bc,d}\hat{e}_d$.) The bound in (6.17) follows from (6.18) and (6.19) using (6.20) and an integration by parts. (Keep in mind with regards to the integration by parts that the derivatives in (6.20) are tangent to the spheres at constant distance from p.

Step 6: This step and Step 7 derive (6.5). To start, introduce S to denote the symmetric section of T*X⊗T*X with components $\{S_{\alpha\beta}\}_{\alpha,\beta \in \{1,2,3,4\}}$ (written with respect to a oriented, orthonormal frame on the whole of $B_r$) given by

$$S_{\alpha\beta} = \tfrac{1}{2}(\langle \nabla_{A,\alpha}a, \nabla_{A,\beta}a\rangle + \langle \nabla_{A,\beta}a, \nabla_{A,\alpha}a\rangle) - \tfrac{1}{2}\delta_{\alpha\beta}|\nabla_A a|^2 .$$

(6.21)

Commute derivatives to see that the components $\{(\nabla^\dagger S)_\beta\}_{\beta \in \{1,2,3,4\}}$ of the 1-form $\nabla^\dagger S$ are given by the formula

$$(\nabla^\dagger S)_\beta = \tfrac{1}{2}(\langle \nabla^\dagger \nabla a, \nabla_{A,\beta}a\rangle + \langle \nabla_{A,\beta}a, \nabla^\dagger \nabla a\rangle) + \tfrac{1}{2}F_{A,\alpha\beta}(\langle \nabla_{A,\alpha}a, a\rangle - \langle a, \nabla_{A,\alpha}a\rangle)$$
$$+ \tfrac{1}{2}(\langle \nabla_{A,\alpha}a, \mathcal{R}_{\alpha\beta}a\rangle + \langle \mathcal{R}_{\alpha\beta}a, \nabla_{A,\alpha}a\rangle$$

(6.22)

with $\{F_{A,\alpha\beta}\}_{\alpha,\beta \in \{1,2,3,4\}}$ denoting the components of $F_A$ with respect to the given basis. The terms from the set $\{\mathcal{R}_{ac\alpha\beta}\}_{a,c \in \{1,2,3\}, \alpha,\beta \in \{1,2,3,4\}}$ that appear in this equation are each linear combinations of components of the Riemann curvature tensor.

Use (2.9) to write the term $\langle \nabla^\dagger \nabla a, \nabla_{A,\beta}a\rangle$ that appears in (6.22) as

$$\tfrac{1}{4}r^2\nabla_\beta|\langle a, \tau a\rangle|^2 + \tfrac{1}{2}(\langle \nabla_{A,\beta}a, \mathfrak{R}a\rangle + \langle \mathfrak{R}a, \nabla_{A,\beta}a\rangle) .$$

(6.23)

The term that is denoted by $\mathfrak{r}$ in (6.5) is the sum of the terms with $\mathcal{R}$ in (6.22) and the terms with $\mathfrak{R}$ in (6.23).

Step 7: Use the identity in (2.14) to write the term with $F_A$ in (6.22) as

$$\tfrac{1}{2}F_{A,\alpha\beta}(\langle \nabla_{A,\alpha}a, a\rangle - \langle a, \nabla_{A,\alpha}a\rangle) = \tfrac{1}{2}r^{-2}F_{A,\alpha\beta}\nabla_\nu F_{A,\alpha\nu} .$$

(6.24)

The right hand side of this can be written in turn as

$$-\tfrac{1}{2} r^{-2} \nabla_\nu (\langle F_{A,(\cdot)\beta}, F_{A,(\cdot)\nu}\rangle - \tfrac{1}{2}\delta_{\beta\nu}|F_A|^2)$$

(6.25)

by using the equation $dF_A = 0$. The expression in (6.25) is the same as

$$- r^{-2} \nabla_\nu \langle F_A^+{}_{,\alpha\beta}, F_A^-{}_{,\alpha\nu}\rangle$$

(6.26)

because $\langle F_{A,(\cdot)\beta} F_{A,(\cdot)\nu}\rangle - \tfrac{1}{2}\delta_{\alpha\beta}|F_A|^2$ is equal to $2\langle F_A^+{}_{,(\cdot)\beta}, F_{A,(\cdot)\nu}\rangle$. The formula in (6.5) follows directly from (6.22), (6.23), (6.24)–(6.26) and the preceding definition of q.

**c) A local averaged version of N**

The second and third bullets of Lemma 6.1 are useful because these say in effect that if N(r) is small for a given r, and if $s \in [\tfrac{1}{2} r, r]$, then $N(s) \leq 4 N(r)$ and thus N(s) is not a great deal larger than N(r). As explained below, a consequence is that N can be replaced by a local average that is easier to work with. Looking ahead, the local average is preferred by virtue of the fact that its derivative can be written as an integral over a 4-dimensional domain in X. This is in contrast to the integrals that appear in Proposition 6.2, which are integrals over hypersurfaces.

The definition of this local average of N requires the a priori choice of a number $\mu \in (0, \tfrac{1}{100}]$. The corresponding average is denoted in what follows by $\underline{N}$ (the number $\mu$ is not notationally indicated.) This $\underline{N}$ is the function on $(0, c_0^{-1}]$ that is defined by the rule

$$r \to \underline{N}(r) = \tfrac{\pi}{2} \int_0^1 \sin(\pi\tau) N((1-\tau\mu) r) d\tau .$$

(6.27)

For the record, a change of variables writes $\underline{N}(r)$ as

$$\underline{N}(r) = \tfrac{\pi}{2\mu} \int_{(1-\mu)r}^{r} \tfrac{s}{r}\sin(\tfrac{\pi}{\mu}(1-\tfrac{s}{r})) N(s) \tfrac{ds}{s} .$$

(6.28)

Note that if N is constant on $[(1-\mu)r, r]$, then $\underline{N}(r)$ is N(r) because the integral of the function $\tau \to \sin(\pi\tau)$ on the interval [0, 1] is $\tfrac{2}{\pi}$. In general $\underline{N}(r)$ is a weighted average of the values N on the interval $[(1-\mu)r, r]$. (There is no special reason for using the sine function to define $\underline{N}$ via (6.28). Any function that vanishes at the end points of the interval $[(1-\mu)r, r]$ should suffice.)

The function $\underline{N}$ is a proxy for the function N when N is small because of the features that are described in the lemma that follows.

**Lemma 6.3**: *There exists $\kappa > 1$ with the following significance: Fix $r > 1$ and suppose that $(A, a)$ is solution to (1.5). Fix a point in $X$ so as to define the function $N$ and fix $\mu \in (0, \frac{1}{100})$ to define the function $\underline{N}$. Supposing that $r \in (0, \kappa^{-1})$, then*

- *If $N(r) \leq 1$, then $\underline{N}(r) \leq \frac{1}{(1-\mu)^5} N(r)$.*
- *If $N(r) \leq 1$ and if $N(r) \leq 4 N((1-\mu)r)$ then $N(r) \leq \frac{4}{(1-\mu)^5} \underline{N}(r)$*
- *If $\underline{N}(r) \leq 1$, then $N(\cdot) \leq 3 \underline{N}(r)$ on $[(1-\mu)r, (1-\frac{1}{2}\mu)r]$ and $N((1-\mu)r) \leq \frac{1}{(1-\mu)^5} \underline{N}(r)$.*
- *If $\underline{N}(r) \leq 1$, then $\underline{N}((1-\mu)r) \leq \frac{1}{(1-\mu)^5} \underline{N}(r)$.*

*Proof of Lemma 6.3*: The top bullet follows from the second bullet of Lemma 6.1 because the latter implies that

$$N(s) \leq \frac{1}{(1-\mu)^5} N(r)$$

(6.29)

if $N(r) \leq 1$ and if $s \in [(1-\mu)r, r]$. To prove the second bullet, first invoke (6.29) to see that $N((1-\mu)r) \leq \frac{1}{(1-\mu)^5}$. With this bound handy, the top bullet of Lemma 6.1 implies that

$$N(s) \geq (1-\mu)^5 N((1-\mu)r)$$

(6.30)

for $s \in [(1-\mu)r, r]$. As a consequence $\underline{N}(r) \geq (1-\mu)^5 N((1-\mu)r)$. This last bound with the assumption that $N((1-\mu)r) \geq \frac{1}{4} N(r)$ leads directly to the assertion of the second bullet in Lemma 6.3. The first inequality in the third bullet also follows from the second bullet of Lemma 6.1 since $N(r)$ must be less than $2\underline{N}(r)$ at some point on the interval $[(1-\frac{1}{2}\mu)r, r]$ if the integral in (6.28) is equal to $\underline{N}(r)$. The second inequality in the third bullet follows from (6.28) and Lemma 6.1 since $\frac{\pi}{2} \sin(\pi(\cdot))$ on $(0,1)$ is positive and its integral is equal to 1. The fourth bullet's inequality follows from (6.29) using the formula for $\underline{N}$ in (6.27).

**d) The derivative of $\underline{N}$**

The derivative of $\underline{N}$ can be written using the formula in (6.28) as

$$r \frac{d}{dr} \underline{N}(r) = \frac{\pi}{2\mu} \int_{(1-\mu)r}^{r} \frac{s}{r} \sin(\frac{\pi}{\mu}(1-\frac{s}{r}))(\frac{d}{ds} N(s)) ds$$

(6.31)

because $r \frac{d}{dr} f(\frac{s}{r}) = -s \frac{d}{ds} f(\frac{s}{r})$ for any differentiable function $f$. Now let $\mathcal{A}(r, \mu)$ denote the spherical shell centered at p with outer radius r and inner radius $e^{-\mu} r$. (This is the complement in $B_r$ of the concentric ball with radius $e^{-\mu} r$.) Proposition 6.2 can be invoked to write (6.31) in turn as

$$r \frac{d}{dr} \underline{N}(r) = \frac{\pi}{2\mu} \int_{\mathcal{A}(r,\mu)} \frac{s}{r} \sin(\frac{\pi}{\mu}(1-\frac{s}{r}))(\frac{2}{s^2 K(s)^2} |\nabla_{A,r} a - \frac{1}{s} N a|^2 + \frac{1}{r^2 s^3 K(s)^2}(|E_A|^2 + \langle E_A, B_A \rangle)) + \hat{q}$$

(6.32)

with $\hat{q}$ denoting a function on $(0, c_0^{-1}]$ with norm bounded by $c_\mu r^2 (1 + \underline{N})$.

The following proposition says more about the derivative of $\underline{N}$. It is the central result from this section.

**Proposition 6.4**: *There exists $\kappa > 1$ and, given $\mu \in (0, \frac{1}{100})$, there exists $\kappa_\mu > \kappa$; and they have the following significance: Fix $r > 1$ and suppose that $(A, a)$ is a solution to (1.5). Fix a point $p \in X$ so as to define the functions $K$ and $N$ as done in Section 3a. Use $\mu$ to define the function $\underline{N}$ as instructed in Section 6c. For $r \in (0, \kappa_\mu^{-1})$, define $\underline{N}_\diamond = \underline{N}(\frac{1}{(1-\mu)^2} r)$. If $r$ is such that $\underline{N}_\diamond < \kappa^{-1}$, then*

$$r \frac{d}{dr} \underline{N} > -\kappa_\mu \frac{1}{(rrK)^{1/\kappa}} \left( \underline{N}_\diamond + r^2 + (\sqrt{\underline{N}_\diamond} + r) (\frac{1}{r^2 r^2 K^2} \int_{B_{r/(1-\mu)}} |F_{\hat{A}}^-|^2 )^{1/2} \right) - \kappa_\mu r^2 (1 + \underline{N}).$$

By way of a parenthetical remark (for now), this inequality is useful because of the appearance of the factor $\frac{1}{(rrK)^{1/\kappa}}$ on its right hand side.

*Proof of Proposition 6.4*: The proof that follows has four parts.

*Part 1*: Reintroduce the connection $\hat{A}$ and the 1-form $\flat$ from Proposition 5.3. The 2-form $F_A$ when written using $\hat{A}$ and $\flat$ is $F_A = F_{\hat{A}} + d\flat$. The contraction of $F_{\hat{A}}$ with the unit length vector field on $B_r - p$ that is tangent (and outward pointing) to the geodesics from $p$ is denoted by $E_{\hat{A}}$; and the corresponding contraction with $d\flat$ is denoted by $E_{d\flat}$. Let $B_{\hat{A}}$ denote the contraction of this vector field with $*F_{\hat{A}}$. Since $F_{\hat{A}}^+ = 0$, the 1-forms $E_{\hat{A}}$ and $B_{\hat{A}}$ obey the identity $E_{\hat{A}} + B_{\hat{A}} = 0$. This implies that the function $|E_A|^2 + \langle E_A, B_A \rangle$ that appears on the right hand side of (6.32) can be written using $\hat{A}$ and $\flat$ without terms that are quadratic functions of the components of $F_{\hat{A}}$:

$$|E_A|^2 + \langle E_A, B_A \rangle = |E_{d\flat}|^2 + + 2\langle E_{d\flat}, E_{\hat{A}} \rangle - \frac{1}{2} *(d\flat \wedge d\flat) - *(d\flat \wedge F_{\hat{A}}) .$$

(6.33)

*Part 2*: This part of the proof rewrites (6.33) as the sum of the manifestly positive term $|E_{d\mathfrak{b}}|^2$ plus a term that is a total derivative plus a term that contains a component of $\mathfrak{b}$ (with no derivatives). To do this, use the fact that $F_{\hat{A}}$ is a closed 2-form to first write

$$|E_A|^2 + \langle E_A, B_A \rangle = |E_{d\mathfrak{b}}|^2 + 2\langle E_{d\mathfrak{b}}, E_{\hat{A}} \rangle - *d(\tfrac{1}{2} \mathfrak{b} \wedge d\mathfrak{b} + \mathfrak{b} \wedge F_{\hat{A}}) \ .$$

(6.34)

The right most term on the right hand side of (6.34) is a total derivative. To rewrite the middle term on the right hand side of (6.34), introduce $\mathcal{L}_r$ to denote the Lie derivative on $B_r-p$ along the unit length vector field that is tangent (and outward pointing) to the geodesics through p; and let $\mathfrak{b}_r$ denote the contraction of $\mathfrak{b}$ with this vector field. Since $E_{d\mathfrak{b}} = \mathcal{L}_r\mathfrak{b} - d\mathfrak{b}_r$, the term $2\langle E_{d\mathfrak{b}}, E_{\hat{A}} \rangle$ on the right hand side of (6.34) is $2\langle (\mathcal{L}_r\mathfrak{b} - d\mathfrak{b}_r), E_{\hat{A}} \rangle$. The product rule for derivatives can now be used to rewrite this term as:

$$2\langle E_{d\mathfrak{b}}, E_{\hat{A}} \rangle = 2\nabla_r \langle \mathfrak{b}, E_{\hat{A}} \rangle - *d*(\mathfrak{b}_r E_{\hat{A}}) + \tfrac{1}{r}\mathcal{R}_0(\mathfrak{b}, F_{\hat{A}}) + \mathcal{R}_1(\mathfrak{b}, \nabla F_{\hat{A}})$$

(6.35)

where $\mathcal{R}_0$ and $\mathcal{R}_1$ are tensors on $B_r$ with norms bounded by $c_0$. The left most two terms on the right hand side of (6.35) are total derivatives; and the right most two terms are linear in the components of $\mathfrak{b}$. To summarize, (6.34) and (6.35) have rewritten (6.33) as

$$|E_A|^2 + \langle E_A, B_A \rangle = |E_{d\mathfrak{b}}|^2 - *d(\tfrac{1}{2} \mathfrak{b} \wedge d\mathfrak{b} + \mathfrak{b} \wedge F_{\hat{A}} - *\mathfrak{b}_r E_{\hat{A}}) + 2\nabla_r \langle \mathfrak{b}, E_{\hat{A}} \rangle$$
$$+ \tfrac{1}{r}\mathcal{R}_0(\mathfrak{b}, F_{\hat{A}}) + \mathcal{R}_1(\mathfrak{b}, \nabla F_{\hat{A}})$$

(6.36)

which has the desired form.

*Part 3*: Replace the left hand side of (6.36) where it appears on the right hand side of (6.32) with the sum on the right hand side of (6.36). Since the formula for the derivative of $\underline{N}$ involves a *volume* integral (as opposed to the surface integral in Proposition 6.3's formula for $\tfrac{d}{dr}$ N), an integrate by parts can be used to write each term from the right hand side of (6.36) but the manifestly positive $|E_{d\mathfrak{b}}|^2$ as the integral of a sum of terms that each have a component of $\mathfrak{b}$ as a factor. In particular, the result of this rewriting leads directly to the following inequality:

$$r\tfrac{d}{dr}\underline{N}(r) \geq -c_\mu \int_{\mathcal{A}(r,\mu)} \tfrac{1}{r^2 s^3 K(s)^2} (1+\underline{N}(s))|\mathfrak{b}|(|d\mathfrak{b}| + |F_{\hat{A}}| + s|\nabla F_{\hat{A}}|) - c_\mu r^2(1+\underline{N}) \ .$$

(6.37)

Given the assumption that $\underline{N}(\tfrac{1}{(1-\mu)^2} r) \leq 1$, it follows from the third and fourth bullets of Lemma 6.3 (with r replaced by $\tfrac{1}{(1-\mu)^2} r$) and from Lemma 6.1 that $N(s) \leq c_0$ on $[(1-\mu)r, \tfrac{1}{1-\mu}r]$. This implies (via (3.9)) that K obeys

- $K(s) \geq c_0^{-1} K(r)$ *for* $s \in [(1-\mu)r, r]$.
- $K(s) \leq c_0 K(r)$ *for* $s \in [r, \frac{1}{1-\mu}r]$.

(6.38)

The bound from the top bullet and the $c_0$ bound for $N$ on the interval $[(1-\mu)r, r]$ can be used in (6.37) to see that

$$r \frac{d}{dr} \underline{N}(r) \geq - c_\mu \frac{1}{r^2 r^3 K^2} \int_{\mathcal{A}(r,\mu)} |\mathfrak{b}|(|d\mathfrak{b}| + |F_{\hat{A}}| + s|\nabla F_{\hat{A}}|) - c_\mu r^2 (1 + \underline{N}) \ .$$

(6.39)

Hölder's inequality can now be brought to bear on (6.39) with the result being

$$r \frac{d}{dr} \underline{N}(r) \geq - c_\mu \frac{1}{r^2 r^3 K^2} \left( \int_{B_r} |\mathfrak{b}|^2 \right)^{1/2} \left( \int_{B_r} (|d\mathfrak{b}|^2 + |F_{\hat{A}}|^2 + s^2 |\nabla F_{\hat{A}}|^2) \right)^{1/2} - c_\mu r^2 (1 + \underline{N})$$

(6.40)

With (6.40) understood, invoke Item b) of the third bullet of Proposition 5.3 to obtain a bound for the $L^2$ norms of $\mathfrak{b}$ on $B_r$. What is denoted by $\tilde{K}$ in Item b) of Proposition 5.3's third bullet can be replaced by $c_0 K$ by virtue of what is said in the second bullet of (6.38); and the appearances of $N$ and $\tilde{N}$ can be replaced by $c_0 \underline{N}_\diamond$ because of what is said in Lemma 6.3. These replacements lead from Item b) of the third bullet in Proposition 5.3 to the bound

$$\int_{B_r} |\mathfrak{b}|^2 \leq c_\mu (rr K)^{2-1/\kappa} (\underline{N}_\diamond + r^2) \, r^2 \ .$$

(6.41)

Meanwhile, the fourth bullet of Proposition 5.3 with r replaced by $\frac{1}{1-\mu} r$ can be used to bound the integral of $|d\mathfrak{b}|^2$ that appears in (6.40) by $c_\mu (rr K)^2 (\underline{N}_\diamond + r^2)$. This again invokes Lemma 6.3 to replace the appearance of $N(\frac{1}{1-\mu} r)$ by $c_\mu \underline{N}_\diamond$ and it uses the second bullet of (6.38) to replace $K(\frac{1}{1-\mu} r)$ by $c_\mu K(r)$.

What was said in the preceding paragraph about the $\mathfrak{b}$ and $d\mathfrak{b}$ integrals in (6.41) take (6.40) to the bound

$$r \frac{d}{dr} \underline{N}(r) \geq - c_\mu \frac{1}{(rr K)^{1/c}} (\underline{N}_\diamond + r^2 + (\sqrt{\underline{N}_\diamond} + r)(\frac{1}{r^2 r^2 K^2} \int_{B_r} |F_{\hat{A}}|^2 + \frac{1}{r^2 K^2} \int_{B_r} |\nabla F_{\hat{A}}|^2)^{1/2}) - c_\mu r^2 (1 + \underline{N}) \ .$$

(6.42)

This last inequality is almost what is asserted by Proposition 6.4; the appearance of the term with $\nabla F_{\hat{A}}$ is the only anomaly.

*Part 4*: As explained in this last part of the proof, the $L^2$ norm of $\nabla F_{\hat{A}}$ on the ball $B_r$ is bounded by $c_\mu \frac{1}{r}$ times the $L^2$ norm of $F_{\hat{A}}$ on the concentric ball with radius $\frac{1}{1-\mu} r$. Granted this bound, then (6.42) leads directly to the assertion made by Proposition 6.4.

Meanwhile, this asserted bound for the $L^2$ norm of $\nabla F_{\hat{A}}$ follows from the fact that $F_{\hat{A}}$ is a closed, anti-self dual 2-form. Thus, it is harmonic and so it obeys a second order equation of the form

$$\nabla^\dagger \nabla F_{\hat{A}} + \mathcal{R} F_{\hat{A}} = 0 \tag{6.43}$$

with $\mathcal{R}$ being an endomorphism that is linear in the components of the Riemann curvature. In particular, $|\mathcal{R}| \leq c_0$. The equation in (6.43) implies in turn inequality

$$\tfrac{1}{2} d^\dagger d |F_{\hat{A}}|^2 + |\nabla F_{\hat{A}}|^2 \leq c_0 |F_{\hat{A}}|^2 . \tag{6.44}$$

Keeping (6.44) on hold for the moment, construct from $\chi$ a non-negative 'bump' function with compact support in the radius $\frac{1}{1-\mu} r$ ball centered at p, equal to 1 on $B_r$ and whose differential has norm bounded by $c_0 \frac{1}{\mu r}$. Multiply both sides of (6.44) by this function and then integrate over the radius $\frac{1}{1-\mu} r$ ball centered at p. A suitable integration by parts leads from the resulting integral identity to the asserted $L^2$ bound for $\nabla F_{\hat{A}}$.

## 7. Power law bounds for κ

The propositions and lemmas from Sections 5 and 6 will be used to derive a power law growth bound for the function κ above any given (small) value of r. This is a bound of the form $\kappa(s) \leq s^\gamma$ for $s \geq r$ but less than some fixed number ρ with γ and ρ determined a priori by a lower bound for $N(r)$. The precise result is summarized by the next proposition. To set the notation, introduce $c_\Diamond$ to denote $10^6$ times the maximum of the versions of κ that appear in Lemmas 6.1 and 6.3 and in Propositions 5.1, 6.2 and 6.4. Given a positive number s which is small enough so that the number $N(s)$ is defined (for example, $s \leq c_\Diamond^{-1}$), let $L(s)$ denote the minimum of $N(s)$ and $\frac{1}{10,000} c_\Diamond^{-1}$.

**Proposition 7.1**: *Given $\mu \in (0, \frac{1}{100})$, there exists $\kappa_\mu > 100$, and given also $\varepsilon \in (0, \kappa_\mu^{-1})$, there exists $\rho_{\mu,\varepsilon} \in (0, \kappa_\mu^{-1})$ with the following significance: Fix $r > 1$ and suppose that the pair $(A, a)$ is a solution to (1.5). Fix $p \in X$ to define the functions $N$ and $\kappa$, and then $L$. If $r \in (0, \rho_{\mu,\varepsilon}]$ is such that $N(r) \geq \varepsilon$, then $\kappa(\cdot)$ on the interval $[r, \rho_{\mu,\varepsilon}]$ obeys $\kappa(s) \leq s^{L(r)^3/\kappa_\mu}$.*

To be sure, Proposition 7.1 does not say that the condition $N(r) \geq \varepsilon$ implies that $N(s)$ for $s \geq r$ is bounded away from zero by a fixed fraction of the smaller of 1 or the square of N(r) (which is an event that would lead to the second bullet's bound via (3.9).) However, the proposition does say, in effect, that intervals in $(r, c_0^{-1}]$ where N is very small are always offset by intervals that are not so far away (when measured using $\frac{ds}{s}$) where N is relatively large. The following lemma describes an elementary example where K has a power law upper bound on a subset with N not (necessarily) relatively large everywhere on that subset.

**Lemma 7.2**: *Let $\alpha$, $\beta$ and $\gamma$ denote positive numbers with $\alpha$, $\beta$ being greater than 1. Suppose that $r \in (0, \kappa^{-1}]$ and $r' > \alpha r$ are such that $K(r) \leq (\frac{r}{r'})^\gamma K(r')$. If $s \in [\beta^{-1} r, r]$ then $K(s) \leq (\frac{s}{r'})^{\gamma \ln(\alpha)/\ln(\alpha\beta)} K(r')$. If, in addition, $K(s) \leq (\frac{s}{r'})^\gamma K(r')$ for all $s \in [r, r']$, then $K(s) \leq (\frac{s}{r'})^{\gamma \ln(\alpha)/\ln(\alpha\beta)} K(r')$ for all $s \in [\beta^{-1} r, r']$.*

*Proof of Lemma 7.2*: The second assertion follows from the first. The first assertion follows from the bound $K(r) \leq (\frac{r}{r'})^\gamma K(r')$ because $K(s) < K(r)$ for $s < r$, and because $\frac{r}{r'} \leq (\frac{s}{r'})^{\ln(\alpha)/\ln(\alpha\beta)}$ when $r' \geq \alpha r$ and $s \geq \beta^{-1} r$.

This lemma is used many times in the upcoming proof of Proposition 7.1.

Proposition 7.1 is seen momentarily to follow from the next two propositions which assert 'local' power law bounds for K.

**Proposition 7.3**: *Given $\mu \in (0, \frac{1}{100})$ there exists $\kappa_\mu > 100$, and given also $\varepsilon \in (0, \kappa_\mu^{-1})$, there exists $\rho_{\mu,\varepsilon} \in (0, \kappa_\mu^{-1})$ with the following significance: Suppose that $r > 1$ and that $(A, a)$ is a solution to $r$'s version of (1.5). Fix $p \in X$ and fix an input value $r \in (0, \rho_{\mu,\varepsilon})$ with $N(r) \geq \varepsilon$ and such that $rrK(r) \geq \kappa_\mu$. Given this input, there are numbers $s_0$ and $s_1$ (the output) with the properties listed below.*
- *$s_0 > 2r$ and $s_1 \geq s_0$, but neither $s_0$ nor $s_1$ is greater $c_\diamond^{-1}$.*
- *Let $s_*$ denote $s_0$ or $s_1$. If $s_* \leq \rho_{\mu,\varepsilon}$, then either $N(s_*) > \frac{1}{10,000} c_\diamond^{-1}$ or $N(s_*) \geq 2N(r)$.*
- *$K(s) \leq (\frac{s}{s_1})^{L(r)/\kappa_\mu} K(s_1)$ for all $s \in [r, s_0]$.*

Notice that $s_0$ and $s_1$ have different roles in the third bullet of the proposition: The point $s_0$ is the upper bound for the set where the inequality holds, but it is $s_1$ that appears in the actual inequality.

Proposition 7.3 has a lower bound requirement on its input r because of the requirement that $rrK(r) \geq \kappa_\mu$. The next proposition considers the case where this constraint is violated. To set the notation, and supposing that $\mu$ has been specified, introduce $\kappa_{\diamond\mu}$ to denote the version of $\kappa_\mu$ that appears in Proposition 7.3.

**Proposition 7.4**: *Given $\mu \in (0, \frac{1}{100})$ there exists $\kappa_\mu > 100$, and given also $\varepsilon \in (0, \kappa_\mu^{-1})$, there exists $\rho_{\mu,\varepsilon} \in (0, \kappa_\mu^{-1})$ with the following significance: Suppose that $r > 1$ and that $(A, a)$ is a solution to $r$'s version of (1.5). Fix $p \in X$ and fix an input value $r \in (0, \rho_{\mu,\varepsilon})$ with $N(r) \geq \varepsilon$ and such that $rrK(r) \leq \kappa_{\diamond\mu}$. Given this input, there are numbers $s_0$ and $s_1$ (the output) with the properties listed below.*

- $s_0 > 2r$ *and* $s_1 \geq s_0$, *but neither* $s_0$ *nor* $s_1$ *is greater* $c_\diamond^{-1}$.
- *Let* $s_*$ *denote* $s_0$ *or* $s_1$. *If* $s_* \leq \rho_{\mu,\varepsilon}$, *then one of the following conditions holds*
  a) $N(s_*) > \frac{1}{10,000} c_\diamond^{-1}$.
  b) $N(s_*) \geq 2N(r)$.
  c) $rs_* K(s_*) \geq \kappa_{\diamond\mu}$ *and* $N(s_*) \geq \kappa_\mu^{-1} N(r)^3$.
- $K(s) \leq (\frac{s}{s_1})^{L(r)^3/\kappa_\mu} K(s_1)$ *for all* $s \in [r, s_0]$.

Proposition 7.3 is proved separately in Section 8 because there are many parts to the proof. Proposition 7.4 is proved separately in Section 9.

***Proof of Proposition 7.1***: The proof has three parts. Parts 1 and 2 assume that the input $r$ for Proposition 7.1 obeys the Proposition 7.3 constraint $rrK(r) \geq \kappa_{\diamond\mu}$; and Part 2 considers the case where this condition is violated. (The argument when this condition is violated requires just one extra comment which is given at the end.)

*Part 1*: As will be evident momentarily, Proposition 7.3 is the central element of an iterative algorithm that will prove Proposition 7.1 when $rrK(r) \geq \kappa_{\diamond\mu}$. To first get a sense of the algorithm, suppose for the moment in this Part 1 that Proposition 7.3 were slightly stronger in the following way: For any input number (to use for Proposition 7.3's version of $r$), the proposition returns $s_0$ and $s_1$ with $s_0$ being *equal* to $s_1$. Assuming this stronger assertion, one could invoke this hypothetical stronger version of Proposition 7.3 with the input being the version of $r$ in Proposition 7.1 to obtain a point $s_1$ with $s_1 > 2r$ such that

$$K(s) \leq (\tfrac{s}{s_1})^{L(r)/\kappa_{\diamond\mu}} K(s_1) \text{ for all } s \in [r, s_1] .$$

(7.1)

This is the first step of the iteration. If $s_1 \geq \rho_{\mu,\varepsilon}$, then one can stop. If $s_1 < \rho_{\mu,\varepsilon}$, then Proposition 7.3 can be invoked again with $s_1$ used for the input (use $s_1$ in lieu of what the proposition calls $r$) to get an output number $s_{11} > 2s_1$ (which is greater than 4 times the number $r$ from Proposition 7.1) such that

$$K(s) \leq \left(\frac{s}{s_{11}}\right)^{L(s_1)/\kappa_{\diamond\mu}} K(s_{11}) \quad \textit{for all} \ \ s \in [s_1, s_{11}] \ .$$
(7.2)

Because $N(s_1) \geq N(r)$ (it is actually at least twice $N(r)$ in the cases when $N(r) < \frac{1}{10,000} c_\diamond^{-1}$), this inequality and (7.1) imply that

$$K(s) \leq \left(\frac{s}{s_{11}}\right)^{L(r)/\kappa_{\diamond\mu}} K(s_1) \quad \textit{for all } s \in [r, s_{11}] \ .$$
(7.3)

Indeed, this follows from $s \in [s_1, s_{11}]$ from (7.2) and it is derived for $s \in [r, s_1]$ by replacing $K(s_1)$ in (8.1) with the $s = s_1$ version of the right hand side of (7.2). One could then repeat all of this with $s_{11}$ used as input to Proposition 7.3 (assuming that $s_{11} < \rho_{\mu,\varepsilon}$) to obtain $s_{111} > 2 s_{11}$ (which is greater than $8r$) so that (7.3) holds with $s_{11}$ replaced by $s_{111}$; and so on leading finally to some number $s_\ddagger = s_{11\cdots1} \geq \rho_{\mu,\varepsilon}$ such that $K(s) \leq \left(\frac{s}{s_\ddagger}\right)^{L(r)/\kappa_{\diamond\mu}} K(s_\ddagger)$ for all $s$ between $r$ and $s_\ddagger$. This is what needs to be proved to verify Proposition 7.1

*Part 2*: This part proves the $rr K(r) \geq \kappa_{\diamond\mu}$ case of Proposition 7.1. As written, Proposition 7.3 takes as input a number it denotes by $r$ and it gives back two numbers $s_0 > 2r$ and $s_1 \geq s_0$ with no guarantee that these are the same. In this case, there is, in any event, an iteration of Proposition 7.3 to verify Proposition 7.1; but it now proceeds along a branching binary tree. To explain, the first application of Proposition 7.3 uses the number $r$ from Proposition 7.1 as its input (the Proposition 7.3 version of $r$), and it outputs numbers $\{s_0, s_1\}$. These are such that $s_0 \geq 2r$ and $s_1 \geq s_0$ and such that

$$K(s) \leq \left(\frac{s}{s_1}\right)^{L(r)/\kappa_{\diamond\mu}} K(s_1) \quad \textit{for all } s \in [r, s_0] \ .$$
(7.4)

Note in particular that this bound for $K(s)$ in the case when $s_1 > s_0$ is missing the values of $s$ between $s_0$ and $s_1$. In any event, one returns to Proposition 7.3 (when $s_0 < \rho_{\mu,\varepsilon}$) with inputs $s_0$ and then $s_1$ (they are used in lieu of what Proposition 7.3 denotes by $r$) to generate respective pairs $(s_{00}, s_{01})$ and $(s_{10}, s_{11})$; and then one returns again with these four inputs, and so on to produce a binary tree labled set $\{s_{ab\cdots c}\}$ in I with labels $a, b, \cdots, c$ from $\{0, 1\}$. Any such binary digit labeled point (denoted here by $s_*$) has $L(s_*) \geq L(r)$; and its progeny (denoted by $s_{*0}$ and $s_{*1}$) obey $s_{*1} \geq s_{*0} \geq 2 s_*$, and they are such that

$$K(s) \leq \left(\frac{s}{s_{*1}}\right)^{L(s_{*1})/\kappa_{\diamond\mu}} K(s_{*1}) \quad \textit{for all } s \in [s_*, s_{*0}] \ .$$
(7.5)

This interation ends along any path in the binary tree starting from $r$ when the binary digit label $*$ is such that $s_* \geq \rho_{\mu,\varepsilon}$. The set so constructed is finite because the respective size of $s_{*0}$ and $s_{*1}$ is 2 or more times that of $s_*$.

The construction is such that each point $s \in [r, \rho_{\mu\varepsilon}]$ is between some $s_*$ and a corresponding $s_{*0}$, even those points that end with the binary digit 1. This fact with (7.5) leads directly to Proposition 7.1's assertion because no $L(s_*)$ is smaller than $L(r)$.

*Part 3*: The preceding analysis assumed at the outset that the input r obeyed the Proposition 7.3 constraint $rrK(r) \geq \kappa_{\Diamond\mu}$. Let $r_\Diamond$ denote the value of r where $rrK(r) = \kappa_{\Diamond\mu}$. (There is a unique such r because the function $r \to rK(r)$ is increasing.) If $r < r_\Diamond$, then the iteration in Part 2 can be repeated almost verbatim using what is said by Proposition 7.4 instead of Proposition 7.3 (and with $\kappa_{\Diamond\mu}$ replaced by Proposition 7.4's version of $\kappa_\mu$) to go from r to $s_0$ and $s_1$; and then to the collection $\{s_{00}, s_{01}, s_{10}, s_{11}\}$; and, in general, from any given $s_*$ (with $*$ a binary number) to $s_{*0}$ and $s_{*1}$. But, in this case, the iteration stops on a given branch of the binary tree when $s_* \geq r_\Diamond$. In this event, the iteration in Part 2 is continued using Proposition 7.3 with input being this value of $s_*$ and its corresponding $N(s_*)$. By virtue of Item c of the second bullet of Proposition 7.4, this $N(s_*)$ is greater than the minimum of $c_\mu^{-1} N(r)^3$ and $c_\Diamond^{-1}$. Therefore, the Proposition 7.3 iteration has the starting value of N being $c_\mu^{-1} N(r)^3$ as opposed to being $N(r)$. Otherwise, this iteration is just as in Part 2.

## 8. Proof of Proposition 7.3

The proof of Proposition 7.3 needs various preliminary observations and lemmas, and so it is deferred to the end of this section (Section 8e). The intervening subsections state and prove these preliminary assertions. The results in Sections 5 and 6 are used for the most part to prove these preliminary items.

**a) The sets I and 𝐼.**

Let $c_\Diamond$ denote $10^6$ times the maximum of the versions of $\kappa$ that appear in Lemmas 6.1 and 6.3 and Propositions 5.1 and 6.4. Fix $r_1 < \frac{1}{100} c_\Diamond^{-1}$ for the moment and write the interval $[r_1, c_\Diamond^{-1}]$ as $I \cup I'$ with I distinguished as follows:

$$\text{If } r \in I, \text{ then } N(r) > \tfrac{1}{10,000} c_\Diamond^{-1} .$$

(8.1)

Since this is an open condition, the set I can be written as a disjoint union of open intervals. If $(s_1, s_2)$ is any such interval, and $s \in [s_1, s_2]$, it follows from (3.9) that

$$K(s) \leq \left(\tfrac{s}{s_2}\right)^{1/(10,000 c_\Diamond)} K(s_2) .$$

(8.2)

Meanwhile, the set $I'$ is contained in the open set $𝐼$ that is defined by the rule

$$\text{If } r \in \mathcal{I}, \text{ then } N(r) < c_0^{-1}.$$

(8.3)

Since $\mathcal{I}$ is an open set, it is a disjoint union of intervals. Let $(s_0, s_1)$ denote one of these intervals. If $(s, s´) \in (s_0, s_1)$, then it follows from (3.9) that

$$K(s) \geq (\tfrac{s}{s´})^{1/c_0} K(s´),$$

(8.4)

which is the reverse direction of the inequality in (8.2).

The points not in $\mathcal{I}$ can be used to illustrate the phenomena in Lemma 7.1 as follows: Suppose that $s_1$ is not in $\mathcal{I}$ (so that $N(s_1) \geq c_0^{-1}$). Then it is a consequence of Lemma 6.1 that there exists $s_2 \geq 8s_1$ such that $N(s) \geq \tfrac{1}{100} c_0^{-1}$ on the whole interval $[s_1, s_2]$. This implies that $K(s) \leq (\tfrac{s}{s_2})^{1/(100c_0)} K(s_2)$ for $s \in [s_1, s_2]$; and since $s_2 \geq 8s_1$, Lemma 7.2 guarantees that

$$K(s) \leq (\tfrac{s}{s_2})^{1/(200c_0)} K(s_2) \quad \text{for all } s \in [\tfrac{1}{8} s_1, s_2]$$

(8.5)

For example, if $(s_0, s_1)$ is a component of the set $\mathcal{I}$, then $s_1$ is not in $\mathcal{I}$ and so (8.5) holds; and it gives an opposite inequality to (8.4) near the $s_1$ end of the component $(s_0, s_1)$.

The bounds in (8.2) and (8.5) are formally summarized by the following lemma.

**Lemma 8.1**: *There exists $\kappa > 1$ with the following significance: Fix $r > 1$ and suppose that $(A, a)$ is a solution to (1.5). Fix $p \in X$ so as to define the functions $K$ and $N$, and the sets $I$ and $\mathcal{I}$.*

- *If $(s_1, s_2)$ is a component of $I$ and if $s \in [s_1, s_2]$, then $K(s) \leq (\tfrac{s}{s_2})^{1/(10,000 c_0)} K(s_2)$.*
- *If $s_1 \notin \mathcal{I}$, then there exists $s_2 \geq 8s_1$ such that $[s_1, s_2) \subset I$. In addition, if $s \in [\tfrac{1}{8} s_1, s_2]$, then $K(s) \leq (\tfrac{s}{s_2})^{1/(200 c_0)} K(s_2)$.*

Not much more will be said directly about I. Most of what is said next is dealing with $\mathcal{I}$.

**b) The set $\mathcal{I}_\mu$**

Fix $\mu \in (0, \tfrac{1}{100})$ and use it to define the function $\underline{N}$ as instructed in Section 6c. Having done so, then define subset $\mathcal{I}_\mu \subset \mathcal{I}$ by the following rule: A point r from $\mathcal{I}$ is in $\mathcal{I}_\mu$ when $rr K(r) \geq 1$ and when the following three conditions are met:

- *The interval $[r, 4r]$ is entirely in $\mathcal{I}$*

- $\underline{N}(\frac{1}{(1-\mu)^2} r) < 10^4 \underline{N}(r)$.

- $\int_{B_{r/(1-\mu)}} |F_A|^2 < c_\diamond^8 \, r^2 r^2 K(r)^2 \underline{N}(r)$ .

(8.6)

If $r \in \mathcal{I}_\mu$, then $\underline{N}(r)$ is small (it is less than $c_\diamond^{-1}$) on the interval $[r, 4r]$ since the interior of this whole interval is in $\mathcal{I}$. This implies, first that $\underline{N}(\frac{1}{(1-\mu)^2} r) \le z^{-1}$ with $z$ denote the version of $\kappa$ that appears in Proposition 6.4. (Remember that $\underline{N}(\frac{1}{(1-\mu)^2} r)$ is denoted by $\underline{N}_\diamond$ in Proposition 6.4.) The bound $\underline{N}(\frac{1}{(1-\mu)^2} r) \le z^{-1}$ puts Proposition 6.4 in play; and, as explained momentarily, this with the second and third bullets of (8.6) leads to an inequality for the derivative of $\underline{N}$:

$$r \frac{d}{dr} \underline{N} > -x \frac{1}{(rrK(r))^{1/z}} \underline{N} - z_\mu r^2 \; ,$$

(8.7)

with $z_\mu$ denoting the version of $\kappa_\mu$ from Proposition 6.1 and with $x$ denoting $c_0 z_\mu c_\diamond^8$. With regards to the derivation of (8.7): A direct quote of Proposition 6.4 gives the analog of (8.1) where the factor of $\underline{N}$ on the right hand side is replaced by $\sqrt{\underline{N}_\diamond}$ times the $L^2$ norm of $F_A^-$ on the ball $B_{r/(1-\mu)}$. The latter can be replaced by $c_\mu \sqrt{\underline{N}}$ by appeal to the bound in the third bullet of (8.6). Meanwhile, the $\sqrt{\underline{N}}$ factor is no greater than $c_0 \sqrt{\underline{N}_\diamond}$ because the bound $\underline{N}(\frac{1}{(1-\mu)^2} r) \le z^{-1}$ implies (with Lemmas 6.1 and 6.3) that $\underline{N}(r)$ is no greater than $100 \underline{N}(\frac{1}{(1-\mu)^2} r)$. Thus, (8.7) holds with $\underline{N}$ on the right replaced by $\underline{N}_\diamond$. But, since $\underline{N}_\diamond$ is $\underline{N}(\frac{1}{(1-\mu)^2} r)$, it is no greater than $100 \underline{N}(r)$ by appeal to the second bullet in (8.6)

With (8.7) understood, let $Q(r) = \underline{N}(r) + z_\mu r^2$. This $Q$ obeys (by virtue of (8.7))

$$r \frac{d}{dr} Q > -x \frac{1}{(rrK(r))^{1/z}} Q \; .$$

(8.8)

This last equation implies, in turn, the weaker assertion

$$r \frac{d}{dr} Q > -x \frac{1}{(rrK(r))^{1/z}} (1 + \underline{N}) Q$$

(8.9)

because $\underline{N}$ is positive. Meanwhile, (8.9) has the happy feature that it can be written as an inequality between two derivatives:

$$\frac{d}{dr} \ln(Q) > zx \left( \frac{d}{dr} \left( \frac{1}{(rrK(r))^{1/z}} \right) \right) .$$

(8.10)

And, this inequality can be integrated: If the whole of an interval $(s_0, s)$ is in $\mathcal{I}_\mu$, then

$$\underline{N}(s) + z_\mu s^2 \geq \frac{\eta(s)}{\eta(s_0)} (\underline{N}(s_0) + z_\mu s_0^2)$$

(8.11)

with $\eta(\cdot)$ denoting the function on $(0, c_\diamond^{-1})$ that is defined by the formula

$$\eta(s) = \exp(zx \frac{1}{(rs\kappa(s))^{1/z}}) \,.$$

(8.12)

The function $\eta$ is a decreasing function of $s$; and it is bounded by $e^{zx}$ when $rs\kappa(s) \geq 1$.

The following lemma summarizes and applies the preceding discussion.

**Lemma 8.2**: *There exists $\kappa > 1$, and given $\mu \in (0, \frac{1}{100})$, there exists $\kappa_\mu > \kappa$; and these numbers have the following significance: Fix $r > 1$ and suppose that $(A, a)$ is a solution to (1.5). Fix $p \in X$ to define the functions $\kappa$ and $N$, and the set $\mathcal{I}$. Then, fix $\mu \in (0, \frac{1}{100})$ to define the set $\mathcal{I}_\mu$. Suppose that $s_0 \in (0, \kappa_\mu^{-1})$ and $s_1 \in (s_0, \kappa_\mu^{-1})$ are such that the whole interval $(s_0, s_1)$ is in $\mathcal{I}_\mu$. Then the assertions below are obeyed when $s \in [s_0, s_1]$.*

- $\underline{N}(s) \geq (1 - \mu)^{10} \frac{\eta(s)}{\eta(s_0)} \underline{N}(s_0) - z_\mu s^2$.
- $\kappa(s) \leq e^{\kappa_\mu (s_1 - s) s_1} (\frac{s}{s_1})^{\delta(s_0)} \kappa(s_1)$ *with* $\delta(s_0) = (1 - \mu)^{10} \eta(s_0)^{-1} \underline{N}(s_0)$.

***Proof of Lemma 8.2***: If $(s_0, s_1)$ is in $\mathcal{I}_\mu$, then (8.11) holds for all $s \in [s_0, s_1]$. From this, the bound in the first bullet of Lemma 6.3 implies the top bullet of Lemma 8.2. The second bullet then follows from the first using (3.9).

### c) The sets $\mathcal{I}_\mu$ and $\mathcal{K}$

Supposing that $s_0 \in \mathcal{I}_\mu$, there will be a largest number in $(s_0, c_\diamond^{-1}]$ such that the whole of the open interval between $s_0$ and this number is in $\mathcal{I}_\mu$. Denote this number by $s_1$. Supposing that $s_1$ is not $c_\diamond^{-1}$ (the upper endpoint), then $s_1$ is described by one (or more) of the following:

- *There exists $s \in [s_1, 4s_1]$ such that $\underline{N}(s) = c_\diamond^{-1}$.*
- $\underline{N}(\frac{1}{(1-\mu)^2} s_1) = 10^4 \underline{N}(s_1)$.
- $\int_{B_{s_1/(1-\mu)}} |F_A|^2 = c_\diamond^8 r^2 s_1^2 \kappa(s_1)^2 \underline{N}(s_1)$.

(8.13)

The top bullet condition in (8.13) implies that $4s_1$ is in the set $I$ (this is the set where $N > \frac{1}{10,000} c_\diamond^{-1}$) (by virtue of Lemma 6.1). The second bullet condition in (8.13) implies that $s_1$ is in a set to be denoted by $\mathcal{J}_\mu$ which is defined as follows:

*A point* $r \in (0, c_\diamond^{-1}]$ *is in* $\mathcal{J}_\mu$ *when* $N(r) < \frac{1}{2}$ *and* $\underline{N}(\frac{1}{(1-\mu)^2} r) \geq 100 \underline{N}(r)$ .

(8.14)

Meanwhile, the third bullet condition in (8.13) implies that $s_1$ is in a set to be denoted by $\mathcal{K}$ which is defined by the rule

*A point* $r \in (0, c_\diamond^{-1}]$ *is in* $\mathcal{K}$ *when* $N(r) < \frac{1}{2}$ *and* $\int_{B_r} |F_A|^2 > 10^3 c_\diamond^4 r^2 r^2 K(r)^2 N(r)$ .

(8.15)

To be sure: The 4'th power of $c_\diamond$ appears in (8.15) whereas the 8'th power appears in the third bullet of (8.13); but keep in mind that $10^3 c_\diamond^4 \ll c_\diamond^8$ because $c_\diamond > 10^6$. The sets $\mathcal{J}_\mu$ and $\mathcal{K}$ are open, but not necessarily pairwise disjoint. Nor are they necessarily disjoint from $I$ or from $\mathcal{I}_\mu$. The rest of this subsection talks about $\mathcal{J}_\mu$; the next subsection talks about $\mathcal{K}$.

The salient property of $\mathcal{J}_\mu$ is that the function $\underline{N}$ is, to all intents and purposes, an increasing function on any of $\mathcal{J}_\mu$'s components (and, as explained momentarily, this must be true of $N$ also). To make a precise statement, let $(s_0, s_1)$ denote an interval that is entirely in $\mathcal{J}_\mu$. Introduce $m$ now to denote the largest integer such that $s_1 > \frac{1}{(1-\mu)^{2m}} s_0$. When $k \in \{0, 1, 2, \ldots, m+1\}$, let $t_k = \frac{1}{(1-\mu)^{2k}} s$ so that $t_0 = s_0$ and $t_m \in [(1-\mu)^2 s_1, s_1]$. The definition of $\mathcal{J}_\mu$ in (8.15) says that $\underline{N}(t_k) \geq 100 \underline{N}(t_{k-1})$ for all $k \in \{1, 2, \ldots, m\}$. This in turn implies that $\underline{N}(t_k) \geq 100^{k-1} \underline{N}(s_0)$. Then, since $N(t_k) \geq (1-\mu)^5 \underline{N}(t_k)$ by virtue of the top bullet of Lemma 6.3, it follows that $N(t_k) \geq (1-\mu)^5 100^k \underline{N}(s_0)$ for all $k \in \{1, 2, \ldots, m\}$. This bound with those in Lemma 6.1 lead to the following:

- $N(s) \geq (1-\mu)^{10} (\frac{s}{s_0})^{(\ln 10)/|\ln(1-\mu)|} \underline{N}(s_0)$ *if* $s \in [\frac{1}{(1-\mu)^2} s_0, \frac{1}{(1-\mu)^2} s_1]$ .
- $N(s) \geq (1-\mu)^5 \underline{N}(s_0)$ *if* $s \in [s_0, \frac{1}{(1-\mu)^2} s_0]$ .

(8.16)

The bounds in (8.16) are used to prove the following lemma:

**Lemma 8.3**: *There exists* $\kappa > 1$, *and given* $\mu \in (0, \frac{1}{100})$, *there exists* $\kappa_\mu > \kappa$; *and these numbers have the following significance: Fix* $r > 1$ *and suppose that* $(A, a)$ *is a solution to (1.5). Fix* $p \in X$ *to define the functions* $K$ *and* $N$, *and the set* $I$. *Then, fix* $\mu \in (0, \frac{1}{100})$ *to* *define the set* $\mathcal{J}_\mu$. *Suppose that* $s_1 \in (0, \kappa_\mu^{-1})$ *is a point in* $\mathcal{J}_\mu$, *that* $s_0 \in (0, s_1]$, *and that all of* $[s_0, s_1]$ *is in* $\mathcal{J}_\mu$. *Let* $\gamma$ *denote the minimum of* $90\underline{N}(s_0)$ *and* $c_\diamond^{-1}$. *Let* $s_2 = \frac{1}{(1-\mu)^4} s_1$. *Then*

- $N(s) \geq \gamma$ *for* $s \in [\frac{1}{(1-\mu)^2} s_0, s_2]$ .
- $K(s) \leq (\frac{s}{s_2})^{\gamma/2} K(s_2)$ *for* $s \in [s_0, s_2]$ .

***Proof of Lemma 8.3***: Let $t_1 = \frac{1}{(1-\mu)^2} s_0$. It follows from (8.16) and from the top bullet of Lemma 6.3 and from Lemma 6.1 that $N(s)$ is larger than the minimum of $(1-\mu)^{10} 100\underline{N}(s_0)$ and $c_\diamond^{-1}$ and if $s \in [t_1, s_2]$. Note in particular that this minimum is no greater than the minimum of $90\underline{N}(s_0)$ and $c_\diamond^{-1}$. Denote the latter minimum by $\gamma$. The $N \geq \gamma$ lower bound on the interval $[t_1, s_2]$ together with (3.9) implies that $K(s) \leq (\frac{s}{s_2})^\gamma K(s_2)$ for all $s \in [t_1, s_1]$. An appeal to Lemma 7.2 then gives the bound $K(s) \leq (\frac{s}{s_2})^{\gamma/2} K(s_2)$ for $s \in [s_0, t_1]$.

**d) The set** $\mathcal{K}$

Let $\kappa_\ddagger$ denote here the version of the number $\kappa$ that appears in Proposition 5.1. Define a function $M$ on $(0, \kappa_\ddagger^{-1})$ by the rule $s \to M(s) = 2\kappa_\ddagger r^2 s^2 K(s)^2 N(s)$. It is a consequence of (3.7) that this function can be used in Proposition 5.1. Since $c_\diamond$ is much greater than $\kappa_\ddagger$, it follows that if $r \in \mathcal{K}$, then

$$\int_{B_r} |F_A|^2 > M(r).$$

(8.17)

Now, if $s \geq r$ then the integral of $|F_A|^2$ over $B_s$ is also greater $10^3 c_\diamond^4 r^2 r^2 K(r)^2 N(r)$. Therefore, if $s \in (r, 4r)$, then the integral of $|F_A|^2$ over $B_s$ is greater $60 c_\diamond^4 r^2 s^2 K^2(s) N(r)$. As a consequence, the integral of $|F_A|^2$ over $B_s$ will be greater than $M(s)$ unless $N(s)$ is greater than $60 c_\diamond^3 N(r)$. In this event, it would follow from Lemma 6.1 that $N(4r)$ is greater than the minimum of $\frac{1}{1000}$ and $c_\diamond^3 N(r)$. By way of a summary, one (or both) of the following conditions holds if $r \in \mathcal{K}$:

- $\int_{B_s} |F_A|^2 > M(s)$ *for all* $s \in (r, 4r)$.
- $N(4r) \geq \min(\frac{1}{1000}, c_\diamond^3 N(r))$ .

(8.18)

More is said about these conditions momentarily.

Write $\mathcal{K}$ as union of two open sets to be denoted by $\mathcal{K}_1$ and $\mathcal{K}_2$. A point $r$ from $\mathcal{K}$ is in $\mathcal{K}_1$ when the top bullet in (8.18) is true; and it is in $\mathcal{K}_2$ when $N(4r)$ is greater than the minimum of $\frac{1}{2000}$ and $\frac{1}{2} c_\diamond^3 N(r)$. Both $\mathcal{K}_1$ and $\mathcal{K}_2$ are open sets (they can intersect). A discussion that is very much like that in Section 7c leads to the observation that N is, to all intents and purposes, an increasing function on any interval in the set $\mathcal{K}_2$. What is meant by this is made precise in the following lemma.

**Lemma 8.4**: *There exists $\kappa > 1$, and given $\mu \in (0, \frac{1}{100})$, there exists $\kappa_\mu > \kappa$; and these numbers have the following significance: Fix $r > 1$ and suppose that $(A, a)$ is a solution to (1.5). Fix $p \in X$ to define the functions $K$ and $N$ and the set $\mathcal{K}_2$. Suppose that $s_1$ is a point from $(0, \kappa^{-1})$ in $\mathcal{K}_2$, that $s_0 \in (0, s_1]$, and that all of $[s_0, s_1]$ is in $\mathcal{K}_2$. Let $\gamma$ denote the minimum of $c_\diamond^2 N(s_0)$ and $\frac{1}{20,000}$. Let $s_2 = 16 s_1$. Then*

- $N(s) \geq \gamma$ *for* $s \in [4 s_0, s_2]$.
- $K(s) \leq (\frac{s}{s_2})^{\gamma/2} K(s_2)$ *for* $s \in [s_0, s_2]$.

***Proof of Lemma 8.4***: If $s_1 > 4 s_0$ and if $s \in [4 s_0, 4 s_1]$, then $N(s)$ is greater than the minimum of $\frac{1}{2000}$ and $2^{-m} c_\diamond^{3m} N(s_0)$ with m being least positive integer such that $s \geq 4^m s_0$. This observation leads to the bound for $N(s)$ in the top bullet for $s \in [4 s_0, 4 s_1]$. The bound in the top bullet then follows also for $s \in [4 s_1, 16 s_1]$ via an appeal to Lemma 6.1. Meanwhile, the bound in the top bullet of the lemma implies (via (3.9)) the bound in the lower bullet for $s \in [4 s_0, 16 s_1]$ (which is $[4 s_0, s_2]$). Given this last bound, then the bound in the lower bullet of the lemma for $s \in [s_0, 4 s_0]$ follows from an appeal to Lemma 7.2.

Suppose now that $r \in \mathcal{K}$ and that the top bullet of (8.18) holds (which means that r is in $\mathcal{K}_1$). As a consequence of Proposition 5.1, there exists $r_1 \in [4r, c_\diamond^{-1}]$ such that

$$\int_{B_s} |F_A|^2 \geq M(s)$$

(8.19)

at each s in $[r, r_1]$ with this being an equality at $s = r_1$. Moreover,

$$\int_{B_s} |F_A|^2 \leq 64 (\tfrac{s}{r_1})^{2 + 1/\kappa_\ddagger} 2 \kappa_\ddagger r^2 r_1^2 K^2(r_1) N(r_1)$$

(8.20)

for each $s \in [r, \tfrac{1}{4} r_1]$ (and in particular $s = r$) unless $r_1 = c_\diamond^{-1}$ in which case

$$\int_{B_s} |F_A|^2 \leq \kappa_\ddagger^2 (\tfrac{s}{r_1})^{1/\kappa_\ddagger} r_1^2 \, r^2.$$

(8.21)

Supposing that (8.20) holds, then (8.20) and the top bullet in (8.18) imply the inequality

$$K(r)^2 N(r) \le c_0^{-3} \left(\tfrac{r}{r_1}\right)^{1/\kappa_\ddagger} K^2(r_1) N(r_1) \ .$$

(8.22)

On the other hand, if $r_1 = c_0^{-1}$ so that (8.21) holds, then (8.21) and (8.18) imply that

$$K^2(r) N(r) \le c_0^{-2} \left(\tfrac{r}{r_1}\right)^{1/\kappa_\ddagger} \ .$$

(8.23)

The appearance of $N(r_1)$ in (8.22) is awkward when $N(r_1)$ is big, and to be definite, when $N(r_1)$ is greater than 1. In this event, the factor of $N(r_1)$ can be removed from (8.22) at the expense of replacing $K(r_1)$ by $K(2r_1)$ and by replacing the factor $c_0^{-3}$ by $\tfrac{4}{\ln(2)} c_0^{-3}$ (which is less than $c_0^{-2}$) by invoking (5.24). Thus,

$$K(r)^2 N(r) \le c_0^{-2} \left(\tfrac{r}{2r_1}\right)^{1/\kappa_\ddagger} K(2r_1)^2 \ ,$$

(8.24)

This last inequality will be used in lieu of (8.22) when $N(r_1) > 1$.

The inequalities in (8.22) and/or (8.24) or (8.23) have immediate consequences that are summarized in the next lemma.

**Lemma 8.5**: *There exists $\kappa > 1$ with the following significance: Fix $r > 1$ and suppose that $(A, a)$ is a solution to (1.5). Fix $p \in X$ to define the functions $K$ and $N$ and the set $\mathcal{K}_1$. Suppose that $r \in (0, \kappa^{-1}]$ is a point from $\mathcal{K}_1$. Define $r_1 \in [4r, c_0^{-1}]$ as described above so that (8.19) and (8.20) hold unless $r_1 = c_0^{-1}$, in which case (8.21) holds. Let $r_\ddagger$ denote either $r_1$ or $2r_1$ depending on whether $N(r_1) < 1$ or not.*

- *Suppose that $r_1 < c_0^{-1}$.*
  a) $K(r)^2 \le \tfrac{1}{N(r) c_0^2} \left(\tfrac{r}{r_\ddagger}\right)^{1/\kappa_\ddagger} K(r_\ddagger)^2 \ .$
  b) *Fix $s \in [r, r_\ddagger]$ such that the interval $(r, s)$ is entirely in $\mathcal{I}$. If $N(r) \ge c_0^{-2}$ or if $\left(\tfrac{s}{r_\ddagger}\right)^{1/3\kappa_\ddagger} < N(r) c_0^2$, then $K(s)^2 \le \left(\tfrac{s}{r_\ddagger}\right)^{1/3\kappa_\ddagger} K(r_\ddagger)^2 \ .$*
- *Suppose that $r_1 = c_0^{-1}$.*
  a) $K(r)^2 \le \tfrac{1}{N(r) c_0^2} \left(\tfrac{r}{r_1}\right)^{1/\kappa_\ddagger} \ .$
  b) *Fix $s \in [r, r_\ddagger]$ such that the interval $(r, s)$ is entirely in $\mathcal{I}$. If $N(r) \ge c_0^{-2}$ or if $\left(\tfrac{s}{r_1}\right)^{1/3\kappa_\ddagger} < N(r) c_0^2$, then $K(s)^2 \le \left(\tfrac{s}{r_1}\right)^{1/3\kappa_\ddagger} \ .$*

*Proof of Lemma 8.5*: The Item a) of the top bullet follows directly from either (8.22) or (8.24) depending on whether $r_\ddagger = r_1$ or $r_\ddagger = 2r_1$ (thus, whether $N(r_1) < 1$ or not). Meanwhile, Item a) of the second bullet follows from (8.22). To prove Item b) of the first bullet, invoke (8.4) with s in (8.23) taken to be r and with s´ in (8.23) taken to be the value of s in Item b) in the top bullet of the lemma. This (8.23) bound with Item a) of the top bullet of the lemma to see that

$$K(s)^2 \leq \frac{1}{N(r) c_\diamond^2} \left(\frac{s}{r}\right)^{1/c_\diamond} \left(\frac{r}{r_\ddagger}\right)^{1/\kappa_\ddagger} K(r_\ddagger)^2$$

(8.25)

when the interval (r, s) is in I. The bound in Item b) of the top bullet follows from this because $c_\diamond > 3\kappa_\ddagger$. Essentially the same argument proves Item b) of the second bullet of the lemma from (8.4) and Item a) of the second bullet of the lemma.

The following lemma talks about the size of the number $N(r_1)$ from (8.22). (It is important in subsequent applications that $N(r_1)$ not be too small.)

**Lemma 8.6**: *There exists $\kappa > 1$ with the following significance: Fix $r > 1$ and suppose that $(A, a)$ is a solution to (1.5). Fix $p \in X$ to define the functions K and N and the set $\mathcal{K}_1$. Suppose that $r \in (0, \kappa^{-1}]$ is a point from $\mathcal{K}_1$. Define $r_1 \in [4r, c_\diamond^{-1}]$ as described above; and assume that (8.19) and (8.20) hold.*
- *If $r_1$ is in r's component of $\mathcal{I}$, then $N(r_1) \geq c_\diamond^3 \left(\frac{r_0}{r}\right)^{1/2\kappa_\ddagger} N(r)$.*
- *If $r_1$ is not in r's component of $\mathcal{I}$, then $N(r_1) \geq \frac{1}{10,000} c_\diamond^{-1}$; as a consequence, $r_1$ is in I.*

*Proof of Lemma 8.6*: Suppose first that $r_1$ is in r's component of $\mathcal{I}$. Granted this, then (8.4) holds with $s = r$ and with $s´ = r_1$. This version of (8.4) and (8.22) are not compatible unless the assertion in the lemma top bullet holds.

Now suppose that $r_1$ is <u>not</u> in r's component of $\mathcal{I}$. If $r_1$ is not in I, then it is in the set $\mathcal{I}$. If $r_1$ is in $\mathcal{I}$ but not in r's component of $\mathcal{I}$, then there is a least $s \in (r, r_1)$ such that the whole of the interval $(s, r_1)$ is in $\mathcal{I}$. Denote this point by $s_1$. This $s_1$ has $N(s_1) = c_\diamond^{-1}$ because $s_1$ is in the closure of $\mathcal{I}$ but not in $\mathcal{I}$. If it turns out that $N(r_1)$ is greater than $\frac{1}{64} N(s_1)$, then $N(r_1) > \frac{1}{10,000} c_\diamond^{-1}$ so $r_1$ is in I. For example, this occurs if $s_1 \geq \frac{1}{4} r_1$ by virtue of Lemma 6.1. This bound $N(r_1) > \frac{1}{64} N(s_1)$ must also hold if $s_1$ is less than $\frac{1}{4} r_1$ and here is why: If $s_1 < \frac{1}{4} r_1$, then (8.19) and (8.20) lead to the bound

$$K^2(s_1) N(s_1) \leq 64 \left(\frac{s_1}{r_1}\right)^{1/\kappa_\ddagger} K^2(r_1) N(r_1) .$$

(8.26)

If $N(r_1) < \frac{1}{64} N(s_1)$, then (8.26) implies that $K(s_1)^2 \leq \left(\frac{s_1}{r_1}\right)^{1/\kappa_\ddagger} K(r_1)^2$. This last bound is nonsense because it runs afoul of the $s = s_1$ and $s´ = r_0$ version of (8.4).

**e) Proof of Proposition 7.3**

This part of the subsection proves the assertion of Proposition 7.3. The proof given below is long and intricate because different arguments are needed for various cases that depend on particulars of the behavior of N. For this reason, the argument is presented in eight parts.

Here is a point to keep in mind with regards to the proof of Proposition 7.3: The second bullet of Proposition 7.3 holds if the points $s_0$ and $s_1$ obey

*Let $s_*$ denote $s_0$ or $s_1$. If $s_* \leq \rho_{\mu,\varepsilon}$, then either $N(s_*) \geq 4N(r)$ or $N(s_*) > \frac{1}{10,000} c_\diamond^{-1}$.*

(8.27)

To be sure, this differs from the second bullet because of the appearance here of $4N(r)$ instead of $2N(r)$. The verification of (8.27) for all but the final case discussed below facilitates the verification of the second bullet of Proposition 7.3 for the last case.

Some other points to keep in mind: The threshold for membership in the closure of the set I is $N(\cdot) \geq \frac{1}{10,000} c_\diamond^{-1}$. In particular, if points $s_0$ and $s_1$ satisfying the first bullet of Proposition 7.3 are in the closure of I, then they also satisfy the conditions in (8.27) and thus in the second bullet of Proposition 7.3. Also keep in mind that membership in $\mathcal{I}$ is the condition $N(\cdot) < c_\diamond^{-1}$, and thus any point on the boundary of $\mathcal{I}$ has $N(\cdot) = c_\diamond^{-1}$. The latter fact is repeatedly used.

*Part 1*: If the whole of the interval $(r, \frac{2}{(1-\mu)^2} r)$ is in I, then $s_0$ can be set equal to $s_1$ with both equal to $\frac{2}{(1-\mu)^2} r$. The first and second bullets of Proposition 8.1 (and (8.6)) follow directly because $N(s_0)$ and $N(s_1)$ are greater than or equal to $\frac{1}{10,000} c_\diamond^{-1}$. Meanwhile the bound that is asserted by the third bullet of Proposition 7.3 holds with any $\kappa_\mu \geq 2$ because the bound in this case is an instance of (8.2).

*Part 2*: Now suppose that some point from the interval $[r, 5r]$ is <u>not</u> in the set $\mathcal{I}$. This point (it will be called $\rho$) must be such that $N(\rho) \geq c_\diamond^{-1}$. Lemma 7.2's second bullet implies the assertions of Proposition 7.3 (and (8.27)) with $s_0 = s_1 = 3\rho$ and any $\kappa_\mu \geq 2$.

With the preceding understood, the unwritten assumption henceforth is that the whole interval $[r, 5r]$ is in $\mathcal{I}$.

*Part 3*: The subsequent parts of the proof consider the point

$$\hat{r} = \frac{1}{(1-\mu)^2} r.$$

(8.28)

This point is a member of $\mathcal{I}$ so it is in one (or more) of the subsets $\mathcal{I}_\mu$ (from Section 7b) or $\mathcal{J}_\mu$ or $\mathcal{K}$ (both from Section 7c); and if it is from $\mathcal{K}$, then it is in one its subsets $\mathcal{K}_1$ or $\mathcal{K}_2$ (from Section 7d). What was just said is also true for r, but, as it turns out, $\hat{r}$ has an advantage over r in the upcoming Part 7. The various cases when $\hat{r}$ is in $\mathcal{I}_\mu, \mathcal{J}_\mu, \mathcal{K}_1$ or $\mathcal{K}_2$, are considered in reverse order starting in Part 4 with the case of $\mathcal{K}_2$.

In all of these case, there are three facts to keep in mind with regards to $\hat{r}$ and r:

- $N(\hat{r}) \geq (1-\mu)^{10} N(r)$.
- $\underline{N}(\hat{r}) \geq (1-\mu)^{10} \underline{N}(r)$.
- *If* $s_1 > 2\hat{r}$ *(which is* $\frac{2}{(1-\mu)^2}$ *r) and* $\sigma > 1$ *and* $K(s) \leq (\frac{s}{s_1})^{1/\sigma} K(s_1)$ *for all* $s \in [\hat{r}, s_1]$, *then* $K(s) \leq (\frac{s}{s_1})^{1/(2\sigma)} K(s_1)$ *for all* $s \in [r, s_1]$.

(8.29)

The assertions of the first two bullets are direct consequences of Lemmas 6.1 and 6.2; and the assertion in the third bullet follows from Lemma 7.2.

*Part 4*: Consider the event that $\hat{r}$ is in the set $\mathcal{K}_2$. Lemma 8.4 can be invoked in this case with its versions of $s_0$ and $s_1$ being equal to $\hat{r}$. The point $s_2$ from Lemma 8.4 can then be used for the Proposition 7.3 version of $s_0$ and the Proposition 7.3 version of $s_1$. It follows from the first bullet of Lemma 8.4 that the number $s_2$ from Lemma 8.4 satisfies the requirements of the second bullet of Proposition 7.3 and so the first two bullets of Proposition 7.3 are satisfied with $s_0$ and $s_1$ both being Lemma 8.4's number $s_2$. Note in this regard that Lemma 8.4's version of $s_2$ is such that $N(s_2)$ in this case is greater than the minimum of $c_\diamond^2 N(\hat{r})$ and $\frac{1}{20,000}$. This guarantees that (8.27) holds (and therefore Proposition 7.3's second bullet) when Proposition 7.3's version of $s_0$ and $s_1$ are Lemma 8.4's version of $s_2$.

The bound in the third bullet of Proposition 7.3 follows if Proposition 7.3's version of $\kappa_\mu$ is any number greater than 4. Indeed, the inequality in the third bullet of Proposition 7.3 holds for $s \in [\hat{r}, s_0]$ with $\kappa_\mu$ greater than 4 because of the second bullet of Lemma 8.4; and this fact with the first and third bullets of (8.7) implies the inequality in Proposition 7.3's third bullet when $\kappa_\mu \geq 4$.

*Part 5*: This part of the proof discusses the case when $\hat{r}$ is in the set $\mathcal{K}_1$. There are four steps to the proof in this case.

<u>Step 1</u>: Use the point $\hat{r}$ for the Lemma 8.5 version of the number r. Let $r_1$ denote the output point that is supplied by Lemma 8.5. Assume here in this step that this number $r_1$ obeys $r_1 = c_\diamond^{-1}$. In this event, take Proposition 7.3's version of $s_1$ to be $c_\diamond^{-1}$. Item b) of the second bullet of Lemma 8.5 with r replaced by $\hat{r}$ holds.

If $N(r) \geq (1-\mu)^{-10} c_\diamond^{-2}$, set $s_0'$ to be the largest s in the interval $(\hat{r}, c_\diamond^{-1}]$ such that the interval (r, s) is entirely in $\mathcal{I}$. (Keep in mind that this largest s is greater than 5r because r is in $\mathcal{I}$.) Item b) of the second bullet of Lemma 8.5 (with r replaced by $\hat{r}$) and the first and third bullets of (8.29) together imply that

$$\kappa(s) \le (\tfrac{s}{s_1})^{1/6\kappa_{\ddagger}} \quad for\ s \in [r, s_0''].$$
(8.30)

If $N(r) \ge \varepsilon$ with $\varepsilon$ less than $(1-\mu)^{-10} c_\diamond^{-2}$, then this same part of Lemma 8.5 and the same bullets of (8.29) also lead to (8.30) provided that the number r is less than the minimum of $\tfrac{1}{100} c_\diamond^{-1}$ and $\tfrac{1}{100} (\varepsilon c_\diamond^2)^{3\kappa_{\ddagger}} c_\diamond^{-1}$; and provided that the number $s_0'$ that appears in (8.30) is redefined so as to be the largest s greater than r but not greater than the minimum of $c_\diamond^{-1}$ and $(\varepsilon c_\diamond^2)^{3\kappa_{\ddagger}} c_\diamond^{-1}$ such that the interval (r, s) entirely in $\mathcal{I}$.

Since $\kappa(s_1) \ge c_\diamond^{-1}$, the inquality in (8.30) implies the inequality in the third bullet of Proposition 7.3. with $\kappa_\mu$ any number greater than 6 if $s_0$ is chosen to be the smaller of $s_0'$ and $c_\mu^{-1} c_\diamond^{-1}$ (which is $c_\mu^{-1} s_1$); and if r is assumed to be less than $\tfrac{1}{100} s_0$.

Step 2: This step and the subsequent steps assume that the number $r_1$ from Lemma 8.5 (with input $\hat{r}$) is strictly less than $c_\diamond^{-1}$. The point $r_1$ is either in $\hat{r}$'s component of $\mathcal{I}$ or not. This step and Step 3 treat that case when $r_1$ is not in $\hat{r}$'s component of $\mathcal{I}$.

By way of a parenthetical remark if $\hat{r}$ is in I, then $r_1$ can not be in $\hat{r}$'s component of $\mathcal{I}$; The reason is this: If $r_1$ and $\hat{r}$ were in the same component of $\mathcal{I}$, then $r_{\ddagger}$ would be $r_1$; and then Item a) of the top bullet of Lemma 8.5 implies that $\kappa(\hat{r}) \le (\tfrac{\hat{r}}{r_1})^{1/\kappa_{\ddagger}} \kappa(r_1)^2$ which runs afoul of (8.4). Note in this regard that $\hat{r}$ is in I if $N(r) \ge (1-\mu)^{-10} \tfrac{1}{10,000} c_\diamond^{-1}$ by virtue of the top bullet in (8.30).

If $r_1$ is not in $\hat{r}$'s component of $\mathcal{I}$, then $r_1$ must be in I according to Lemma 8.6. Moreover, if by chance $N(r_1) \ge 1$, then $N(2r_1) \ge c_\diamond^{-1}$ (this is implied by Lemma 6.1). Thus, the number $r_{\ddagger}$ from Lemma 8.5 is in I whether $r_{\ddagger}$ is $r_1$ or $2r_1$. With the preceding understood, take Proposition 7.3's version of $s_1$ to be $r_{\ddagger}$. Let $\rho$ denote the largest number $s \in (\hat{r}, s_1)$ such that $(\hat{r}, s)$ is in $\mathcal{I}$ and take Proposition 7.3's version of $s_0$ to equal $\rho$. This $s_0$ has $N(s_0) = c_\diamond^{-1}$ (because $\rho$ is on the boundary of $\mathcal{I}$) and so $s_0$ and $s_1$ obey the requirements of the first two bullets of Proposition 7.3 and (8.27).

If $N(r) \ge (1-\mu)^{-10} c_\diamond^{-2}$, then the inequality of the third bullet in Proposition 7.3 with any $\kappa_\mu \ge 6$ is a consequence of the first and third bullets of (8.27) and Item b) of the first bullet of Lemma 8.5 (with r replaced therein by $\hat{r}$). If $N(r) \le (1-\mu)^{-10} c_\diamond^{-2}$ with r obeying

$$r \le ((1-\mu)^{10} N(r) c_\diamond^2)^{3\kappa_{\ddagger}} r_{\ddagger} \ ,$$
(8.31)

then the inequality in Proposition 7.3's third bullet still holds (with $s_1$ and $s_0$ and $\kappa_\mu$ as just described). The proposition's third bullet in this subcase is again a direct consequence of the first and third bullet of (8.27) and Item b) of the first bullet of Lemma 8.5.

Step 3: This step assumes that $N(r) \leq (1-\mu)^{-10} c_\diamond^{-2}$ and that the inequality in (8.31) is reversed. In this case, define $\rho$ as in the previous step (it is the largest s in $(r, r_\ddagger)$ so that $(r, s)$ is entirely in $\mathcal{I}$), and take $s_0$ and $s_1$ in Proposition 7.3 to be the number $4\rho$. Since $N(\rho) = c_\diamond^{-1}$ (because $\rho$ is on the boundary of $\mathcal{I}$), the whole of $[\rho, 4\rho]$ is in $I$ (appeal to Lemma 6.1). This last observation has two consequences: First, $N(4\rho) > \frac{1}{10,000} c_\diamond^{-1}$ which implies that if both $s_0$ and $s_1$ equal $4\rho$, then they obey the requirements of the first two bullets in Proposition 7.3 and also (8.27). It also implies that

$$K(s) \leq (\tfrac{s}{s_1})^{1/(10,000 c_\diamond)} K(s_1) \quad \textit{for all } s \in [\tfrac{1}{4} s_1, s_1]$$

(8.32)

which is an instance of (8.2). With (8.32) in hand, now invoke Lemma 7.2 to see that

$$K(s) \leq (\tfrac{s}{s_1})^{\gamma_\ddagger} K(s_1) \quad \textit{for all } s \in [r, s_1]$$

(8.33)

with $\gamma_\ddagger = \frac{1}{10,000} c_\diamond^{-1} \ln(4)/\ln(4\tfrac{\rho}{r})$. The third bullet of Proposition 7.3 follows from (8.33) because the reverse inequality in (8.31) implies that $\gamma_\ddagger \geq c_\mu^{-1} N(r)$. Indeed, this is so for the following reason: The reverse inequality in (8.31) makes $r \geq ((1-\mu)^{10} N(r) c_\diamond^2)^{3\kappa_\ddagger} \rho$ because $r_\ddagger > \rho$. This implies that $\ln(\tfrac{\rho}{r}) \leq 3\kappa_\ddagger(|\ln N(r)| + 1)$ and thus, that $1/\ln(4\tfrac{\rho}{r})$ is greater than $\frac{1}{100} \kappa_\ddagger^{-1} N(r)$.

Step 4: This step assumes that $r_1 < c_\diamond^{-1}$ and that $r_1$ is in $\hat{r}$'s component of $\mathcal{I}$. Since $r_1$ is in $\mathcal{I}$, the number $N(r_1)$ is less than 1 (it is less than $c_\diamond^{-1}$) which implies that $r_\ddagger$ from Lemma 8.5 is equal to $r_1$. This is the first observation.

The second observation is that (8.22) and (8.4) are compatible only if

$$N(r_1) \geq c_\diamond^3 (\tfrac{r_1}{r})^{1/3\kappa_\ddagger} N(\hat{r}) ;$$

(8.34)

and this implies, in turn, that

$$N(r_1) \geq (1-\mu)^{10} c_\diamond^3 (\tfrac{r_1}{r})^{1/3\kappa_\ddagger} N(r)$$

(8.35)

since $N(\hat{r}) \geq (1-\mu)^{-10} N(r)$ (look at (8.29)). Moreover, since $N(r_1) \leq c_\diamond^{-1}$ (because $r_1$ is in $\mathcal{I}$), the inequality in (8.35) implies that

$$r \geq (\tfrac{1}{2} c_\diamond^{-4} N(r))^{3\kappa_\ddagger} r_1$$

(8.36)

so r can not be arbitrarilly small relative to $r_1$. This fact is going to be used momentarily.

Meanwhile, it is a consequence of (8.35) and Lemma 6.1 that $N(s)$ for $s \in [r_1, 4r_1]$ is greater than the minimum $\tfrac{1}{32} c_\diamond^3 (\tfrac{r_1}{r})^{1/3\kappa_\ddagger} N(r)$ and $c_\diamond^{-1}$. Therefore, the conditions of the first and second bullets of Proposition 7.3 and (8.27) are met if $s_0$ and $s_1$ are both $4r_1$.

To see about the third bullet of Proposition 7.3, introduce by way of notation $\gamma$ to denote the minimum $\tfrac{1}{32} c_\diamond^3 (\tfrac{r_1}{r})^{1/3\kappa_\ddagger}$ and $c_\diamond^{-1} N(r)^{-1}$. This notation is used to write the bound from the preceding paragraph for $N$ on $[\tfrac{1}{4} s_1, s_1]$ as $N(s) \geq \gamma N(r)$. Granted this lower bound for $N$, then (3.9) implies that

$$K(s) \leq (\tfrac{s}{s_1})^{\gamma N(r)} K(s_1)$$

(8.37)

for all $s \in [\tfrac{1}{4} s_1, s_1]$. It then follows from Lemma 7.2 that

$$K(s) \leq (\tfrac{s}{s_1})^{\gamma' N(r)} K(s_1)$$

(8.38)

for all $s \in [r, s_1]$ with $\gamma'$ defined by the rule

$$\gamma' = \gamma \ln(4)/\ln(4\tfrac{r_1}{r}).$$

(8.39)

Thus, the third bullet of Proposition 7.3 is obeyed because, as explained next, $\gamma' \geq c_\mu^{-1}$. Here is why this bound holds: If $\gamma$ is $\tfrac{1}{32} c_\diamond^3 (\tfrac{r_1}{r})^{1/3\kappa_\ddagger}$, then this divided by the number $\ln(\tfrac{r_1}{r})$ is greater than $c_\diamond^{-1}$. On the other hand, if $\gamma = c_\diamond^{-1} N(r)^{-1}$, then $\gamma/\ln(\tfrac{r_1}{r})$ is no smaller than $c_\mu^{-1} (\tfrac{1}{N(r)|\ln N(r)|})$ because of (8.36), and this is no smaller than $c_\mu^{-1}$.

*Part 6*: This step considers the case where $\hat{r}$ from (8.26) is in the set $\mathcal{J}_\mu$. The plan for what follows is to first invoke Lemma 8.3 with its versions of $s_0$ and $s_1$ both equal to $4\hat{r}$. It is a consequence of the second bullet of (8.29) and the top bullet of Lemma 8.3 that the value of $N$ at the point $\tfrac{1}{(1-\mu)^4} \hat{r}$ is no smaller than the minimum of $80 N(r)$ and $c_\diamond^{-1}$. It follows from this and Lemma 6.1 that the values of $N$ on the interval $[\tfrac{1}{(1-\mu)^4} \hat{r}, \tfrac{4}{(1-\mu)^4} \hat{r}]$ are greater than the minimum of $4 N(r)$ and $\tfrac{1}{20} c_\diamond^{-1}$. Therefore, the first and second bullets of Proposition 7.3 and (8.27) hold with the Proposition 7.3 versions of the numbers $s_0$ and $s_1$ being $\tfrac{4}{(1-\mu)^4} \hat{r}$. Meanwhile, this lower bound for $N$ and (3.9) lead to the bound

$$K(s) \le (\tfrac{s}{s_1})^{\gamma N(r)} K(s_1) \quad for \quad s \in [\tfrac{1}{4} s_1, s_1]$$

(8.40)

with $\gamma$ here denoting the minimum of 4 and $\tfrac{1}{20} c_\diamond^{-1} N(r)^{-1}$. The inequality in (8.40) and the third bullet of (8.29) lead directly to the inequality in the third bullet of Proposition 7.3 with any choice of $\kappa_\mu$ greater than 100.

*Part 7*: The last case to resolve (which is done here and in Part 8) is the case where $\hat{r}$ is not in $\mathcal{J}_\mu$ nor in $\mathcal{K}$ (which is the union of $\mathcal{K}_1$ and $\mathcal{K}_2$). In this event, $\hat{r}$ is in the set $\mathcal{I}_\mu$. Since $\mathcal{I}_\mu$ is an open set, there is a largest $s \in (\hat{r}, c_\diamond^{-1}]$ such that the whole of the interval $(\hat{r}, s)$ is in $\mathcal{I}_\mu$. Denote this number by $\rho$. The conclusions of Lemma 8.2 hold with the Lemma 8.2 version of $s_0$ being $\hat{r}$ and with the Lemma 8.2 version of $s_1$ being $\rho$. The first bullet of Lemma 8.2 in this instance implies that N obeys

$$N(s) \ge (1-\mu)^{20} \tfrac{1}{\eta(r)} N(r) - z_\mu s^2 \quad for\ all \quad s \in [r, \rho]$$

(8.41)

This is so because $\underline{N}(\hat{r})$ is greater than $(1-\mu)^{20} N(r)$ because of the second bullet in (8.29).

By way of a parenthetical remark, it is at this point in the proof of Proposition 7.3 where $\hat{r}$ has the advantage over r. If r were used here, then (8.41) could only be guaranteed if N(r) was replaced by $\tfrac{1}{10,000} N(r)$ which is a factor of $\tfrac{1}{10,000}$ smaller. (This extra factor of $\tfrac{1}{10,000}$ comes from the second bullet of (8.6).) It is important in the final steps that there is no such extra factor because then N(s) for $s \in [r, \rho]$ is very close (when $\eta(r) \sim 1$) to N(r); and this property of N(s) is needed to find $s_0$ and $s_1$ that obey the second bullet of Proposition 7.3. This is also where the factor of 4 in (8.27) is needed. Looking ahead, the key point is that 4 times the right hand side of (8.41) is greater than 2N(r) when $\eta(r)$ is close to 1 and s is small.

The second bullet of Lemma 8.2 and (8.29) imply that if $s_1 \in (r, \rho]$ is given, then

$$K(s) \le e^{z_\mu (s_1-s)/c_\diamond} (\tfrac{s}{s_1})^{\delta(r) N(r)} K(s_1) \quad for\ all\ s \in [r, s_1].$$

(8.42)

with $\delta(r)$ denoting the number $(1-\mu)^{20} \tfrac{1}{\eta(r)}$. It is important to note that this number $\delta(r)$ can not be arbitrarilly small when $rrK(r) \ge 1$. This is because $\eta(r) \le e^{zx}$ when $rrK(r) > 1$ (look at the definition of $\eta$ in (8.12)). In particular, $\delta(r)$ is no less than

$$\delta_\mu = (1-\mu)^{20} e^{-zx}$$

(8.43)

when the condition $rrK(r) > 1$ is enforced.

The rest of this step completes the proof of Proposition 7.3 when $\rho \geq c_\diamond^{-2}$. In this case, the number $s_1$ for Proposition 7.3 can be taken to be $c_\diamond^{-2}$. Then, (8.42) implies that

$$K(s) \leq \left(\frac{s}{s_1}\right)^{\frac{1}{2}\delta_\mu N(r)} K(s_1) \quad \text{for all } s \in [r, s_0]$$

(8.44)

if $s_0$ obeys the condition

$$e^{z_\mu s_1/c_\diamond} \left(\frac{s_0}{s_1}\right)^{\frac{1}{2}\delta_\mu \varepsilon} \leq 1.$$

(8.45)

In particular, (8.45) is obeyed for any choice of $s_0$ less than $e^{-2z_\mu/(c_\diamond^3 \varepsilon \delta_\mu)} c_\diamond^{-2}$.

Granted the preceding, assume henceforth that $\rho < c_\diamond^{-2}$.

*Part 8*: The $s_1 = \rho$ version of (8.13) must hold because $\rho$ is the largest $s \in (\hat{r}, c_\diamond^{-1}]$ such that the whole of the interval $(\hat{r}, s)$ is in $\mathcal{I}_\mu$. This is to say that $\rho$ is described by one or more of the three bullets in (8.13). In each case, one of Parts 2-6 can be invoked with $\rho$ used in lieu of r to see that the conclusions of Proposition 7.3 and also (8.27) hold for the point $\rho$. Indeed, Part 2 can be invoked with $\rho$ replacing r if the first bullet of (8.13) describes $\rho$; and Part 6 can be invoked with $\rho$ used instead of r if the second bullet of (8.13) describes $\rho$, and Parts 4-5 can be invoked if the third bullet in (8.13) describes $\rho$.

With the preceding understood, let $s_0$ and $s_1$ denote the points supplied by Proposition 7.3 and (8.27) with its input being $\rho$ instead of r. These obey $\hat{s}_0 \geq 2r$ and $s_1 \geq s_0$ because $\rho > r$. Therefore, they obey the first bullet of the r version of Proposition 7.3. To see about the second bullet of the r version of Proposition 7.3, note that $s_0$ and $s_1$ obey (8.27) with $\rho$ used in lieu of r. What with (8.31), the latter version of (8.27) implies the following: If $s_*$ is $s_0$ or $s_1$, then

$$N(s_*) \geq 4((1-\mu)^{20} \frac{1}{\eta(r)} N(r) - z_\mu \rho^2) \quad \text{if} \quad N(s_*) \leq \frac{1}{10{,}000} c_\diamond^{-1} \text{ and } s_* \leq \hat{\rho}_{\mu,\varepsilon}$$

(8.46)

with $\hat{\rho}_{\mu,\varepsilon}$ being a number that depends only on $\mu$ and $\varepsilon$. Now, the key observation is that $\eta(r) \geq \frac{999}{1000}$ if $rrK(r) \geq c_\mu$; and also $z_\mu \rho^2 < \frac{1}{10{,}000} \varepsilon$ if $\rho < c_\mu^{-1}\sqrt{\varepsilon}$. Granted these conditions, then (8.46) implies that $N(s_*) \geq 2N(r)$ if $N(s_*) \leq \frac{1}{10{,}000} c_\diamond^{-1}$ and if $s_0$ is less than the minimum of the numbers $\hat{\rho}_{\mu,\varepsilon}$ and $c_\mu^{-1}\sqrt{\varepsilon}$. And, this implies what asserted by the second bullet of Proposition 7.3 for a suitable choice of $\rho_{\mu,\varepsilon}$.

The final task is to verify the third bullet of Proposition 7.3. This follows directly from the two observations. First, the $s_1 = \rho$ version of (8.42) and (8.43) imply that

$$K(s) \leq e^{z_\mu (\rho - s)/c_\diamond} \left(\tfrac{s}{\rho}\right)^{\delta_\mu N(r)} K(\rho) \quad \textit{for all } s \in [r, \rho]$$

(8.47)

when $rrK(r) \geq 1$. And, Proposition 7.3's third bullet applied to $\rho$ and (8.41) imply that

$$K(s) \leq \left(\tfrac{s}{s_1}\right)^{L(r)/c_\mu} \quad \textit{for all } s \in [\rho, s_0] \, .$$

(8.48)

These imply imply the third bullet of Proposition 7.3 with $\rho_{\mu,\varepsilon}$ being a number that depends only on $\mu$ and $\varepsilon$ as required.

## 9. Proof of Proposition 7.4

The proof of Proposition 7.4 is given momentarily in Section 9b. The intervening sections establish some preliminary results that are needed in Section 9b. To set the notation, and supposing that $\mu \in (0, \tfrac{1}{100})$ has been chosen, let $\kappa_{\diamond\mu}$ again denote the value of $\kappa_\mu$ from Proposition 7.3. Then, let $r_\diamond$ denote the value of r where $rrK(r) = \kappa_{\diamond\mu}$. (This number is well defined because the derivative of the function $r \to rrK(r)$ is positive where r is positive.)

### a) The set $\mathcal{L}$

The sets I and $\mathcal{I}$ were defined in Section 8, the former being the subset of $(0, c_\diamond^{-1}]$ where N is greater than $\tfrac{1}{10,000} c_\diamond^{-1}$, and the latter being the subset of this same interval where N is less than $c_\diamond^{-1}$. A subset $\mathcal{L}$ of $\mathcal{I}$ is defined by the rule whereby $r \in \mathcal{L}$ when

$$\int_{B_r} |F_A|^2 < c_\diamond^8 r^2 r^2 K(r)^2 N(r)$$

(9.1)

Note that the set $\mathcal{L}$ is an open subset of $\mathcal{I}$. The following lemma describes the functions K and N on $\mathcal{L}$.

**Lemma 9.1**: *There exists $\kappa > 1$; and given $\mu \in (0, \tfrac{1}{100}]$, there exists $\kappa_\mu > \kappa$; and these numbers have the following significance: Fix $r > 1$ and a pair $(A, a)$ solving (1.5). Fix $p \in X$ and use $(A, a)$ and p to define K, N, and the set $\mathcal{L}$. Fix a point $r \in (0, r_\diamond]$ from $\mathcal{L}$ such that $r < \kappa_\mu^{-1} N(r) r_\diamond$. Then, let $\rho$ denote the largest $s \in [r, \kappa_\mu^{-1} N(r) r_\diamond]$ such that the whole interval $[r, \rho)$ is in $\mathcal{L}$ (and thus in $\mathcal{I}$). If $s \in [r, \rho]$, then*
- $N(s) \geq \tfrac{7}{8} N(r)$ .
- $K(s) \leq \left(\tfrac{s}{\rho}\right)^{\tfrac{7}{8} N(r)} K(\rho)$ .

***Proof of Lemma 9.1***: As explained momentarily, if $s \in (0, (1-\mu)r_\diamond)$ and $\frac{1}{1-\mu} s \in \mathcal{L}$, then the derivative of the function $N$ at $s$ obeys

$$\frac{d}{ds} N \geq -c_\mu \frac{1}{r_\diamond} .$$
(9.2)

As a consequence of this, if $s \in (r, (1-\mu)r_\diamond)$ and if all of the interval $[r, \frac{1}{1-\mu} s)$ is in $\mathcal{L}$, then

$$N(s) \geq N(r) - c_\mu \frac{s}{r_\diamond} .$$
(9.3)

Therefore, if $r < c_\mu^{-1} N(r) r_\diamond$ and if $\rho$ is greater than $r$ but less than $c_\mu^{-1} N(r) r_\diamond$, and if the whole of the interval $[r, \rho)$ is in $\mathcal{L}$, then

$$N(s) \geq \tfrac{15}{16} N(r) \quad \textit{if } s \in [r, (1-\mu)\rho].$$
(9.4)

This bound with Lemma 6.1 imply the assertion of the top bullet of the lemma. The second bullet follows from the top bullet using the identity in (3.9).

The proof of (9.2) has three steps.

Step 1: Return to the equation in Proposition 6.2 for the derivative of $N$. Let $s$ denote a point from $[r, (1-\mu)r_\diamond]$. Since $\frac{1}{1-\mu} s$ is in $\mathcal{I}$, and since $|E_A + B_A|$ is no greater than $c_0 r^2 |a|^2$, it follows from Proposition 3.1 and (3.9) that $|E_A + B_A|$ on $B_s$ is bounded by $c_\mu r^2 K^2(\frac{1}{1-\mu} s)$. This, in turn, is bounded by $c_\mu r^2 K(s)^2$ for the following reason: The assumption that $\frac{1}{1-\mu} s$ is in $\mathcal{L}$ implies that $\frac{1}{1-\mu} s$ is in $\mathcal{I}$ and so $N(\frac{1}{1-\mu} s) < c_\diamond^{-1}$. Then, by virtue of Lemma 6.1, the function $N$ on $[s, \frac{1}{1-\mu} s]$ is at most $100 c_\diamond^{-1}$; and with this understood, then (3.9) bounds $K(\frac{1}{1-\mu} s)$ by $c_0 K(s)$. Anyway, because of this $r^2 K(s)^2$ bound for $|E_A + B_A|$, the equation in Proposition 6.2 requires that

$$\frac{d}{ds} N \geq -c_\mu \frac{1}{s^2} \int_{\partial B_s} |F_A| - c_\mu s$$
(9.5)

at any point $s \in [r, (1-\mu)r_\diamond]$ with $\frac{1}{1-\mu} s$ being in $\mathcal{L}$.

Step 2: This step and Step 3 derive a useful bound for the integral of $|F_A|$ over $\partial B_s$ from (9.5). To begin the derivation, reintroduce the $\mu = \tfrac{1}{4}$ version of $\chi_\mu$, thus $\chi_{1/4}$. Use this function to write $|F_A|$ on $\partial B_s$ as $(1-\chi_{1/4})|F_A|$ where it appears in (9.5); and then use the fundamental theorem of calculus (in reverse) to obtain the bound

$$\int_{\partial B_s} |F_A| \leq c_0 \left( \int_{B_s} |\nabla F_A| + \frac{1}{s} \int_{B_s} |F_A| \right) .$$

(9.6)

The latter bound (with the Cauchy-Schwarz inequality) leads to the bound

$$\int_{\partial B_s} |F_A| \leq c_0 s^2 \left( \int_{B_s} |\nabla F_A|^2 \right)^{1/2} + c_0 s \left( \int_{B_s} |F_A|^2 \right)^{1/2} .$$

(9.7)

To procede from here, invoke Proposition 5.4 for the integral of $|\nabla F_A|^2$ to obtain a bound for the integral of $|F_A|$ on $\partial B_s$ that reads

$$\int_{\partial B_s} |F_A| \leq c_\mu s \left( \int_{B_{s/(1-\mu)}} |F_A|^2 + r^4 s^4 K(\tfrac{1}{1-\mu} s)^4 (s^2 + N(\tfrac{1}{1-\mu} s)) \right)^{1/2} .$$

(9.8)

Keep in mind now that $N(\tfrac{1}{1-\mu} s) \leq c_\diamond^{-1}$ because $\tfrac{1}{1-\mu} s$ is in $\mathcal{T}$; and that $K(\tfrac{1}{1-\mu} s) \leq c_0 K(s)$ (which was explained in Step 1). Also keep in mind that the inequality in (9.1) holds for $\tfrac{1}{1-\mu} s$. And, keep in mind that $r^2 s^2 K(s)^2 \leq \kappa_{\diamond\mu}$ (because $\tfrac{1}{1-\mu} s \leq r_\diamond$). This list of facts leads directly from (9.8) to the bound

$$\int_{\partial B_s} |F_A| \leq c_\mu s (r s K(s)) .$$

(9.9)

More is said about the left hand side of this inequality in the next step.

<u>Step 3</u>: To exploit (9.9), write the identity $r r_\diamond K(r_\diamond) = \sqrt{\kappa_{\diamond\mu}}$ as

$$r s K(s) = \sqrt{\kappa_{\diamond\mu}} \, \tfrac{s}{r_\diamond} \, \tfrac{K(s)}{K(r_\diamond)} .$$

(9.10)

This leads to the bound $r s K(s) \leq c_\mu \tfrac{s}{r_\diamond}$. Use the latter bound with (9.9) to bound the $L^1$ norm of $|F_A|$ on $\partial B_s$ by $c_\mu \tfrac{s^2}{r_\diamond}$. Using this last bound for the integral of $|F_A|$ in (9.5) leads directly to (9.2)

**b) Proof of Proposition 7.4**

The proof has four parts. With regards to the proof, note that if (8.27) holds, then so does the second bullet of the proposition.

*Part 1*: If the whole of the set [r, 3r] is in $\mathcal{I}$, then what is said in Part 1 of Section 8e can be repeated almost word for word to obtain points $s_0$ and $s_1$ that obey the bullets of Proposition 7.4 and (8.27) (with $\rho_{\mu,\varepsilon} = c_\diamond^{-1}$). With the preceding understood, assume henceforth that r is in the set $\mathcal{I}$.

Define a subset $\mathcal{L}_\ddagger \subset \mathcal{I}$ by the rule where by $r \in \mathcal{L}_\ddagger$ when the whole interval [r, 4r] is in $\mathcal{L}$. To be sure, this definition says that the following conditions are obeyed by r:

- *The interval [r, 4r] is in $\mathcal{I}$.*
- $\int_{B_s} |F_A|^2 < c_\diamond^8 r^2 s^2 K(s)^2 N(s)$ *for all* $s \in [r, 4r]$.

(9.12)

If r is not in $\mathcal{L}_\ddagger$, then at least one of the bullets in (9.12) are violated. If it is the top bullet, then what is said in Part 2 of Section 8e can be repeated to obtain points $s_0$ and $s_1$ that obey the bullets of Proposition 7.4 and (8.27) (with $\rho_{\mu,\varepsilon} = c_\diamond^{-1}$). By the same token, what is said in Parts 4-6 of Section 8e (with Lemma 7.2 or the third bullet of (8.29)) can be used to verify the conclusions of Proposition 7.4 and (8.27) for suitable $\rho_{\mu,\varepsilon}$ when the top bullet of (9.12) holds, but the second bullet does not.

*Part 2*: Part 1 discussed the cases when one or both of the bullets in (9.12) are violated. This part and Parts 3 and 4 consider the remaining case, which is when both bullets in (9.12) hold and so r is in the set $\mathcal{L}_\ddagger$. To set the notation for what is to come, let $\kappa_\ddagger$ denote now the version of $\kappa_\mu$ that appears in Lemma 9.1.

This part and Parts 3 and 4 of the proof consider the event that r is a point from $\mathcal{L}_\ddagger$ and that r is strictly less than $\kappa_\ddagger^{-1} N(r) r_\diamond$. Let $\rho$ denote the largest number s in the closed interval between r and $\kappa_\ddagger^{-1} N(r) r_\diamond$ such that the whole of the interval [r, s) is in $\mathcal{L}$. Lemma 9.1 can be invoked to get bounds for K on [r, $\rho$] and to conclude that $N(\rho) \geq \frac{7}{8} N(r)$.

Now, there are three possibilities for the point $\rho$: Either the top bullet in (9.12) with $\rho$ replacing r is false, or the top bullet in (9.12) holds with $\rho$ replacing r but the lower bullet in (9.12) is false, or else $\rho = \kappa_\ddagger^{-1} N(r) r_\diamond$. These cases are considered in Parts 3, and 4 of the proof.

*Part 3*: If either the top bullet fails with $\rho$ replacing r, or the top bullet is true and the lower bullet is false (both with $\rho$ replacing r) then Part 1 of the proof can be applied with $\rho$ replacing r. And, this version of Part 1 finds two points $s_0 \geq 2\rho$ and $s_1 \geq s_0$ which are such that N is greater than the minimum of $\frac{1}{10,000} c_\diamond^{-1}$ and $4N(\rho)$; and

$$\kappa(s) \leq \left(\tfrac{s}{s_1}\right)^{L(s_1)/\kappa_{\diamond\mu}} \kappa(s_1) \quad \textit{for all } s \in [\rho, s_0].$$

(9.13)

The points $s_0$ and $s_1$ obey the first two bullets of Proposition 7.4 with input r. Indeed, the first bullet holds because $s_0$ is greater than 2r (because $\rho > r$); and the second bullet holds because neither $N(s_0)$ nor $N(s_1)$ is less than $4N(\rho)$ and $N(\rho)$ is not less than $\tfrac{7}{8} N(r)$. The points $s_0$ and $s_1$ also obey the third bullet of Proposition 7.1. This follows from (9.13) for $s \in [\rho, s_0]$); and it follows for points between r and $\rho$ by virtue of the second bullet of Lemma 9.1 and the $s = \rho$ version of (9.13).

*Part 4*: Supposing that $r < \kappa_{\ddagger}^{-1} N(r) r_\diamond$, define the number $\rho$ as in the Part 2, but assume here that $\rho = \kappa_{\ddagger}^{-1} N(r) r_\diamond$. If $r \geq \kappa_{\ddagger}^{-1} N(r) r_\diamond$, set $\rho = r$. In any event, the two bullets of Lemma 9.1 hold when $\rho > r$. These bullets also hold when $\rho = r$ being that they are tautological in this case.
Lemma 6.1 implies that

$$N(s) \geq \left(\tfrac{\rho}{s}\right)^2 N(\rho) \quad \textit{for } s \in [\rho, 4r_\diamond].$$

(9.13)

Noting that $\tfrac{\rho}{s} \geq \tfrac{1}{4} \kappa_{\ddagger}^{-1} N(r)$ for $s \in [\rho, 4r_\diamond]$, and noting that $N(\rho) \geq \tfrac{7}{8} N(r)$, it follows from (9.13) that $N(s) \geq c_\mu^{-1} N(r)^3$ for $s \in [\rho, 4r_\diamond]$. This is also obeyed for $s \in [r, \rho]$, by virtue of the second bullet in Lemma 9.1. Therefore, the first and second bullets of Proposition 7.4 are obeyed if $s_0$ and $s_1$ are taken to be $4r_\diamond$. The third bullet of Proposition 7.4 then follows from the second by virtue of (3.9).

## 10. Proof of Theorem 1.1

The tools to prove Theorem 1.1 have now been assembled; and so it is time to give the proof. By way of a reminder, the context is this: A sequence $\{(r_n, (A_n, a_n))\}_{n \in \mathbb{N}}$ is given with the sequence $\{r_n\}_{n \in \mathbb{N}} \subset (1, \infty)$, and with each $(A, a) \in \{(A_n, a_n)\}$ obeying the corresponding $r \in \{r_n\}_{n \in \mathbb{N}}$ version of (1.5). (Here and in what follows, $\mathbb{N}$ denotes the positive integers.) Given these sequences as input, the task is to find a subsequence $\Theta \subset \mathbb{N}$ so that the sequence $\{(A_n, a_n)\}_{n \in \Theta}$ converges in the manner described by Theorem 1.1. The event that the sequence $\{r_n\}_{n \in \mathbb{N}}$ has a bounded subsequence is restated as Proposition 2.2 and proved in Section 2c (thus reproducing a result by Ben Mares [Ma].) All that follows here assumes (implicitly for the most part) that the sequence $\{r_n\}_{n \in \mathbb{N}}$ is increasing and unbounded.

## a) $L^2_1$ and $L^\infty$ limits

Various aspects of Proposition 2.1 in [T1] have analogs for sequences of solutions to (1.5). They are restated as the following proposition:

**Proposition 10.1**: *There exists $\kappa > 1$ with the following significance: Let $\{r_n\}_{n\in\mathbb{N}} \subset (1,\infty)$ denote an increasing and unbounded sequence; and for each $n \in \mathbb{N}$, let $(A_n, a_n)$ denote a solution to the $r = r_n$ version of (1.5). There exists a subsequence $\Lambda \subset \mathbb{N}$ such that the bulleted items listed below hold.*

- *The sequences $\{\int_X (|d|a_n||^2 + |a_n|^2)\}_{n\in\Lambda}$ and $\{\sup_X |a_n|\}_{n\in\Lambda}$ are bounded by $\kappa$.*
- *The sequence $\{|a_n|\}_{n\in\Lambda}$ converges weakly in the $L^2_1$ topology and strongly in all $p < \infty$ versions of the $L^p$ topology.*
- *The limit function (it is denoted by $|v|$) is an $L^\infty$ function whose value can be defined at each point in X by the rule whereby $|v|(p) = \limsup_{n\in\Lambda} |a_n|(p)$ for each $p \in X$.*
- *The sequence $\{a_n \otimes a_n^\dagger\}_{n\in\Lambda}$ converges strongly in any $q < \infty$ version of the $L^q$ topology on the space of sections of $\mathbb{S}^+ \otimes \mathbb{S}^{+\dagger}$. The limit section is denoted by $v \otimes v^\dagger$ and its trace is the function $|v|^2$.*
- *Use $f$ to denote a given $C^0$ function.*
  i) *The sequences $\{\int_X f |F_{A_n}|^2\}_{n\in\Lambda}$ converge.*
  ii) *The sequences $\{\int_X f |\nabla_{A_n} a_n|^2\}_{n\in\Lambda}$ and $\{r_n^2 \int_X f |\langle a_n, \tau a_n\rangle|^2\}_{n\in\Lambda}$ converge. The limit of the first sequence is denoted by $Q_{\nabla, f}$ and that of the second by $Q_{\tau, f}$. These are such that*

$$\tfrac{1}{2}\int_X d^* df \, |v|^2 + Q_{\nabla, f} + Q_{\tau, f} + \int_X f \langle \mathfrak{R}, v \otimes v^\dagger \rangle = 0,$$

- *Fix $p \in X$ and let $G_p$ denote the Green's function with pole at $p$ for the operator $d^\dagger d + 1$. The sequence that is indexed by $\Lambda$ with n'th term being the integral of $G_p(|\nabla_{A_n} a_n|^2 + r_n^2 |\langle a_n, \tau a_n\rangle|^2)$ is bounded by $\kappa$. Let $Q_{\diamond, p}$ denote the lim-inf of this sequence of integrals. The function $|v|^2$ obeys the equation*

$$\tfrac{1}{2} |v|^2(p) + Q_{\diamond, p} = \int_X G_p (\tfrac{1}{2} |v|^2 - \langle \mathfrak{R}, v \otimes v^\dagger \rangle).$$

The proof of this proposition differs only cosmetically from Part II of the proof of Proposition 2.1 in [T] to which the reader is directed. The proof is not given here.

## b) Hölder continuity across the zero locus

The third bullet of the Proposition 10.1 defines the values of $|v|$ at each point in X. With this definition understood, let Z now denote the set of points where $|v| = 0$. This is to say that $p \in Z$ if and only if:

$$\limsup_{n \in \Lambda} |a_n|(p) = 0 \ .$$

(10.1)

The proposition below makes the formal assertion that Z is a closed set in X and that $|v|$ is uniformly Hölder continuous across Z.

**Proposition 10.2**: *There exists $\kappa > 1$ with the following significance: Let $\{r_n, (A_n, a_n)\}_{n \in \mathbb{N}}$ denote the input sequence for Proposition 10.1, and let $\Lambda$ denote the subsequence of $\mathbb{N}$ given in Proposition 10.1. Use the corresponding sequence $\{a_n\}_{n \in \Lambda}$ to define the function $|v|$ as instructed by Proposition 10.1; and then define Z to be its zero locus which is the set of points in X where (10.1) holds. The set Z is a closed set in X. Furthermore, if p is in Z, then $|v|$ on the radius $\kappa^{-1}$ ball centered at p obeys $|v| \le \kappa \mathrm{dist}(p, \cdot)^{1/\kappa}$.*

This proposition is seen momentarily to follow from the upcoming Proposition 10.3
    To set notation for Proposition 10.3, suppose for the moment that a point $p \in X$ has been chosen. Define the non-negative function $\mathcal{K}$ on $(0, c_0^{-1}]$ by the rule

$$r \to \mathcal{K}(r)^2 = \frac{1}{r^3} \int_{\partial B_r} |v|^2 \ .$$

(10.2)

This function is bounded because $|v|$ is bounded (see the third bullet in Proposition 10.1.)

**Proposition 10.3**: *There exists $\kappa > 1$ with the following significance: Let $\{r_n, (A_n, a_n)\}_{n \in \mathbb{N}}$ denote the input sequence for Proposition 10.1, and let $\Lambda$ denote the subsequence of $\mathbb{N}$ given in Proposition 10.1. Use the corresponding sequence $\{a_n\}_{n \in \Lambda}$ to define the function $|v|$ as instructed by Proposition 10.1; and then define Z to be its zero locus. If $p \in Z$ and if $r \in (0, \kappa^{-1}]$, then $\mathcal{K}(r) \le \kappa r^{1/\kappa}$.*

This proposition is proved in the Section 10c. Assume it to be true for the time being.

*Proof of Proposition 10.2*: The fact that Z is closed follows from the asserted Hölder norm bound for $|v|$ along Z. To prove this bound, first fix $p \in Z$, and then $r \in (0, c_0^{-1}]$. Having done this, let $\chi_{1/4}$ denote the $\mu = \frac{1}{4}$ version of the function $\chi_\mu$ that was introduced in Section 3c. (This version of $\chi_\mu$ is equal to 0 where the distance to p is greater than $\frac{27}{32} r$ and it is equal to 0 where the distance to $\mu$ is less than $\frac{26}{32}$ r). Let q denote a chosen point from $B_{r/2}$ and let $G_q$ denote the Greens function for the operator $d^*d+1$ on X with pole at the point q. Set $f$ to be the function $(1-\chi_{1/4})G_q$. Use this in Item ii) of the fifth bullet of Proposition 10.1 and use q in lieu of p in the sixth bullet of Proposition 10.1. Subtracting the fifth bullet's identity from sixth bullet's identity leads to an inequality for $|v|^2(q)$ that takes the form

$$\tfrac{1}{2}|v|^2(q) \le \int_X \chi_{1/4} G_q (\tfrac{1}{2}|v|^2 - \langle \mathfrak{R}, v\otimes v^\dagger\rangle) + \int_X (\tfrac{1}{2}d^*d\chi_{1/4}G_p - \langle d\chi_{1/4}, dG_q\rangle)|v|^2 .$$

(10.3)

This in turn leads to the bound

$$|v|^2(q) \le c_0 r^2 + c_0 \tfrac{1}{r^4} \int_{B_r} |v|^2 .$$

(10.4)

By way of an explanation, the $c_0 r^2$ term in (10.4) accounts for the left most integral on the right hand side of (10.3). This is because $|v|\le c_0$ and because

$$|G_q| \le c_0 \tfrac{1}{\mathrm{dist}(\cdot,q)^2} .$$

(10.5)

Meanwhile, the right most integral on the right hand side of (10.3) is accounted for by the product of $c_0 \tfrac{1}{r^4}$ times the integral of $|v|^2$ that appears on the right hand side of (10.5). Indeed, the factor of $\tfrac{1}{r^4}$ accounts for the fact that $|d\chi_{1/4}| \le c_0 \tfrac{1}{r}$ and that $|\nabla d\chi_{1/4}| \le c_0 \tfrac{1}{r^2}$; and the fact that $G_q \le c_0 \tfrac{1}{r^2}$ and $|\nabla G_q| \le c_0 \tfrac{1}{r^3}$ on the support of $|d\chi_{1/4}|$.

It is a consequence of Proposition 10.3 that the right hand side of (10.4) is bounded by $c_0 r^{1/c_0}$ when $r \le c_0^{-1}$; and the latter bound implies the asserted Hölder norm bound in Proposition 10.2

**b) Lemmas about N and K as defined by sequence in X**

To set some notation, fix $p \in X$ and a number $s \in (0, c_0^{-1}]$ Supposing that an integer $n \in \Lambda$ have been chosen, then $K_n(p; \cdot)$ and $N_n(p; \cdot)$ are used to denote the versions of the respective functions K and N that are defined by the data $(A_n, a_n)$ and the point p. In the event that there is little likelihood of confusion with regards to the point from X, these are denoted by $K_n$ and $N_n$.

**Lemma 10.4**: *There exists $\kappa > 1$ with the following significance: Let $\{r_n, (A_n, a_n)\}_{n\in\Lambda}$ denote the input data for Proposition 10.1. Fix $p \in X$ and suppose there are numbers* $\gamma \in (0, 1)$ and $\rho \in (0, \kappa^{-1}]$ and $r_\diamond \in [0, \kappa^{-1}\rho)$; a subsequence $\Xi \subset \Lambda$; and sequences $\{q_n\}_{n \in \Xi} \subset X$ and $\{r_n\}_{n \in \Xi} \subset [0, \kappa^{-1}]$; all with the following properties:
- The sequence $\{q_n\}_{n \in \Xi}$ converges to $p$
- The sequence $\{r_n\}_{n \in \Xi}$ converges to $r_\diamond$.
- For each $n \in \Xi$, and $\kappa_n(q_n; s) \leq s^\gamma$ for all $s \in (r_n, \rho)$.

Then the point $p$'s version of the function $\mathcal{K}$ obeys $\mathcal{K}(s) \leq s^\gamma$ for all $s \in (r_\diamond, \kappa^{-1}\rho)$.

Lemma 10.4 will be proved momentarily. Here is a corollary:

**Lemma 10.5**: *Given $\mu \in (0, \frac{1}{100})$, there exists $\kappa_\mu > 100$, and given also $\varepsilon \in (0, \kappa_\mu^{-1})$, there exists $\rho_{\mu,\varepsilon} \in (0, \kappa_\mu^{-1})$ with the following significance: Let $\{r_n, (A_n, a_n)\}_{n \in \Lambda}$ denote the input data for Proposition 10.1. Fix $p \in X$ and suppose there exits $r_\diamond \in [0, \rho_{\mu,\varepsilon})$; a subsequence $\Xi \subset \Lambda$; and sequences $\{q_n\}_{n \in \Xi} \subset X$ and $\{r_n\}_{n \in \Xi} \subset [0, \rho_{\mu,\varepsilon}]$; all with the following properties*:
- *The sequence $\{q_n\}_{n \in \Xi}$ converges to $p$*
- *The sequence $\{r_n\}_{n \in \Xi}$ converges to $r_\diamond$.*
- *For each $n \in \Xi$, the number $N_n(q_n; r_n)$ obeys $N_n(q_n; r_n) \geq \varepsilon$.*

*Then the point $p$'s version of the function $\mathcal{K}$ obeys $\mathcal{K}(s) \leq s^{\varepsilon^3/\kappa_\mu}$ for all $s \in (r_\diamond, \kappa_\mu^{-1}\rho_{\mu,\varepsilon})$.*

Keep in mind that both Lemmas 9.4 and 9.5 can be applied when $r_\diamond = 0$ and when some or all of the $r_n$ are also 0. They also hold in the simple case when all $q_n$ are equal to $p$. Note in addition that the bound in these lemmas for $\mathcal{K}$ holds for the endpoint $r_\diamond$ (if $r_\diamond > 0$) and at the top endpoint of the relevant interval because the function $\mathfrak{h}(f, \cdot)$ is continuous on the interval $(0, c_0^{-1})$ for any given function $f \in L^2_1(X)$; and, in particular, for $f = |\nu|$.

*Proof of Lemma 10.5*: An appeal to the point $q_n$ version of Proposition 7.1 can be made using as input $r = r_n$. The result is a bound of the form

$$\kappa_n(q_n; s) \leq s^{\varepsilon^3/\kappa_\mu} \quad \text{for all } s \in [r_n, \rho_{\mu,\varepsilon}]$$

(10.6)

Now invoke Lemma 9.4.

*Proof of Lemma 10.4*: The proof of has two parts.

    *Part 1*: This part of the proof is a digression for some background. To start, fix a point $p \in X$ and $r \in (0, c_0^{-1}]$. Let $B_r(p)$ denote the ball of radius $r$ centered at $p$ and let $\partial B_r(p)$ denote its boundary sphere. A quadratic function on $C^\infty(X)$ (the Frechet space of smooth functions on $X$) is defined by the rule

$$f \to \mathfrak{h}_{(p,r)}(f) = \int_{\partial B_r(p)} f^2 \ .$$

(10.7)

It is a consequence of standard Sobolev theorems that $\mathfrak{h}_{(p,r)}$ extends to the Hilbert space of $L^2_1$ functions as a bounded, compact function. The extension as a bounded function follows from (2.1) using an integration by parts (in reverse) to write $\mathfrak{h}_{(p,r)}(f)$ as an integral over $B_r$. The fact that $\mathfrak{h}_{(p,r)}$ is a compact is harder to prove (see e.g. [AF].)

To say that $\mathfrak{h}_{(p,r)}$ is compact means the following: If $\Lambda \subset \mathbb{N}$ is an unbounded set and if $\{f_n\}_{n\in\Lambda}$ converges weakly in the $L^2_1$ topology to a function (to be denoted by $f$), then the sequence of numbers $\{\mathfrak{h}_{(p,r)}(f_n)\}_{n\in\Lambda}$ converges to $\mathfrak{h}_{(p,r)}(f)$. Moreover, this convergence is uniform with respect to p and r in the following sense: Given positive numbers $\varepsilon$ and $\rho$ (with $\rho < c_0^{-1}$], there exists $N_{(\varepsilon,\rho)}$ such that when $n > N_{(\varepsilon,\rho)}$, then

$$|\mathfrak{h}_{(p,r)}(f) - \mathfrak{h}_{(p,r)}(f_n)| < \varepsilon \quad \text{for all } p \in X \text{ and } r \in (\rho, c_0^{-1}].$$

(10.8)

There is a related point to make about $\mathfrak{h}_{(p,r)}$ (which is also a consequence of the Sobolev inequalities in (2.1)): The function $\mathfrak{h}_{(p,r)}$ is jointly continuous when viewed as a function on the product space $X \times (0, c_0^{-1}] \times L^2_1(X)$. To say more about the joint continuity, fix $r \in (0, c_0^{-1}]$ and $\tau \in [-\frac{1}{8}, \frac{1}{8}]$. Suppose that p and q are points in X with distance at most $\frac{1}{100} r$ between them. Let $\Delta$ denote this distance from p to q. Under these circumstances, the point p is well inside the radius $(1+\tau)r$ ball centered at q (and q is well inside the radius r ball centered at p.) Therefore, pushing points in or out along the geodesic arcs from p starting on the radius r sphere centered at p defines an isotopy that moves the radius r sphere centered at p to the radius $(1+\tau)r$ sphere centered at q. Using the coordinates of this isotopy with an integration by parts leads to the following observation: Supposing that $f$ is an $L^2_1$ function on X with $L^2_1$ norm equal to 1, then

$$\left| \int_{\text{dist}(p,\cdot)=r} f^2 - \int_{\text{dist}(q,\cdot)=(1+\tau)r} f^2 \right| \le c_0 r^{-1}(\Delta + |\tau|) \ .$$

(10.9)

*Part 2*: Appling the inequality in (10.9) with $f = |a_n|$ and $q = q_n$ and $s \in (r_\diamond, \rho]$ leads to the conclusion that

$$K_n(p; s)^2 \le K_n(q_n, s)^2 + c_0 r^{-4} \text{dist}(p, q_n) \ .$$

(10.10)

Let $B_s$ denote the radius s ball centered at p. The inequality in (10.10) and the assumption in the third bullet of the lemma lead to the bound

$$\tfrac{1}{s^3} \int_{\partial B_s} |a_n|^2 \leq (1+c_0 s^2) s^{2\gamma} + c_0 \tfrac{1}{s^4} \operatorname{dist}(p, q_n) \ .$$

(10.11)

Granted this bound, it then follows that

$$\lim_{n \in \Xi} \tfrac{1}{s^3} \int_{\partial B_s} |a_n|^2 \leq (1+c_0 s^2) s^{2\gamma}$$

(10.12)

because of the assumption in the top bullet if the lemma.

The assertion of Lemma 10.4 follows because (10.12) holds for any $s \in (r_\diamond, \rho)$.

### c) Proof of Proposition 10.3

The proof of the proposition is given in 3 parts.

*Part 1*: Fix $\mu \in (0, \tfrac{1}{100})$ and then fix $c > c_\mu$ so that the conclusions of Propositions 3.1-3.3 and 4.1 and 7.1 hold for these choices of $\mu$ and $c$. Let $\kappa_\ddagger$ denote in what follows 10,000 times the maximum of the versions of $\kappa_\mu$ that appear in Propositions 3.1-3.3 and 4.1 and 7.1, and in Lemma 10.5.

Fix some small $\delta \in (0, \tfrac{1}{1,000,000})$ about which more will be said below when a specific choice on the order of $\mu^3$ is made. Given $\delta$ and given an integer $n \in \Lambda$ and given $q \in X$, let $r_{\ddagger n}(q)$ denote greatest lower bound of the numbers $s \in (0, \kappa_\ddagger^{-1} c^{-1} \delta^2]$ such that $N_n(q; s) = \kappa_\ddagger^{-1} c^{-2} \delta^3$ (if such a number exists). If $N_n(q; s)$ is less than $\kappa_\ddagger^{-1} c^{-2} \delta^3$ on all of the interval $(0, \kappa_\ddagger^{-1} c^{-1} \delta^2]$, then set $r_{\ddagger n}(q)$ to be the upper bound for the interval, thus $\kappa_\ddagger^{-1} c^{-1} \delta^2$. Denote this special value $\kappa_\ddagger^{-1} c^{-2} \delta^3$ for the function $N_n$ by $\varepsilon$; its dependence on $\mu$, $c$ and $\delta$ is implicit. (Looking ahead, these definitions will facilitate an appeal to Proposition 4.1.)

It is a consequence of Lemma 10.5 that the conclusions of Proposition 10.3 hold for the chosen point p (with the number $\kappa$ being $\varepsilon^{-3} \kappa_\mu$ with $\kappa_\mu$ from Lemma 10.5) if there is a subsequence $\Xi \subset \Lambda$ such that $\lim_{n \in \Xi} r_{\ddagger n}(p) = 0$. (Take each $q_n$ in the lemma to be p and take each $r_n$ to be $r_{\ddagger n}$.) In fact, slightly more is true: It is also a consequence of Lemma 10.5 that the conclusions of Proposition 10.3 hold (with the same $\kappa$) if there is a subsequence $\Xi \subset \Lambda$ and a convergent sequence $\{q_n\}_{n \in \Xi}$ with limit p such that $\lim_{n \to \infty} r_{\ddagger n}(q_n) = 0$. (One can take $r_n$ in Lemma 10.5 to be $r_{\ddagger n}(q_n)$.)

With the preceding understood, the remaining parts of the proof make the following assumption:

*Given $\delta \in (0,1)$, there exists $\rho_\delta > 0$ such that if $q \in X$ has $\operatorname{dist}(q,p) < \rho_\delta$ then $r_{\ddagger n}(q) > \varepsilon$ for all $n \in \Lambda$ (equivalently, $N_n(q; \cdot) < \varepsilon$ on $(0, \rho_\delta)$ for all $n \in \Lambda$).*

(10.13)

(If this assumption is violated, then there is a subsequence $\Xi$ and a sequence $\{q_n\}_{n \in \Xi}$ of the sort described in the preceding paragraph.)

*Part 2*: Supposing that (10.13) holds, then (of course) the sequence $\{r_{\ddag n}(p)\}_{n \in \Lambda}$ is bounded away from zero. Let $r_{\ddag}$ denote the lim-inf of this sequence, a positive number. For each $n \in \Theta$, let $r_{cFn}$ denote the $(A_n, a_n)$ version of the number $r_{cF}$ that is defined in Section 3a. Suppose in this part of the proof that $\Lambda$ has a subsequence that labels a corresponding subsequence in $\{r_{cFn}\}_{n \in \Lambda}$ with a positive lower bound. Fix such a subsequence (to be denoted by $\Xi$) with the property that $\{r_{cFn}\}_{n \in \Xi}$ converges to its non-zero limit. Let $r_{cF\ddag}$ denote this positive limit; and then let $r_\diamond$ denote the minimum of $r_{\ddag}$ and $r_{cF\ddag}$. This number $r_\diamond$ is such that

$$\lim_{n \in \Xi_p} K_n((1 - 2\mu)r_\diamond) = 0 \tag{10.14}$$

because there would otherwise be conflict between (10.1) and Proposition 3.3. Meanwhile, if (10.14) is to occur, then $\lim_{n \in \Xi} K_n(r_{\ddag n}(p)) = 0$ because of (3.9) and because $N_n(p; s) \le \varepsilon$ for $s \le r_{\ddag n}(p)$. This implies that $\mathcal{K}(s) = 0$ for $s \le r_{\ddag}$. Meanwhile, Lemma 10.5 can be invoked with all $q_n$ equal to $p$ and with each $r_n$ equal to $r_{\ddag n}(p)$ to see that

$$\mathcal{K}(s) \le s^{\varepsilon^3/\kappa_\ddag} \text{ for all } s \in [r_\ddag, \kappa_\ddag^{-1}\rho_{\mu,\varepsilon}] \tag{10.15}$$

with $\rho_{\mu,\varepsilon}$ coming from Lemma 10.5. Since $\mathcal{K}(s) = 0$ for $s \le r_\ddag$, the inequality in (10.15) holds for all $s \in (0, c_\ddag^{-1}]$, and thus Proposition 10.3 holds for p if $\kappa$ is greater than $\kappa_\ddag \delta^{-3}$.

*Part 3*: The last case to imagine is the case when $\{r_{\ddag n}(p)\}_{n \in \Lambda}$ has a positive lower bound and $\{r_{cFn}\}_{n \in \Lambda}$ has limit zero. Although this last case can be imagined, it can not occur. The next lemma makes an official statement to this effect. To set the notation, suppose for the moment that q is any given point in X. Let $r_{cFn}(q)$ denote the version of the number $r_{cF}$ that is defined using the point q and the pair $(A_n, a_n)$ (see (3.1).)

**Lemma 10.6**: *There exists $\kappa > 100$, and given $\mu \in (0, \frac{1}{100})$ there exists $\kappa_\mu > \kappa$ with the following significance: Let $\{r_n, (A_n, a_n)\}_{n \in \Lambda}$ denote the input data for Proposition 10.1. Fix $c > \kappa_\mu$ and $\delta \in (0, \kappa^{-1}\mu^3)$; and, given $n \in \Lambda$, use these numbers with $(A_n, a_n)$ to define the functions $r_{\ddag n}(\cdot)$ and $r_{cFn}(\cdot)$. Suppose that $p \in X$ and that there exists $\rho > 0$ such that each $n \in \Lambda$ version of $r_{\ddag n}(\cdot)$ is greater than $\rho$ on the ball of radius $\rho$ centered at p. Then, the sequence $\{r_{cFn}(p)\}_{n \in \Lambda}$ is bounded away from zero.*

*Proof of Lemma 10.6*: Suppose that the assertion is false to generate nonsense. To this end, assume that there is a subsequence $\Xi \subset \Lambda$ such that $\{r_{cFn}\}_{n\in\Xi}$ converges to zero. The argument to generate nonsense from this assumption has five steps.

Step 1: According to the assumptions of the lemma, if $q \in X$ and $\mathrm{dist}(q, p) < \rho$, then $r_{\ddagger n}(q) \geq \rho$. Thus, supposing $\{q_n\}_{n\in\Xi}$ is a convergent sequence with limit p, the sequence $\{r_{\ddagger n}(q)\}$ is bounded away from zero. By way of a contrast, the sequence $\{r_{cFn}(q_n)\}_{n\in\Xi}$ must converge to zero, for otherwise the sequence $\{r_{cFn}\}_{n\in\Xi}$ would have a positive lower bound. To be precise,

$$r_{cFn}(q) \leq r_{cFn} + \mathrm{dist}(p, q_n)$$

(10.16)

so as not to foul the definition of $r_{cFn}$ (which is $r_{cFn}(p)$). Therefore, a finite number elements from $\Xi$ can be ejected if necessary with the result being a subsequence (to be denoted henceforth by $\Xi$ also) obeying

$$r_{cFn}(q_n) < \tfrac{1}{100} r_{\ddagger n}(q_n) \;\; \text{for all } n \in \Xi .$$

(10.17)

Step 2: Because $r_{cFn}(q_n)$ is so small, Proposition 4.1 can be brought to bear when n is large using the number $(1-\mu) r_{cFn}(q_n)$ for the value of $r_c$ and using $\delta$ for Proposition 4.1's version of $\varepsilon$. By way of a reminder, this version of Proposition 4.1 says that

$$\int_{\mathrm{dist}(\cdot, q_n) \leq (1-\mu)^2 r_{cFn}(q_n)} |F_{A_n}|^2 < \delta c^{-2}.$$

(10.18)

Meanwhile, the whole of the radius $r_{cFn}(q_n)$ ball centered at $q_n$ can covered by a union of $N+1$ balls with $N \leq c_0 \mu^{-3}$; the first ball being the radius $(1-\mu)^2 r_{cFn}(q_n)$ ball centered at $q_n$ (this is the domain of integration in (10.18)) and the other N balls having radius $3\mu \, r_{cFn}(q_n)$ and centers on the boundary of the radius $r_{cFn}(q_n)$ ball centered at $q_n$. Let $\vartheta_{n1}$ denote the set centers of these N other balls. If Proposition 4.1 can be applied to every point in $\vartheta_{n1}$ with the input number $r_c$ being larger than $4\mu \, r_{cFn}(q_n)$, then the integral of $|F_{A_n}|^2$ over each radius $4\mu r_{cFn}(q_n)$ ball with center in $\vartheta_{n1}$ will likewise be less than $\delta c^{-2}$. This event would lead to the bound

$$\int_{\mathrm{dist}(\cdot, q_n) \leq r_{cFn}(q_n)} |F_{A_n}|^2 \leq c_0 \mu^{-3} \delta c^{-2}$$

(10.19)

which is nonsense if $\delta < \frac{1}{100} c_0^{-1} \mu^3$ because it runs afoul of the definition in (3.1) of the number $r_{cF}$ for the point $q_n$.

Step 3: Therefore, assuming henceforth that $\delta$ is less than $\frac{1}{100} c_0^{-1} \mu^3$, then the following is true: The set $\vartheta_{n1}$ as defined in Step 2 has at least one point whose version of Proposition 4.1 can not use the number $4\mu r_{cFn}(q_n)$ as the input value of r. Supposing that $q_{1n}$ is such a point, this means that one or both of the following occur:

- *The value of $q_{1n}$'s version of the function $N_n$ at $r = 4\mu r_{cFn}(q_n)$ is not less than ε.*
- *The value of $q_{n1}$'s version of $r_{cF}$ is less than $4\mu r_{cFn}(q_n)$.*

(10.20)

Now, the top bullet in (10.20) can not occur as soon as n is large enough so that

$$r_{cFn} < \tfrac{1}{100} \rho \quad and \quad dist(p, q_n) < \tfrac{1}{100} \rho .$$

(10.21)

This is because $dist(q_n, q_{n+1}) = r_{cFn}(q_n)$ and (10.15) and the assumption about $r_{\pm n}(\cdot)$ being greater than ρ on the radius ρ ball centered at p. Therefore, it must be that the lower bullet occurs for large n. This lower bullet is says that

$$r_{cFn}(q_{1n}) < 4\mu r_{cFn}(q_n)$$

(10.22)

(which in turn is less than $4\mu$ time the sum of $r_{cFn}$ and the distance from p to $q_n$).

Step 4: Let $\{q_n\}_{n\in\Xi}$ denote the sequence from the previous steps. Each $q_n$ has a corresponding point $q_{n1}$, so there is a new sequence $\{q_{1n}\}_{n\in\Xi}$ which converges to p also. This is because

$$dist(q_{1n}, p) \le r_{cFn}(q_n) \le r_{cFn} + dist(p, q_n) \quad \textit{for each } n \in \Xi.$$

(10.23)

But, now the arguments in the previous steps can be applied to the sequence $\{q_{1n}\}_{n\in\Xi}$ with the result being yet another sequence $\{q_{2n}\}_{n\in\Xi}$ with $r_{cFn}(q_{2n}) < (4\mu)^2 r_{cFn}(q_n)$ and with

$$dist(q_{2n}, p) \le (1+4\mu) r_{cFn} + dist(p, q_n) \quad \textit{for each } n \in \Xi .$$

(10.24)

And, these same arguments can be again invoked with the $\{q_{2n}\}_{n\in\Xi}$ sequence to generate a sequence $\{q_{3n}\}_{n\in\Xi}$, and so on generating successively sequences $\{q_{4n}\}_{n\in\Xi}, \ldots, \{q_{kn}\}_{n\in\Xi}, \ldots$ for ever that obey, for each n and k,

- $r_{cFn}(q_{kn}) \le (4\mu)^k r_{cFn}(q_n)$ ,
- $\text{dist}(q_{kn}, p) \le (1 + 4\mu + \cdots + (4\mu)^k) r_{cFn} + \text{dist}(p, q_n) \le \frac{1}{1-4\mu} r_{cFn} + \text{dist}(p, q_n)$ .

(10.25)

Step 5: As explained directly, the implication that the iteration process can go from $q_{kn}$ to $q_{(k+1)n}$ for any $k \ge 1$ such that (10.25) holds is the desired nonsense that proves Lemma 10.6. Indeed, this iteration can't go on forever because $A_n$ is, after all, a smooth connection on X. This implies, in particular, there is some *positive* $\mathfrak{r}_n$ such that $r_{cFn}(q) > \mathfrak{r}_n$ for all $q \in X$. As a consequence of this and the top bullet in (10.25), the step going from $q_{(k+1)n}$ from $q_{kn}$ is nonsensical when k is such that $(4\mu)^{k+1} r_{cFn}(q_n) < \mathfrak{r}_n$.

### d) Pointwise convergence of $\{|a_n|\}_{n \in \Lambda}$ on X−Z

The pointwise convergence sequence of Proposition 10.1's sequence $\{|a_n|\}_{n \in \Lambda}$ is the topic of the following proposition.

**Proposition 10.7**: *Let $\{r_n, (A_n, a_n)\}_{n \in \Lambda}$ denote the input sequence from Proposition 10.1. The resulting function $|v|$ from Proposition 10.1 is continuous on X and the corresponding sequence of functions $\{|a_n|\}_{n \in \Lambda}$ converges to $|v|$ in the $C^0$ topology.*

The remainder of this subsection is dedicated to the proof of this proposition. This will come momentarily. The lemma that follows directly plays a role in the proof, and it plays a role in subsequent arguments also. To set the stage for the lemma, fix $\mu \in (0, \frac{1}{100})$ and then let $\kappa_{\ddagger}$ denote the minimum of the versions of $\kappa_\mu$ that appear in Propositions 3.1-3.3 and 4.1 and in Lemmas 10.5 and 10.6. Fix $c > \kappa_{\ddagger}^2$ and $\delta \in (0, \kappa_{\ddagger}^{-1}\mu^3)$. With $\mu$, $c$ and $\delta$ chosen, and supposing that $n \in \Lambda$, define the functions $r_{\ddagger n}(\cdot)$ and $r_{cFn}(\cdot)$ using the instructions from Parts 1 and 3 of the previous subsection.

**Lemma 10.8**: *Let $\{r_n, (A_n, a_n)\}_{n \in \Lambda}$ denote the sequence from Proposition 10.1. The corresponding sequence of functions $\{r_{\ddagger n}(\cdot)\}_{n \in \Lambda}$ and $\{r_{cFn}(\cdot)\}_{n \in \Lambda}$ are uniformly bounded away from zero on compact sets in X−Z.*

*Proof of Lemma 10.8*: It is a consequence of Lemma 10.5 that $\{r_{\ddagger n}(\cdot)\}_{n \in \Lambda}$ is bounded from below on any compact set in X−Z. Granted that this is so, then Lemma 10.6 implies that $\{r_{cFn}(\cdot)\}_{n \in \Lambda}$ is also bounded from below on any compact set in X−Z.

*Proof of Proposition 10.7*: The proof has six parts. The first 5 parts prove that $\{|a_n|\}_{n \in \Lambda}$ converges pointwise to $|\nu|$ on X−Z (and it follows from Proposition 10.3 that this is true on Z). The remaining part of the proof explains why $|\nu|$ is continuous and why the convergence of $\{|a_n|\}_{n \in \Lambda}$ to $|\nu|$ is in the $C^0$ topology.

*Part 1*: Given $p \in$ X−Z, fix a ball centered at p so that its closure is in X−Z. Let $r_{\ddagger\ddagger}$ denote a positive lower bound for all of the functions $\{r_{\ddagger n}(\cdot)\}_{n \in \Lambda}$ on this ball; and let $r_{cF\ddagger}$ denote a positive lower bound for all of the functions $\{r_{cFn}(\cdot)\}_{n \in \Lambda}$ on this ball. (These bounds exist courtesy of Lemma 10.8.) Fix a number to be denoted by $r_c$ that is less than half the radius of this same ball, and less than half smallest of the numbers $c^{-1}$, $r_{\ddagger\ddagger}$ and $r_{cF\ddagger}$. It follows from the definition of $r_{cFn}(p)$ that the sequence of curvatures $\{F_{A_n}\}_{n \in \Lambda}$ has $L^2$ norm bounded by $c^{-1}$ on the radius $2r_c$ ball centered at p. (Moreover, Proposition 4.1 says that these $L^2$ norms on the radius $2(1-\mu)r_c$ ball centered at p are bounded by $\delta^{1/2} c^{-1}$.) Meanwhile, the definition of $r_{\ddagger n}$ is such that the function $N_n(p; \cdot)$ is less than $\kappa_{\ddagger}^{-1} c^{-2} \delta^3$ on the interval $(0, 2r_c)$.

*Part 2*: If $q \in X$, then the sequence $\{K_n(q; r)\}_{n \in \Lambda}$ converges for any fixed $r > 0$. Let $K_*(q; r)$ denote the limit. Keep in mind that this limit is uniform given a positive lower bound on r; and it is also uniform with respect to variations of q. (See Part 1 of the proof of Lemma 10.4). For future reference, this means the following: Given numbers $\varepsilon > 0$ and $\rho \in (0, c^{-1}]$, there exists a positive integer $N_{(\varepsilon, \rho)}$ such that if $n \geq N_{(\varepsilon, \rho)}$, then

$$|K_n(q; r) - K_*(q; r)| < \varepsilon \text{ for any } q \in X \text{ and } r \in [\rho, c_0^{-1}].$$

(10.26)

Supposing that $p \in$ X−Z, let $\Xi_p \subset \Lambda$ denote a subsequence with the property that

$$\lim_{n \in \Xi_p} |a_n|(p) = \text{lim-sup}_{n \in \Lambda} |a_n|(p) \quad \textit{(which is } |\nu|(p).)$$

(10.27)

Since any $n \in \Lambda$ version of $K_n(p; \cdot)$ is increasing, it follows as a consequence of (10.27) that the number $K_*(p; r)$ for any $r \in (0, c_0^{-1}]$ is no smaller than $\sqrt{2\pi}|\nu|(p)$. Meanwhile, it follows from Proposition 3.3 applied to the pairs (A, *a*) from $\{(A_n, a_n)\}_{n \in \Xi_p}$ that

$$|\nu|(p) \geq (1 - \kappa_{\ddagger} c^{-1}) \frac{1}{\sqrt{2\pi}} K_*(p; r) \quad \textit{if } r \in (0, r_c]$$

(10.28)

These bound on $K_*(p; r)$ for $r \in (0, r_c]$ and (10.26) lead to the following observation: Given any small $\varepsilon > 0$ and $r \in (0, r_c]$, there exists a positive $N_{(\varepsilon, r)}$ such that the $n \geq N_{(\varepsilon,r)}$ versions of $K_n(p;r)$ differ from $\sqrt{2\pi}|\nu|(p)$ by at most $c_0 \kappa_{\ddagger} c^{-1}$. Since $c^{-1}$ is (by assumption) less that $\kappa_{\ddagger}^{-2}$, this says that $K_n(p;r)$ differs from $\sqrt{2\pi}|\nu|(p)$ by a small fraction of $|\nu|(p)$.

The bounds from the preceding paragraph lead in turn to a positive lower bound for any given $n \in \Lambda$ version of $|a_n|(p)$ when n is large because these bounds and Proposition 3.3 when $(A, a) = (A_n, a_n)$ and $n \geq N_{\varepsilon,r}$ imply the following: If q is a point in X with distance less than $(1-\mu)r$ from p, then

$$\left||a_n|(q) - |\nu|(p)\right| \leq c_0(\varepsilon + \kappa_{\ddagger} c^{-1})|\nu|(p)$$

(10.29)

In particular, this holds for $q = p$.

*Part 3*: Fix for the moment $\sigma \in (0, \tfrac{1}{1000}]$. If $n \in \Lambda$. If $s \in (0, r_c]$ and if $N_n(p,s) \geq \sigma$, then (according to Proposition 7.1)

$$K_n(p;s) \leq (\tfrac{s}{r_c})^{\sigma^3/\kappa_{\ddagger}} K_n(p;r_c) \ .$$

(10.30)

Noting that $K_n(p,s) \geq \sqrt{2\pi}|a_n|(p)$, and noting (10.29) for $r = r_c$, this inequality implies that

$$N_n(p;s) \leq c_0^{-1} \kappa_{\ddagger}^{1/3} (\ln \tfrac{r_c}{s})^{-1/3} \ .$$

(10.31)

This says in particular, that the function $N_n(p;\cdot)$ limits to zero as s limits to zero uniformly with regards to the sequence label n.

*Part 4*: Consider, in light of (10.31), what is said by Proposition 4.1 about the $L^2$ norm of $F_{A_n}$ on the radius s ball centered at p for $s \in (0, r_c]$ In particular, supposing that s is such that the right hand side of (10.31) is less than $\kappa_{\ddagger}^{-1} c^{-2}$, then Proposition 4.1 can be invoked with $\tfrac{1}{1-\mu} s$ used for its version of $r_c$ and with its version of $\varepsilon$ given by

$$\varepsilon(s) = c_0 \kappa_{\ddagger}^{4/9} c^{2/3} (\ln \tfrac{r_c}{s})^{-1/9}$$

(10.32)

because (10.31) says that $N_n(p;s) \leq \kappa_{\ddagger}^{-1} c^{-2} \varepsilon(s)^3$. This results in an $L^2$ bound saying that

$$\int_{B_s} |F_{A_n}|^2 < c_0 \kappa_{\ddagger}^{4/9} c^{-4/3} (\ln \tfrac{r_c}{s})^{-1/9}$$

(10.33)

This says in particular that the function on $(0, c_\mu^{-1} r_c)$ whose value at s is the $L^2$ norm of $F_{A_n}$ on the ball of radius s centered at p limits to zero as s limits to zero, and that this limiting process is uniform with respect to the label $n \in \Lambda$.

*Part* 5: Denote the right hand side of (10.33) by $c(s)^{-2}$ thus defining a function, $c(\cdot)$ on the radius $c_\mu^{-1} r_c$ ball centered at p. Note in particular that $\lim_{s \to 0} c(s)^{-1} = 0$. With $c(s)$ understood, what is said leading up to (10.29) can be repeated with $c$ replaced by $c(s)$ to obtain the following: Given $\varepsilon > 0$ and supposing that $n \geq N_{(\varepsilon,s)}$, then

$$\left||a_n|(q) - |v|(p)\right| \leq c_0 (\varepsilon + \kappa_\ddagger c(s)^{-1}) |v|(p) \quad \text{if } q \in B_{(1-\mu)s} .$$
(10.34)

Since $c(s)^{-1}$ limits to zero as s limits to zero, and since any value of s less than $c_0 r_c$ is fair game, this imples that $\lim_{n \in \Lambda} |a_n|(p) = |v|(p)$. Thus, the lim-inf and the lim-sup of the sequence $\{|a_n|(p)\}_{n \in \Lambda}$ are identical and the sequence $\{|a_n|(p)\}_{n \in \Lambda}$ converges to p.

*Part* 6: Having established the pointwise convergence of $|a_n|$ to p, it remains yet to prove that $|v|$ is continuous and that $\{|a_n|\}_{n \in \Lambda}$ converges to $|v|$ in the $C^0$ topology. The continuity of $|v|$ is implied by (10.34) and its counterpart with the roles of p and q reversed. (If q is in the radius $(1-\mu)s$ ball centered at p, then p is in the radius $(1-\mu)s$ ball centered at q and therefore so (10.34) holds with the roles of p and q reversed. Note in this regard that then number $N_{(\varepsilon,r)}$ does not depend on either p or q.)

These various versions of (10.34) also imply that the convergence of $\{|a_n|\}_{n \in \Lambda}$ to $|v|$ is $C^0$ convergence on compact subsets of X–Z. What with the Hölder bounds for $|v|$ along Z, this implies the $C^0$ convergence of $\{|a_n|\}_{n \in \Lambda}$ to $|v|$ on the whole of X.

**e) Convergence on small balls in X–Z**

The convergence of the input sequence for Proposition 10.1 on X–Z is the topic of this subsection. The upcoming Proposition 10.9 talks about the convergence on a given, small radius ball in X–Z of a sequence of pairs of connection on the product $\mathbb{C}$ bundle and section of $\mathbb{S}^+$ that is obtained from Proposition 10.1's sequence using a corresponding sequence of isomorphisms between E and the product $\mathbb{C}$ bundle. (Keep in mind that X–Z is an open subset of X (because Z is closed) which means that each point in X–Z is contained in a ball that lies entirely in X–Z.) The notation in what follows uses $\theta_0$ to denote the product connection on the product line bundle. It also uses notation that is used in Equation (2.14) from Section 3d.

**Proposition 10.9**: *Let $\{r_n, (A_n, a_n)\}_{n \in \Lambda}$ denote the sequence from Proposition 10.1. Fix a point $p \in X-Z$. There is a ball $B \subset X-Z$ centered at $p$ and, given a subsequence $\Lambda_p \subset \Lambda$, there exists the following data on $B$:*

- *A smooth section section of $\mathbb{S}^+$ over $B$ to be denoted by $\nu_B$ that obeys the following:*
  i) $|\nu_B| = |\nu|$.
  ii) $\langle \nu_B, \tau \nu_B \rangle = 0$.
  iii) $D_{\theta_0 + \hat{A}_B} \nu_B = 0$ *with $\hat{A}_B$ denoting the $i\mathbb{R}$-valued 1-form given by*

$$\hat{A}_B = -\tfrac{1}{2} \tfrac{1}{|\nu_B|^2} (\langle \nu_B, \nabla_{\theta_0} \nu_B \rangle - \langle \nabla_{\theta_0} \nu_B, \nu_B \rangle) .$$

- *A sequence, denoted by $\{g_n\}_{n \in \Lambda}$, of isomorphisms from product bundle $B \times \mathbb{C}$ to $E_B$.*
- *A subsequence $\Theta \subset \Lambda_p$.*

*These are such that the following is true: For any $n \in \Lambda$, write $g_n{}^*A_n$ (which is a connection on the product $\mathbb{C}$-bundle over $B$) as $\theta_0 + A_n$ with $A_n$ being an $i\mathbb{R}$ valued 1-form on $B$. The sequence $\{A_n\}_{n \in \Theta}$ converges to $\hat{A}_B$ in the $L^2_1$ weak topology; and the sequence $\{g_n{}^*a_n\}_{n \in \Theta}$ converges to $\nu_B$ in the $L^2_2$ weak topology to $\nu_B$.*

Note that the formula for $\hat{A}_B$ given in the proposition is equivalent to the assertion that

$$\langle \nu_B, \nabla_{A_B} \nu_B \rangle - \langle \nabla_{A_B} \nu_B, \nu_B \rangle = 0 .$$

(10.35)

This identity implies that the curvature 2-form of the connection $\hat{A}_B = \theta_0 + \hat{A}_B$ obeys

$$F_{\hat{A}_B} = -\tfrac{1}{|\nu_B|^2} (\langle \nabla_{A_B} \nu_B \wedge \nabla_{A_B} \nu_B \rangle + \langle \nu_B, \mathcal{R}^\nabla \nu_B \rangle) .$$

(10.36)

(The notation $\langle \cdot \wedge \cdot \rangle$ is defined in (1.14).)

The rest of this subsection contains the proof of this proposition.

*Proof of Proposition 10.9*: The proof of the proposition has six parts.

*Part 1*: Having fixed $\mu \in (0, \tfrac{1}{100})$, let $\kappa_\ddagger$ denote the minimum of the versions of $\kappa_\mu$ that appear in Propositions 3.1-3.3 and 4.1 and in Lemmas 10.5 and 10.6. Fix $c > \kappa_\ddagger$ and $\delta \in (0, \kappa_\ddagger^{-1} \mu^3)$. Supposing that $n \in \Lambda$, define the functions $r_{\ddagger n}(\cdot)$ and $r_{cFn}(\cdot)$ using $\mu$, $c$ and $\delta$ using the instructions from Parts 1 and 3 of the previous subsection. As noted in Lemma 10.8, the corresponding sequences $\{r_{\ddagger n}(\cdot)\}_{n \in \Lambda}$ and $\{r_{cFn}(\cdot)\}_{n \in \Lambda}$ are uniformly bounded away form zero on compact subsets of $X-Z$.

Given $p \in X-Z$, fix a ball centered at $p$ so that its closure is in $X-Z$. Let $r_{\ddagger\ddagger}$ denote a positive lower bound for all of the functions $\{r_{\ddagger n}(\cdot)\}_{n \in \Lambda}$ on this ball; and let $r_{cF\ddagger}$ denote a positive lower bound for all of the functions $\{r_{cFn}(\cdot)\}_{n\in\Lambda}$ on this ball. As in Part 1 of the proof of Proposition 10.7, fix a number $r_c$ to be less than half the radius of this same ball, and less than half minimum of the numbers $c^{-1}$, $r_{\ddagger\ddagger}$ and $r_{cF\ddagger}$. Keep in mind that the sequence of curvatures $\{F_{A_n}\}_{n\in\Lambda}$ has $L^2$ norm bounded by $c^{-1}$ on the radius $2r_c$ ball centered at p. (Also, (10.33) holds for $s \leq c_\mu^{-1} r_c$.) One other upper bound constraint on $r_c$ is needed for what is done in the next subsection. To give this constraint, note that there exists a positive number (to be denoted by $\rho_*$) with the property that any ball in X of radius $\rho_*$ or less is geodesically convex. Make sure that $r_c$ is less than this number $\rho_*$. The ball B in Proposition 10.9 can be any ball with radius $(1-2\mu)r_c$ or less centered at p.

*Part 2*: Having fixed $n \in \Lambda$, there is an isomorphism from the product bundle over $B_{r_c}$ (the radius $r_c$ ball centered at p) to the bundle E that pulls back the connection $A_n$ to give a connection on the product $\mathbb{C}$ bundle that has the form $\theta_0 + A_n$ with $\theta_0$ denoting the product connection and with $A_n$ denoting an $i\mathbb{R}$ valued 1-form that obeys

- $d*A_n = 0$ .
- *The pull-back of $*A_n$ to $\partial B_{r_c}$ is zero*.
- *The $L^2_1$ norm of $A_n$ on $B_{r_c}$ is no greater than $c_0$ times the $L^2$ norm of $F_{A_n}$ on $B_{r_c}$*.

(10.37)

The proof of (10.37) amounts to solving the Laplace equation on $B_{r_c}$ with suitable boundary values. See, e.g. [U].

Since the $L^2$ norm on $B_{r_c}$ of any $n \in \Lambda$ version of $F_{A_n}$ is at most $c^{-1}$, it follows from the third bullet of (10.37) that the sequence $\{A_n\}_{n\in\Lambda}$ is uniformly bounded in the $L^2_1$ topology on the space of $i\mathbb{R}$ valued 1-forms defined on $B_{r_c}$. As a consequence any subsequence in $\Lambda$ has, itself, a subsequence (to be denoted here by $\vartheta$) such that the corresponding subsequence of $\{A_n\}_{n\in\vartheta}$ converges weakly in the $L^2_1$ topology on the space of $i\mathbb{R}$ valued 1-forms on $B_{r_c}$. Let $\hat{A}_B$ denote the weak limit.

*Part 3*: For each $n \in \vartheta$, let $\hat{a}_n$ denote the pull-back of $a_n$ by the isomorphism in Part 2 from the product $\mathbb{C}$ bundle over $B_{r_c}$ that writes $A_n$ as $\theta_0 + A_n$. (This is a section of $\mathbb{S}^+$ over $B_{r_c}$.) Since the sequence $\{A_n\}_{n\in\vartheta}$ is bounded in the $L^2_1$ topology, it follows from the second bullet of Proposition 3.2 with the Sobolev inequalities in (2.1) that the sequence $\{\hat{a}_n\}_{n\in\Theta}$ is bounded in the $L^2_2$ topology on the space of sections of $\mathbb{S}^+$ over the ball of radius $(1-\mu)r_c$ centered at p. Therefore, there exists a subsequence in $\vartheta$, this denoted by $\Theta$, such that $\{\hat{a}_n\}_{n\in\Theta}$ converges weakly in the $L^2_2$ topology on the space of sections of $\mathbb{S}^+$ over $B_{(1-\mu)r_c}$. Denote this limit by $\nu_B$. The claim that $|\nu_B| = |\nu|$ follows from the fact that $|\hat{a}_n|$ converges strongly in the $L^2$ topology.

*Part 4*: This part considers the limits $\hat{A}_B$ and $\nu_B$. The first point is that $\nu_B$ obeys

$$\langle \nu_B, \tau \nu_B \rangle = 0 \tag{10.38}$$

because of the top bullet in (1.5), and $\nu_B$ (with $\hat{A}_B$) obey the covariant Dirac equation

$$D_{\theta_0 + \hat{A}_B} \nu_B = 0 . \tag{10.39}$$

because of second bullet in (1.5).

There is, in addition, an algebraic equation to relating $\hat{A}_n$ to the covariant derivative of $\nu_B$. To derive this equation, fix any smooth, $i\mathbb{R}$ valued 1-form with compact support on $B_{(1-\mu)r_c}$. Denote this 1-form by $s$. Supposing that $n \in \Theta$, take the inner product of the $r = r_n$ and $(A, a) = (A_n, a_n)$ version of both sides of (2.14) with $r_n^{-2} s$ and then integrate the resulting identity over $B_{(1-\mu)r_c}$. An integration by parts gives

$$- r_n^{-2} \int_{B_{r_c}} \langle ds, dA_n \rangle = \int_{B_{r_c}} \langle s, \langle \hat{a}_n, \nabla_{\theta_0 + A_n} \hat{a}_n \rangle - \langle \nabla_{\theta_0 + A_n} \hat{a}_n, \hat{a}_n \rangle \rangle + r_n^{-2} \int_{B_{r_c}} \langle s, *d * \omega \rangle . \tag{10.40}$$

Take $n \in \Theta$ ever large in (10.40) and use the weak convergence from Parts 2 and 3 to see that the $\hat{A}_B$ and $\nu_B$ obey (10.35). Since $|\nu_B|$ is nowhere zero (it is equal to $|\nu|$), the equation in (10.35) can be used to write $\hat{A}_B$ as

$$\hat{A}_B = -\tfrac{1}{2} \frac{1}{|\nu_B|^2} (\langle \nu_B, \nabla_{\theta_0} \nu_B \rangle - \langle \nabla_{\theta_0} \nu_B, \nu_B \rangle) \tag{10.41}$$

This depiction of $\hat{A}_B$ is used in Part 5 to prove that $\nu_B$ and $A_B$ are smooth.

*Part 5*: The Dirac equation in (10.39) leads to a version of (2.8) that has the schematic form

$$\nabla_{\theta_0}^{\dagger} \nabla_{\theta_0} \nu_B + \frac{1}{|\nu_B|^4} \mathcal{E}(\nabla_{\theta_0} \nu_B, \nu_B) \nu_B + \mathfrak{R} \nu_B = 0 \tag{10.42}$$

with $\mathcal{E}(\cdot, \cdot)$ being a smooth fiber preserving map from $(\mathbb{S}^+ \otimes T^*X) \otimes \mathbb{S}^+$ to $\text{End}(\mathbb{S}^+)$ over B that is quadratic in its first entry. (This is to say that $\mathcal{E}(\flat, \cdot)$ is a quadratic function of the components of $\mathfrak{b}$.)  Note that each term in this equation is an $L^2$ section of $\mathbb{S}^+$ by virtue of (2.1) and the fact that $\nu_B$ is nowhere zero, bounded, $L^2_2$ section of $\mathbb{S}^+$. Elliptic regularity techniques such as those in Chapter 6 of [Mor] can be brought to bear using (10.42) to prove that $\nu_B$ is $C^\infty$ on B.  This implies in turn that $A_B$ is also smooth on B.

By way of an explanatory remark, this elliptic regularity buisness can be put into a perhaps more familiar context as follows:  Write $\mathfrak{b} = \nabla_{\theta_0} \nu_B$. The equation in (10.42) and the identity $d_{\theta_0} \nabla_{\theta_0} \nu_B = \mathcal{R} \nu_B$ (with $\mathcal{R}$ denoting an endomorphism that is defined from the data giving the Clifford module) can be written in the schematic form

$$\mathcal{L}\mathfrak{b} + \mathfrak{E}(\mathfrak{b} \otimes \mathfrak{b}) = \mathfrak{g}$$

(10.43)

with $\mathcal{L}$ being an elliptic, first order operator with smooth coefficients, with $\mathfrak{g}$ in the $L^2_2$ Sobolev space, and with $\mathfrak{E}$ being a homomorphism with domain $\otimes^2(\mathbb{S}^+ \otimes T^*X)$ having $L^2_2$ coefficients.  In dimension 4, this sort of schematic equation with $\mathfrak{b}$ in the $L^2_1$ Sobolev space is analogous to the equations for a self-dual connection on a 4-manifold (in the Hodge gauge) with $\mathfrak{b}$ playing the role of the connection.  The fact that the smoothness of solutions to the latter equation can be proved using the standard bootstrapping arguments is well known to gauge theory afficianados (see Corollary 1.4 in [U]).

**f) Global convergence**

Proposition 10.9 describes the local convergence on X–Z of a sequence from Proposition 10.1.  The next proposition puts Proposition 10.9 into a larger context.

**Proposition 10.10**: *Let $\{r_n, (A_n, a_n)\}_{n \in \Lambda}$ denote the sequence from Proposition 10.1. There exists the following data:*
- *A smooth section section of $\mathbb{S}^+ \otimes E$ over X–Z to be denoted by $\nu$ and a smooth connection on $E|_{X-Z}$ to be denoted by $\hat{A}$ that obey the following:*
  i)   *The norm of $\nu$ is the function $|\nu|$ from Proposition 10.1.*
  ii)  *The function $|\nabla_{\hat{A}} \nu|$ extends to X as an $L^2$ function.*
  iii) *$\langle \nu, \tau \nu \rangle = 0$.*
  iv)  *$D_{\hat{A}} \nu = 0$.*
  v)   *$\langle \nu, \nabla_{\hat{A}} \nu \rangle - \langle \nabla_{\hat{A}} \nu, \nu \rangle = 0$  and  $F_{\hat{A}} = -\frac{1}{|\nu|^2}(\langle \nabla_{\hat{A}} \nu \wedge \nabla_{\hat{A}} \nu \rangle + \langle \nu, \mathcal{R}^\nabla \nu \rangle)$.*
- *A sequence, denoted by $\{g_n\}_{n \in \Lambda}$, of automorphisms of E over X–Z.*
- *A subsequence $\Theta \subset \Lambda$.*

*These are such that the following is true:  The sequence $\{g_n^* A_n\}_{n \in \Theta}$ converges to $\hat{A}$ in the $L^2_1$ weak topology on compact subsets of X–Z; and the sequence $\{g_n^* a_n\}_{n \in \Theta}$ converges to* $v$ in the $L^2_2$ weak topology on compact subsets of X–Z. Moreover, $\{|a_n|\}_{n\in\Theta}$ converges to $|v|$ in the $L^2_1$ and the $C^0$ topologies on the whole of X and $\{|\nabla_{A_n} a_n|\}_{n\in\Theta}$ converges to $|\nabla_{\hat{A}} v|$ in the $L^2$ topology on the whole of X.

**Proof of Proposition 10.10**: The convergence of $\{|a_n|\}_{n\in\Theta}$ to $|v|$ on the whole of X in the $C^0$ topology has already been verified. This fact is used in what follows The proof of the rest of the proposition has six parts.

*Part 1*: Given $p \in$ X–Z, let $B_p \subset$ X–Z denote the point p's version of Proposition 10.9's ball B. Let $B_{\mu,p}$ denote the ball with center p whose radius is smaller by the factor $1-\mu$. Since the set $\{B_{\mu,p}\}_{p\in X-Z}$ is an open cover of X–Z, it has a locally finite, countable subcover. Fix such a subcover and let $\mathcal{V}$ denote the set whose elements are the corresponding versions of $B_p$ (thus, if $B_p \in \mathcal{V}$, then $B_{\mu,p}$ is in the subcover). Label the balls in $\mathcal{V}$ by consective, positive integers starting from 1.

Invoke Proposition 10.9 with p being the center point of $B_1$ and with $\Lambda_p$ being the whole sequence $\Lambda$. Let $\Theta_1$ denote the corresponding subsequence of $\Lambda$ that is supplied by the third bullet of Proposition 10.9. Having started with $B_1$, then invoke Proposition 10.9 with p being the center point of the ball $B_2$ and with $\Lambda_p$ being the subsequence $\Theta_1$. Use $\Theta_2$ to denote the subsequence of $\Theta_1$ that is supplied by the third bullet of Proposition 10.9. Continue in this vein sequentially to generate nested subsequences $\Theta_1 \supset \Theta_2 \supset ...$ with $\Theta_k$ for any integer $k \geq 2$ being the subsequence of $\Theta_{k-1}$ that is obtained from the version of Proposition 10.9 with p being the center point of $B_k$ and with $\Lambda_p$ being $\Theta_{k-1}$. Define now $\Theta = \{n_k\}_{k=1,2,...}$ to be the (increasing) subsequence of $\Lambda$ that is defined so that $n_1 \in \Theta_1$, $n_2$ is a point from $\Theta_2$ that is greater than $n_1$, then $n_3$ is a point from $\Theta_3$ that is greater than $n_2$, and so on.

*Part 2*: Each $B \in \mathcal{V}$ has its corresponding ball and on this ball, the Proposition 10.9 data $(v_B, \hat{A}_B)$ consisting of a section of $\mathbb{S}^+$ and an $i\mathbb{R}$ valued 1-form. The first task is to prove that there is a complex line bundle $E' \to$ X–Z such that $\{(\hat{A}_B, v_B)\}_{B\in\mathcal{V}}$ is the cocycle data for a pair of section, $v$, over X–Z of $\mathbb{S}^+ \otimes E'$ and connection, $\hat{A}$, on $E'$ over X–Z. This will be the case if the following is true: Let B and B´ denote any two intersecting pairs from $\mathcal{V}$. There is a smooth map $u: B \cap B' \to S^1$ such that $(\hat{A}_B, v_B)$ and $(\hat{A}_{B'}, v_{B'})$ on $B \cap B'$ are related by the rule whereby

$$\hat{A}_B = \hat{A}_{B'} - u^{-1}du \quad and \quad v_B = u v_{B'}$$

(10.44)

(It follows from Item iii) of the third bullet in Proposition 10.9 that the right most identity implies the left most identity.) Note in this regard that if (10.44) holds for each pair of intersecting balls from $\mathcal{V}$ then the cocycle condition on triple intersections (which is required for the various $u$'s to define a complex line bundle) is *automatically* obeyed. This is because the left most identity in (10.44) can be used write any given B, B´ version of $u$ as $\langle v_{B'}, v_B \rangle$.

To summarize, if (10.44) holds for each pair of intersecting balls from $\mathcal{V}$, then there is a complex line bundle E´ → X–Z and the data $\{(\hat{A}_B, v_B)\}_{B \in \mathcal{V}}$ defines a pair of connection, $\hat{A}$, on E´ and section, $v$, on $\mathbb{S}^+ \otimes E´$. Granted this, then the conditions in Items i) and iii) and iv) of the top bullet of Proposition 10.10 follow directly from the conditions in Items i) and ii) and iii) of Proposition 10.9. Item ii) of the top bullet of Proposition 10.10 is discussed in the next paragraph. Item v) of Proposition 10.10's top bullet follows from (10.35) and (10.36).

The function $|\nabla_{\hat{A}} v|$ is defined a priori on X–Z because $\hat{A}$ and $v$ are smooth. This function is extended to the whole of X by declaring it to be zero on Z. This extension is a square integrable function on X because of the convergence assertion in Propositions 10.8 (remember that $|\nabla_{\hat{A}} v|$ is $|\nabla_{\hat{A}_B} v_B|$ on any $B \in \mathcal{V}$ because of (10.44)) and because of Item ii) of Proposition 10.1's fifth bullet..

*Part 3*: To obtain $u$ in (10.44), fix $n \in \Theta$ for the moment and define $\hat{a}_n$ on B from $a_n$ as done in Part 3 of the proof of Proposition 10.9. Define $\hat{a}_n´$ by replacing the ball B with B´ in this same part of Proposition 10.9's proof. Define $u_n$ to by the rule

$$u_n = \frac{1}{|a_n|^2} \langle \hat{a}_n´, \hat{a}_n \rangle.$$

(10.45)

This is a map from B∩B´ to $S^1$ with the property that

$$\hat{a}_n = u_n \hat{a}_n´ \quad and \quad A_n = A_n´ - u_n^{-1} du_n.$$

(10.46)

The sequence $\{u_n\}_{n \in \Theta}$ is bounded in the $L^2_2$ topology on the space of maps from B∩B´ to $S^1$ by virtue of the fact that $|a_n|$ is bounded away from zero on B and B´, and by virtue of the fact that sequences $\{\hat{a}_n\}_{n \in \Theta}$ and $\{\hat{a}_{n'}\}_{n' \in \Theta}$ have uniformly bounded $L^2_2$ norms. It is a consequence of this convergence (and the $C^0$ convergence of $\{|a_n|\}_{n \in \Theta}$ to $|v|$) that the sequence $\{u_n\}_{n \in \Theta}$ converges weakly in the $L^2_2$ topology to the map

$$u = \frac{1}{|v|^2} \langle v_{B'}, v_B \rangle.$$

(10.47)

The latter map to $S^1$ is necessarily smooth (because $v_B$ and $v_{B'}$ are smooth and $|v|$ is not zero on $B \cap B'$). The right most equation in (10.44) is automatically obeyed with this version of $u$, and, as noted previously, then the right most equation must also be obeyed.

*Part 4*: This part of the proof explains why the bundle E´ is isomorphic to the bundle E. This is the case in particular if the following is true: If balls B and B´ from $\mathcal{V}$ intersect, then the corresponding sequence $\{u_n\}_{n\in\Theta}$ from (10.45) converges to $u$ in the $C^\infty$ topology on compact subsets of $B \cap B'$. (See, e.g. [U] or the appendix to [Ta].) Moreover, arguments in that differ little from what is said in these two references can also be used to construct the bundle isomorphisms $\{g_n\}_{n\in\Theta}$ for the second bullet of Proposition 10.10 so that $\{(g_n{}^*A_n, g_n{}^*a_n)\}_{n\in\Theta}$ converge to $(\hat{A}, v)$ as required by Proposition 10.10. (Each $n \in \Theta$ version of the map $g_n$ is constructed from the various B, B´ $\in \mathcal{V}$ versions of $\{u^{-1}u_n\}$ using cut-off functions and the fact that any given $u^{-1}u_n$ differs little from the constant function 1 when n is large.) The proof of the $C^\infty$ convergence property starts in the next paragraph.

Since the balls B and B´ are geodesically convex (see Part 1 of the proof of Proposition 10.9), their intersection is also geodesically convex. This implies in particular that the $S^1$ valued function $u^{-1}u_n$ has a logarithm. This is to say that there is a smooth function, $\mathbb{R}$ valued on $B \cap B'$ to be denoted by $w_n$ such that $u^{-1}u_n = e^{iw_n}$. This function $w_n$ is uniquely defined modulo the addition of constant multiples of $2\pi$; and it proves useful to choose the version with the property that its average,

$$\underline{w}_n = \frac{1}{\text{vol}(B\cap B')} \int_{B\cap B'} w_n ,$$

(10.48)

obeys $|\underline{w}_n| \in [-\pi, \pi]$.

There are three key features of $w_n$ to keep in mind: The first is that $\nabla w_n$ converges weakly to zero in the $L^2_1$ topology. This is because

$$\nabla w_n = i(A_n - \hat{A}_B) - i(A_n' - \hat{A}_{B'}) .$$

(10.49)

The second point (which is ultimately the reason for the constraint on $\underline{w}_n$) is this: The balls B and B´ determine a positive number (denoted by $\lambda$) such that

$$\int_{B\cap B'} |\nabla w_n|^2 \geq \lambda \int_{B\cap B'} |w_n - \underline{w}_n|^2 .$$

(10.50)

I will talk about the third point momentarily. These first two points imply that the sequence $\{w_n - \underline{w}_n\}_{n \in \Theta}$ converges weakly to zero in the $L^2_2$ topology and thus it converges strongly to zero in the $L^2_1$ topology.

The third key point is that $w_n$ is a *harmonic* function on $B \cap B'$. This is because $A_n$ and $A_n'$ and $v_B$ and $v_{B'}$ are all coclosed 1-forms. Indeed, the top bullet of (10.37) says this about $A_n$ and its primed counterpart says it about $A_n'$; and the convergence of $\{\hat{A}_n\}_{n \in \Theta}$ to $\hat{A}_B$ in the weak $L^2_1$ topology implies that $\hat{A}_B$ is coclosed (and similarly for $\hat{A}_{B'}$). The coclosed identity $d*\hat{A}_B = 0$ can also be derived directly from the equations in Item iii) of the first bullet of Proposition 10.9. Since $w_n$ is harmonic on $B \cap B'$, so is $w_n - \underline{w}_n$.

*Part 5*: By way of a summary, the sequence $\{w_n - \underline{w}_n\}_{n \in \Theta}$ is a sequence of harmonic functions that converges strongly to zero in the $L^2_1$ topology for $C^\infty(B \cap B')$. Because each of these functions is harmonic, the $L^2$ bounds imply $C^\infty$ bounds on compact subsets. Therefore the sequence converges strongly to zero in the $C^\infty$ topology on any open set in $B \cap B'$ with compact closure in $B \cap B'$. It follows from this that the sequence $\{u^{-1}u_n\}_{n \in \Theta}$ converges to the constant function 1 in the $C^\infty$ topology on open sets in $B \cap B'$ with compact closure if the sequence of numbers $\{\underline{w}_n\}_{n \in \Theta}$ converges to zero.

To see about the the convergence to zero of $\{\underline{w}_n\}_{n \in \Theta}$, note first that this is a bounded sequence of numbers because all of its terms are between $-\pi$ and $\pi$. There are, therefore, convergent subsequences. Fix one and denoted it by $\{\underline{w}_n\}_{n \in \Xi}$ and let $\underline{w}$ denote its limit. The corresponding sequence $\{w_n\}_{n \in \Xi}$ converges in the $C^\infty$ topology on compact subsets of $B \cap B'$ to $\underline{w}$. (This is a priori a number from the interval $[-\pi, \pi]$.) As a consequence, the sequence $\{u^{-1}u_n\}_{n \in \Xi}$ converges in the $C^\infty$ topology on compact subsets of $B \cap B'$ to $e^{i\underline{w}}$. What with (10.45) and (10.47), this means that $\{\langle \hat{a}_n', \hat{a}_n \rangle\}_{n \in \Xi}$ converges in the $C^\infty$ topology on compact subsets of $B \cap B'$ to $e^{i\underline{w}} \langle v_{B'}, v_B \rangle$. But such an event requires that $\underline{w}$ be zero because $\{\langle \hat{a}_n', \hat{a}_n \rangle\}_{n \in \Theta}$ converges in the $L^2$ topology on $B \cap B'$ to $\langle v_{B'}, v_B \rangle$.

*Part 6*: This last part of the proof explains why the sequence $\{|a_n|\}_{n \in \Theta}$ converges to $|v|$ in the $L^2_1$ topology on the whole of X; and why the sequence $\{|\nabla_{A_n} a_n|\}_{n \in \Theta}$ converges to $|\nabla_{\hat{A}} v|$ in the $L^2$ topology on the whole of X. To this end, consider first the asserted $L^2$ convergence of $\{|\nabla_{A_n} a_n|\}_{n \in \Theta}$. Fix for the moment $\delta \in (0, 1)$ and let $\mathcal{U}_\delta$ denote the subset of X where $|v| < \frac{1}{2}\delta$. If $n \in \Theta$ is sufficiently large, then $|a_n| < \delta$ on $\mathcal{U}_\delta$. This understood, then what is said in Parts 1-3 of the proof of Proposition 5.3 using the pair $(A = A_n, a = a_n)$ imply that

$$\int_{\mathcal{U}_\delta} |\nabla_{A_n} a_n|^2 \leq c_0\, \delta^{1/c_0} .$$

(10.51)

This bound and Item ii) of the second bullet of Proposition 10.1 imply in turn that

$$\int_{\mathcal{U}_\delta} |\nabla_{\hat{A}} v|^2 \leq c_0\, \delta^{1/c_0} .$$

(10.52)

To prove the $L^2$ convergence assertion for $\{|\nabla_{A_n} a_n|\}_{n \in \Theta}$, fix $\varepsilon \in (0, 1)$ and take $\delta$ in (10.51) and (10.52) so that the right hand side of both equations are less than $\frac{1}{16}\varepsilon^2$ (this will be the case if $\delta \leq c_0^{-1}\varepsilon^{c_0}$.) Given $\varepsilon$ and $\delta$, it follows from what was said in Parts 1-5 that there exists a positive integer N such that the $L^2$ norm of $|\nabla_{A_n} a_n| - |\nabla_{\hat{A}} v|$ on $X-\mathcal{U}_\delta$ is less than $\frac{1}{4}\varepsilon$ when $n \in \Theta$ is greater than N. This fact and the fact that the right hand sides of (10.51) and (10.52) are less than $\frac{1}{16}\varepsilon^2$ imply that the $L^2$ norm of $|\nabla_{A_n} a_n| - |\nabla_{\hat{A}} v|$ on the whole of X is less than $\varepsilon$ when $n > N$. This is what is needed to prove that $\{|\nabla_{A_n} a_n|\}_{n \in \Theta}$ converges to $|\nabla_{\hat{A}} v|$ on X in the $L^2$ topology. Almost the same argument proves the $L^2_1$ convergence of $\{|a_n|\}_{n \in \Theta}$ to $|v|$ on X in the $L^2_1$ topology because the respective integrals in (10.51) and (10.52) are no greater than those of $|\nabla |a_n||^2$ and $|\nabla |v||^2$. (The $L^2$ convergence of $\{|a_n|\}_{n \in \Theta}$ to $|v|$ on X follows from the $C^0$ convergence of $\{|a_n|\}_{n \in \Theta}$ to $|v|$ on X.)

**g) The question of density**

This subsection proves that the set Z is nowhere dense. The formal assertion that this is so is an instance of the upcoming Proposition 10.11. To set the stage for this proposition, let Z denote a closed set in X, let $E \to X-Z$ denote a complex, Hermitian vector bundle. Let $\hat{A}$ denote a smooth connection on $E|_{X-Z}$ and let $v$ denote a section of $\mathbb{S}^+ \otimes E$ over $X-Z$. Make the following additional assumptions about $v$ and $\hat{A}$:

- *The norm of $|v|$ extends to the whole of X as a continuous and also $L^2_1$ function that vanishes on Z.*
- *$|\nabla_{\hat{A}} v|$ extends over Z to define an $L^2$ function on X.*
- *$D_{\hat{A}} v = 0$ .*
- *If $\omega$ is any given 2-form, then $\langle v, \langle \omega, F_{\hat{A}} \rangle \mathrm{cl}(\omega) v \rangle = 0$ .*

(10.53)

With regards to the fourth bullet, note that $\langle \omega, F_{\hat{A}} \rangle$ is an endomorphism of E. If $\hat{A}$ is flat or anti-self dual, then the fourth bullet is automatically satisfied. Likewise, if $\langle v, \tau v \rangle = 0$, then the fourth bullet is automatically satisfied. In particular, the data supplied by Proposition 10.10 is an instance of $(\hat{A}, v)$ that obeys (10.53).

***Proposition 10.11***: *If the closed set* $Z \subset X$ *and the pair* $(\hat{A},v))$ *is as described above (and, in particular, obey (10.53)), then Z is nowhere dense.*

***Proof of Proposition 10.11***: Fix for the moment a point $p \in X$ and define the function $\mathcal{K}$ on $(0, \infty)$ by the rule in (10.2). Since Z is closed, the point $p \in X$ is a point of density for Z only in the event that $|v|$ is identically zero in some small radius ball centered at p. In turn, this will occur if and only if the function $\mathcal{K}$ is zero on some interval $(0, r_0] \subset (0, c_0^{-1}]$ The proof that $\mathcal{K} > 0$ on $(0, c_0^{-1}]$ is had by mimicking almost word for word the arguments in Sections 2 of [T2] and those in Section 3a in [T2] up through the statement of Lemma 3.2 and those in Section 3b of [T2]. Note that these sections of [T2] prove the proposition when E is a line bundle and $\hat{A}$ is a flat connection on E with $\mathbb{Z}/2$ holonomy. The only change to what is said in these parts of [T2] is a justification for the analog that follows of Equation (2.1) in [T2]:

$$\tfrac{1}{2} d^\dagger d |v|^2 + |\nabla_{\hat{A}} v|^2 + \langle v, \mathfrak{R} v \rangle = 0.$$

(10.54)

The proof of (10.54) starts with the $(A = \hat{A}, a = v)$ version of (2.8). Take the inner product of both sides with $v$. Since the left hand side of the resulting identity is zero (by virtue of the third bullet of (10.53)), the identity gives an equation that is identical to (10.54) except for term $\tfrac{1}{2} \langle v, \text{cl}(F_{\hat{A}}) v \rangle$. But, this term is actually zero because of the fourth bullet in (10.53).

## 11. Proof of Proposition 1.2

The proof of Proposition 1.2 has four parts.

*Part 1*: To prove the assertion in the first bullet of the proposition, let $\{\omega^a\}_{a=1,2,3}$ denote an oriented, orthonormal basis for $\Lambda^+$ over a given ball in X−Z; and then, for each index $a \in \{1, 2, 3\}$, let $\tau^a = \tfrac{1}{\sqrt{2}} \text{cl}(\omega^a)$. The condition $\text{cl}^\dagger(v \otimes v^\dagger) = 0$ is equivalent to the assertion that the vectors $\{v, \tau^1 v, \tau^2 v, \tau^3 v\}$ define an orthogonal basis for $\mathbb{S}^+$ (which is orthonormal if $|v| = 1$) over the given ball. Use this basis to write $\nabla_{\hat{A}} v$ as

$$\nabla_{\hat{A}} v = \alpha^0 v + \sum_{c=1,2,3} \alpha^c \tau^c v$$

(11.1)

with $\{\alpha^0, \alpha^1, \alpha^2, \alpha^3\}$ being $\mathbb{C}$-valued 1-forms on the given ball. Writing $\nabla_{\hat{A}} v$ in this way leads to the identities

- $\langle v, \nabla_{\hat{A}} v \rangle = \alpha^0 |v|^2$.
- $\langle v, \tau^c \nabla_{\hat{A}} v \rangle = -\alpha^c |v|^2$.

(11.2)

Now, $\langle v, \nabla_{\hat{A}} v \rangle$ is real (see Item c) of the second bullet of Theorem 1.1) so $\alpha^0$ is real. Likewise, each $\alpha^c$ is real. This can be seen by differentiating the identity $\langle v, \tau v \rangle = 0$ (it is the identity $cl^\dagger(v \otimes v^\dagger) = 0$ from Item c) of the second bullet of Theorem 1.1.) Since the 1-forms $\{\alpha^0, \alpha^1, \alpha^2, \alpha^3\}$ are real, and since $\langle v, \tau v \rangle = 0$, it follows from (11.1) that

$$\langle \nabla_{\hat{A}} v \wedge \nabla_{\hat{A}} v \rangle = 0 .$$

(11.3)

This last observation and Item d) of the second bullet of Theorem 1.1 lead directly to the assertion of the first bullet in Proposition 1.2.

*Part 2*: Assume in this part and henceforth that $\mathbb{S}^+$ has a complex conjugation (which is denoted by $C$). The Clifford module over X–Z that is defined by the data $\{\mathbb{S}^+ \otimes E, \mathbb{S}^- \otimes E, \nabla_{\hat{A}}, \mathfrak{c}\}$ inherits a complex conjugation from $(\mathbb{S}^+, \mathbb{S}^-, \nabla, \mathfrak{c})$, which is also denoted by $C$; it maps $\mathbb{S}^+ \otimes E$ to $\mathbb{S}^+ \otimes \mathcal{L} \otimes E^{-1}$. (This occurs without constraint on $\hat{A}$ but for it being a Hermitian connection on E.)

Fix $p \in X$ and suppose that $v \in (\mathbb{S}^+ \otimes E)|_p$ and that $Cv = uv$ with $u \in (\mathcal{L} \otimes E^{-2})|_p$. Thus, $v$ is an eigenvector of sorts for $C$. Since $C$ preserves norms (the third bullet of (1.8)), the element $u$ has norm 1. Here is a second observation: If $Cv = uv$ and if $\varsigma$ is a self-dual 2-form, then

$$\langle v, cl(\varsigma) v \rangle = 0 .$$

(11.4)

This is because $\langle v, cl(\varsigma) v \rangle = \langle C(cl(\varsigma) v), Cv \rangle$; which is the same as $\langle cl(\varsigma) C(v), Cv \rangle$ because $C$ commutes with $cl(\cdot)$; which is the same as $\langle cl(\varsigma) v, v \rangle$ because $|u| = 1$. Meanwhile, $\langle cl(\varsigma) v, v \rangle = -\langle v, cl(\varsigma) v \rangle$ because $cl(\cdot)$ is anti-Hermitian.

It follows from (11.4) that the $i\mathbb{R}$ valued, self-dual 2-form $cl^\dagger(v \otimes v^\dagger)$ vanishes identically on any open set where there is a section $u$ of $\mathcal{L} \otimes E^{-2}$ such that $Cv = uv$. (Note that this is true when $\dim_{\mathbb{C}}(\mathbb{S}^+) > 4$ also; the argument didn't refer at all to $\dim_{\mathbb{C}}(\mathbb{S}^+)$.)

There is a converse to the preceding observation if $\mathbb{S}^+$ has dimension 4 (over $\mathbb{C}$) which is as follows: Fix $p \in X$.

*If $v \in (\mathbb{S}^+ \otimes E)|_p$ and if $\mathrm{cl}^\dagger(v \otimes v^\dagger) = 0$, then $Cv = uv$ for some $u \in (\mathcal{L} \otimes E^{-2})|_p$.*

(11.5)

To prove (11.5), let $\{\omega^a\}_{a=1,2,3}$ again denote an oriented, orthonormal basis for $\Lambda^+$ at p; and, again, let $\tau^a = \frac{1}{\sqrt{2}} \mathrm{cl}(\omega^a)$. As noted previously, the condition $\mathrm{cl}^\dagger(v \otimes v^\dagger) = 0$ is equivalent to the assertion that the vectors $\{v, \tau^1 v, \tau^2 v, \tau^3 v\}$ define an orthogonal basis for $\mathbb{S}^+$ (which is orthonormal if $|v| = 1$). Therefore, $Cv$ can be written as $uv + \sum_{a=1,2,3} u^a \tau^a v$ on the set where the frame $\{\omega^a\}_{a=1,2,3}$ is defined. The task is to prove that $u^1$, $u^2$ and $u^3$ are zero. To do this, fix $c \in \{1, 2, 3\}$ and note, on the one hand, that $\langle v, \tau^c C(v) \rangle = -u^c$. On the other hand, $\langle v, \tau^c C(v) \rangle$ is also equal to $\langle C(\tau^c Cv), C(v) \rangle$ by virtue of the third bullet of (1.8). And, the latter expression is $-\langle v, \tau^c C(v) \rangle$ because of the second bullet in (1.8) and because $\tau^c$ is anti-Hermitian. Thus, $\langle v, \tau^c C(v) \rangle = 0$ and so $u^c = 0$. Since this is so for all $c \in \{1,2,3\}$, so $Cv = uv$ as claimed.

*Part 3*: Suppose that $\hat{A}$ is a smooth connection on E over a given open set in X and that $\mathcal{A}$ is likewise a smooth connection on $\mathcal{L}$. They induce, together, a Hermitian connection on $\mathcal{L} \otimes E^{-2}$. This induced connection is denoted by $\hat{A}^{(-2)} + \mathcal{A}$. (The open set will be X–Z, but denote it by U for the moment.)

Now, suppose that $v$ is a section of $\mathbb{S}^+ \otimes E$ over a given open set in X that obeys the equation $D_{\hat{A}} v = 0$. Then $D_{\hat{A}^{(-2)} + \mathcal{A}}(Cv) = 0$ also because of the first and second bullets in (1.8) (Use the identity in (2.26) which is equivalent to the assertion that $D_{\hat{A}} v = 0$.) If it is also the case $Cv = uv$ with $u$ being a section of $\mathcal{L} \otimes E^{-2}$ over the given open set, then $D_{\hat{A}} v$ and $D_{\hat{A}^{(-2)} + \mathcal{A}}(Cv)$ can vanish simultaneously (where $v \neq 0$) only in the event that $\nabla_{\hat{A}^{(-2)} + \mathcal{A}} u = 0$. This implies in particular that $\hat{A}^{(-2)} + \mathcal{A}$ is identified with the product connection when $u$ is used to identify $(\mathcal{L} \otimes E^{-2})|_U$ with $U \times \mathbb{C}$.

The homomorphism $u$ can be viewed equivalently as an isometric isomorphism from $\mathcal{L}$ to $E^2$ (to be denoted by $\hat{u}$); and the preceding observation about $\hat{A}^{(-2)} + \mathcal{A}$ being equivalent to the product connection implies that $\hat{u}$ identifies $\mathcal{A}$ with $\hat{A}^{(2)}$. It follows directly from this that E is a square root of $\mathcal{L}$ (over U) and that $F_{\hat{A}} = \frac{1}{2} F_{\mathcal{A}}$.

The observations in the preceding paragraph when applied to the connection $\hat{A}$ from Proposition 1.2 with U being X–Z are what is asserted the second bullet of Proposition 11.2.

*Part 4*: The third bullet of Proposition 11.2 will be seen momentarily to follow from what is said in [T2] about $\mathbb{Z}/2$ harmonic spinors. To see this, fix an open set $U \subset X$ where the bundles $\mathcal{L}$ and E have isomorphisms with the product bundle $U \times \mathbb{C}$. Having chosen the isomorphisms, use them to identify the respective connections $\mathcal{A}$ and $\hat{\mathcal{A}}$ with connections $\theta_0 + \hat{a}_\mathcal{A}$ and $\theta_0 + \hat{a}_{\hat{\mathcal{A}}}$ with $\theta_0$ denoting the product connection on $U \times \mathbb{C}$. (Note that $\hat{a}_\mathcal{A}$ is an $i\mathbb{R}$ valued 1-form on U whereas $\hat{a}_{\hat{\mathcal{A}}}$ is an $i\mathbb{R}$ valued 1-form that is defined a priori only on U−Z). Use the same isomorphisms to identify the homomorphism $\hat{u}$ on U−Z with a map (to be denoted by $\hat{u}$ also) from U−Z to $S^1$. It follows from what was said in Part 3 that

$$2\hat{a}_{\hat{\mathcal{A}}} + \hat{u}^{-1}d\hat{u} = \hat{a}_\mathcal{A} \quad on \quad U-Z$$

(11.6)

With the preceding on hold for a moment, note that the map $\hat{u}$ has two square roots on any given ball in U−Z and these square roots define a real line bundle $\mathcal{I} \to$ U−Z. Supposing that $B \subset$ U−Z is a given ball, let $\hat{o}_B$ denote a square root of $\hat{u}$ over B; and let $\hat{v}_B = \hat{o}_B^{-1} v$. The versions of $\{\hat{v}_B\}$ on the various balls in U−Z define a section of $(\mathbb{S}^+ \otimes \mathcal{I})|_{U-Z}$ which is denoted below by $\hat{v}$.

Meanwhile, the connection $\theta_0 + \frac{1}{2}\hat{a}_\mathcal{A}$ on the product bundle and the given covariant derivatives on $\mathbb{S}^+$ and $\mathbb{S}^-$ (these were denoted by $\nabla$) define a second covariant derivative on $\mathbb{S}^+|_{U-Z}$ and on $\mathbb{S}^-|_{U-Z}$ (this is denoted by $\nabla_{\theta_0 + \frac{1}{2}\hat{a}_\mathcal{A}}$) such that $\{\mathbb{S}^+|_{U-Z}, \mathbb{S}^+|_{U-Z}, \nabla_{\theta_0 + \frac{1}{2}\hat{a}_\mathcal{A}}, \mathcal{C}\}$ defines a Clifford module over U. Granted this, it follows from (11.6) and Item c) of the second bullet of Theorem 1.1 that $\hat{v}$ is annihilated by the Dirac operator $D_{\theta_0 + \frac{1}{2}\hat{a}_\mathcal{A}}$.

With the preceding understood, and since $|v|$ is Hölder continuous along its zero locus (which is asserted by Item a of Theorem 1.1's second bullet), it follows that the data set $\{Z \cap U, \mathcal{I}, \hat{v}\}$ on the given set U meets the requirements in [T2] that define a $\mathbb{Z}/2$ harmonic spinor. Therefore, the theorems in [T2] can be brought to bear to say more about Z and $v$; and Theorem 1.3 in [T2] in particular says that Z has Hausdorff dimension at most 2.